\documentclass[11pt, reqno]{amsart}
\usepackage{amscd,amsmath,amsopn,amssymb,amsthm,multicol}
\usepackage[color,matrix, all, 2cell]{xy}
\usepackage[new]{old-arrows}
\usepackage{bbold}
 \usepackage{graphicx}
\usepackage{slashed}
\usepackage{lscape}
 \usepackage{cancel}
\usepackage{mathrsfs}
\usepackage{euscript}
\usepackage{graphicx}
\usepackage{upgreek}
\usepackage{textgreek}
\usepackage{enumerate}
\usepackage{color}
\usepackage{tikz}
\usepackage{multirow}
\usepackage{cancel}
\usepackage{comment}
\usepackage{booktabs}
\usepackage{amsmath}

\newcommand{\bec}{\begin{center}}
\newcommand{\ecc}{\end{center}}

\DeclareMathOperator{\ad}{ad}

\DeclareMathOperator{\Ker}{\mathsf{Ker}}

\DeclareMathOperator{\Lie}{\mathsf{Lie}}

\DeclareMathOperator{\Id}{\mathsf{Id}}

\DeclareMathOperator{\Ric}{\mathsf{Ric}}

\DeclareMathOperator{\hol}{\mathsf{hol}}

\DeclareMathOperator{\Scal}{\mathsf{Scal}}
\DeclareMathOperator{\R}{\mathbb{R}}

\DeclareMathOperator{\Imm}{\mathsf{Im}}

 \DeclareMathOperator{\SU}{\mathsf{SU}}
 
\DeclareMathOperator{\U}{\mathsf{U}}

\DeclareMathOperator{\Gg}{\mathsf{G}}

\DeclareMathOperator{\Tg}{\mathsf{T}}

\DeclareMathOperator{\Ss}{\mathsf{S}}

\DeclareMathOperator{\Ed}{\mathsf{End}}

\definecolor{dblue}{rgb}{0.01,0.01,0.42}
\definecolor{red}{rgb}{0.57,0.11,0.15}
\newcommand{\dd}{\mathrm{d}}
\newcommand{\fr}{\mathfrak}

\newcommand{\bb}{\mathbb}
\newcommand{\mc}{\mathcal}

\newcommand{\ep}{\varepsilon}

\DeclareFontFamily{U}{mathx}{}
\DeclareFontShape{U}{mathx}{m}{n}{<-> mathx10}{}
\DeclareSymbolFont{mathx}{U}{mathx}{m}{n}
\DeclareMathAccent{\widehat}{0}{mathx}{"70}
\DeclareMathAccent{\widecheck}{0}{mathx}{"71}
\DeclareMathAlphabet{\mathscrbf}{OMS}{mdugm}{b}{n}

\newcommand{\thickline}{\noalign{\hrule height 1pt}}

 \numberwithin{equation}{section}

\newtheorem{theorem}{Theorem}[section]
\newtheorem{lem}[theorem]{Lemma}
\newtheorem{prop}[theorem]{Proposition}
\newtheorem{corol}[theorem]{Corollary}

\theoremstyle{definition}
\newtheorem{defi}[theorem]{Definition}
\newtheorem{example}[theorem]{Example}
 \newtheorem{rem}[theorem]{Remark}
 
 \theoremstyle{remark}

 \def\bd{\begin{defi}}
\def\ed{\end{defi}}
\def\bt{\begin{theorem}}
\def\et{\end{theorem}}
\def\bl{\begin{lem}}
\def\el{\end{lem}}
\def\bp{\begin{prop}}
\def\ep{\end{prop}}
\def\br{\begin{rem}}
\def\er{\end{rem}}
\def\bc{\begin{corol}}
\def\ec{\end{corol}}
\def\bex{\begin{example}}
\def\eex{\end{example}}
\def\pr{\begin{proof}}
\def\pro{\end{proof}}
\def\eqna{\begin{eqnarray*}}
\def\eqnaa{\begin{eqnarray}}
\def\deqna{\end{eqnarray*}}
\def\deqnaa{\end{eqnarray}}

\definecolor{dark}{rgb}{0.18,0.18,0.68}
\definecolor{mydark}{rgb}{0.78,0.08,0.08}
\definecolor{crew}{rgb}{0.2,0.5,0.2}
\definecolor{mmg}{rgb}{0.31,0.50,0.23}
\definecolor{dblue}{rgb}{0.01,0.01,0.44}
\definecolor{red}{rgb}{0.57,0.11,0.15}
\definecolor{cobalt}{rgb}{0.08,0.14,0.28}
\usepackage[colorlinks,citecolor=cobalt,linkcolor=cobalt,urlcolor=cobalt,pdfpagemode=UseNone,backref = page]{hyperref}

 \input ulem.sty

\begin{document}
\title[Adapted connections with skew-torsion on metric $f$-manifolds]{Adapted connections with skew-torsion on metric $f$-manifolds}

\author{Aleksandra Bor\'owka} 
\address{Institute of Mathematics, Jagiellonian University in Krak\'ow, Poland}
\email{aleksandra.borowka@uj.edu.pl}

\author{Ioannis Chrysikos} 
\address{Department of Mathematics and Statistics,
Faculty of Science, Masaryk University, Kotl\' a\v rsk\' a 2, 611 37 Brno, Czech Republic}
\email{chrysikos@math.muni.cz}

\maketitle
\begin{abstract}
We  provide a natural higher-dimensional generalization of the adapted connections with skew-torsion on almost Hermitian  and almost contact metric manifolds presented in \cite{FrIv}. We prove that  a metric $f$-manifold $(M^{2n+s}, \phi, \xi_i, \eta_j, g)$ with commuting  characteristic vector fields admits a metric connection $\nabla$ with skew-torsion $T$ preserving the structure if and only if each Reeb vector field $\xi_i$ is Killing and the associated  Nijenhuis tensor   is totally skew-symmetric.
This connection is   uniquely determined by its  torsion 3-form $T$,  for which we derive  its explicit formula.   We further  establish necessary and sufficient conditions for  contact metric $f$-manifolds  to admit such a connection. To this end,  for $s\geq 2$ we construct a   broad new class of geometries with parallel skew-torsion in terms of $\mc{S}$-manifolds, i.e., higher-dimensional analogues of Sasakian manifolds.
 We illustrate our results by presenting examples based on low-dimensional   Lie groups.  
  \end{abstract}

 \medskip

 \smallskip


 \tableofcontents

%%%%%%%%%%%%%%%%%%%%%%
%%%%%%%Intro%%%%%%%%%%%%%
%%%%%%%%%%%%%%%%%%%%%%%  

\section{Introduction}
\subsection{Geometries with skew-torsion} 
Geometries with skew-torsion have  attracted considerable attention in recent years. 
Such geometries arise on (pseudo) Riemannian manifolds  that admit a metric connection $\nabla$ whose torsion  $T$ is totally skew-symmetric, i.e., a  3-form. 
 According to Cartan, such geometries constitute one of the three fundamental classes of structures that extend the classical framework of Riemannian geometry by incorporating torsion, yet still retaining enough symmetry to exhibit rich geometric and algebraic properties.
It is now known that connections with skew-torsion play a significant role in the theory of non-integrable  $\Gg$-structures, as they often preserve the  corresponding geometric structure and thus serve as  a natural replacement of the Levi-Civita connection. 
 
 When the  torsion 3-form $T$ is parallel with respect to the connection $\nabla$, we speak for \textsf{geometries with parallel skew-torsion}. Many well-known examples fall into this category, such as naturally reductive homogeneous spaces, nearly K\"ahler manifolds,  nearly parallel $\Gg_2$-manifolds,  Sasakian manifolds, 3-Sasakian manifolds and other; see also the summary in \cite[Section 2.2]{CMS21}.
The $\nabla$-parallelism of $T$ makes the geometry more rigid and tractable, and moreover allows the interplay with   holonomy  theory, see for example \cite{FrIv, CS04, AF04,  AFF15, AD20, CMS21} and for further details on the theory of connections with totally skew-symmetric torsion tensor we refer to \cite{Srni}.   According to  \cite{AFF15, S18},  geometries with parallel skew-torsion are also deeply connected with the  hard classification problem of naturally reductive spaces.   
It is finally established that such geometries  naturally  arise  in theoretical physics and in particular in string theory,  where (parallel) skew-torsion structures yield supersymmetric solutions to the string equations of motion
 (see  \cite{IvP01, FrIv, Srni} and the references therein).

\subsection{Metric $f$-structures} \label{metric-f-structure} In this work we investigate \textsf{adapted metric connections with skew-torsion} on \textsf{Riemannian  metric \negthinspace$f$-manifolds}. 
Recall that  an  \negthinspace\textsf{$f$-structure}  on a  $(2n+s)$-dimensional smooth manifold $M$    is a  non-vanishing  tensor field $\phi$ of type $(1, 1)$ satisfying the identity $\phi^3+\phi=0$ and having constant rank $2n$.
   Any such  tensor field $\phi$  induces a  natural splitting of the tangent bundle,  
 \[
 TM=\Imm(\phi)\oplus\Ker(\phi)
 \]
 and the restriction of $\phi$ to the subbundle $\mc{D}=\Imm(\phi)$ defines an almost complex structure on $\mc{D}$. This means that $M$  naturally exhibits the structure of an almost CR-manifold of CR-dimension $n$ and  CR-codimension $s=\dim\Ker(\phi)$. 
 The notion  of   $f$-structures was introduced by  Yano  \cite{Yano63} and  provides a natural higher-dimensional analogue of the notion of  almost complex  and  almost contact structures,  corresponding  to cases $s=0$ and $s=1$, respectively, see also \cite{R85}. 

 An interesting situation emerges when the subbundle $\Ker(\phi)$ is assumed to be parallelizable, 
   meaning  that there exist global vector fields $\xi_1, \ldots, \xi_s\in\fr{X}(M)$, called the \textsf{characteristic vector fields}, and  1-forms $\eta_1, \ldots, \eta_s$  satisfying 
\[
\Ker(\phi)=\langle \xi_1, \ldots, \xi_s\rangle\,,\quad   \eta_{i}(\xi_j)=\delta_{ij}\,, \quad \eta_{i}\circ\phi=0\,,\quad  \phi^2=-\Id+\sum_{i=1}^{s}\eta_i\otimes\xi_i\,,
\]
where $i, j\in\{1, \ldots, s\}$.
This   is equivalent with a reduction of   the frame bundle of $M^{2n+s}$   to the Lie group $\U(n)\times\Id_{s}$.  Hence in this article we  focus on  $(2n+s)$-dimensional smooth manifolds $M$  endowed with a $\Gg$-structure for the Lie group $\Gg=\U(n)\times\Id_{s}$.   Such manifolds, denoted by  $(M^{2n+s}, \phi, \xi_i, \eta_j)$ and commonly referred to as \textsf{globally framed $f$-manifolds},  have been  widely studied; see for~example \cite{IsY64, BL69, GY70, Blair70, KT72, G72,  BLY73, HOA86, CFF90, DIP01, LP04, DL05, TK07, BP08}. 
In the literature they also appear under alternative names, including  \textsf{$f$-manifolds with complemented frames} (cf. \cite{Blair70, KT72, CFF90}), and \textsf{$f$-manifolds with parallelizable kernel} (cf. \cite{DIP01, BP08}).  

Any globally framed $f$-manifold $(M^{2n+s},  \phi, \xi_i, \eta_j)$ admits an adapted Riemannian metric, that is, a Riemmanian metric $g$ 
satisfying the relation  
 \begin{equation}\label{gXgY}
g(\phi X, \phi Y)=g(X, Y)-\sum_{i=1}^{s} \eta_{i}(X)\eta_{i}(Y)\,,\quad   \ X, Y\in\fr{X}(M)\,.
 \end{equation}
A \textsf{metric \negthinspace $f$-manifold}   is then defined as  a globally framed \negthinspace$f$-manifold $(M^{2n+s}, \phi, \xi_i, \eta_j)$ equipped with  such a Riemannian metric $g$. 
%Note that for $s=1$  a metric $f$-manifold is precisely an  almost contact metric manifold $(M^{2n+1}, \phi, \xi, \eta, g)$, see  \cite{Blair}. 
Metric $f$-structures  serve as  higher-dimensional analogues of  almost Hermitian structures and  almost contact metric structures, a well-established  fact due to the works of Goldberg and Yano, see \cite{GY70,  G72}. 
In 70's,  Blair \cite{Blair70} 
 introduced the notion of \textsf{$\mc{K}$-structures}, as those  metric $f$-structures with closed fundamental 2-form $F(X, Y)=g(X, \phi(Y))$   satisfying a normality condition, $N^{(1)}=0$ identically, where    $N^{(1)}$ is the Nijenhuis tensor.  Obviously,  for $s=0$  we obtain the class of  K\"ahler manifolds,  while  $\mc{K}$-manifolds  for $s=1$ coincide with quasi-Sasaki manifolds (cf. \cite{DIP01}). 
 There are two special subclasses of $\mc{K}$-manifolds,   the  so-called  \textsf{$\mc{S}$-structures} and \textsf{$\mc{C}$-structures}, which  provide higher-dimensional generalizations of Sasakian and cosymplectic manifolds, respectively, 
 see  Section \ref{Definitions}.

%\vskip -1cm

\subsection{Outline}  
In this work, we study general Riemannian metric \negthinspace$f$-structures, with emphasis on their  adapted connections.
In particular, we present
   connections with skew-torsion adapted  to such   Riemannian structures, which, to the best of our knowledge,    were previously unknown  for  $s>1$.
   We then explore some fundamental properties of the induced geometry. 
 
To be more precise,   let $(M^{2n+s},  \phi, \xi_i, \eta_j, g)$  be a metric \negthinspace$f$-manifold     
 whose characteristic
vector fields  $\xi_i$ commute, i.e., $[\xi_i, \xi_j]=0$ for all $i, j\in\{1, \ldots, s\}$.  In Section \ref{Section2} we show that $M^{2n+s}$ admits  
a connection  $\nabla=\nabla^{g}+\frac{1}{2}T$  with skew-torsion $T\in\Gamma(\Lambda^{3}T^*M)$ preserving the structure, i.e.,
\[
\nabla g=0=\nabla\phi\,,\quad \nabla \xi_i=0=\nabla\eta_i\,,\quad i=1, \ldots, s\,,
\]
if and only if each $\xi_i$ is a Killing vector field  and the Nijenhuis tensor $N^{(1)}$ is totally skew-symmetric. If this is the case,  then the
 connection $\nabla$ is  uniquely determined and its   torsion 3-form $T$  is given by
 \[
T=\sum_{i=1}^{s}\eta_{i}\wedge\dd\eta_i+\dd^{\phi}F+N^{(1)}-\sum_{i=1}^{s}(\eta_{i}\wedge(\xi_i\lrcorner N^{(1)}))\,,
\]
with $\dd^{\phi}F:=-\dd F\circ\phi$. This is the content of  Corollary \ref{mainCorol}.

Note that the condition $[\xi_i, \xi_j]=0$ for all $i, j\in\{1, \ldots, s\}$ is not very restrictive, as it holds for \textsf{normal metric $f$-manifolds}, \textsf{contact metric $f$-manifolds},  and in particular \textsf{$\mc{S}$-manifolds}, see below.  Our result should be viewed as an extension of 
 the   well-known statements of Friedrich and Ivanov  \cite{FrIv}, concerning  connections with skew-torsion on
almost Hermitian manifolds and almost contact metric manifolds in dimensions $2n$ and $2n+1$, respectively (see Theorem \ref{FrIvTHM}).   
Our approach for proving this extension  generalizes the methods applied in \cite{FrIv}, however, the computations involved for $s>1$ are in many cases more intricate,
see   the proofs of  
Theorems \ref{Main_Thm_1} and \ref{mainThm2}.

Next,  in Section \ref{Section3}, we  study  several special classes of metric $f$-manifolds,  including  normal metric $f$-manifolds,  $\mc{K}$-manifolds and $\mc{S}$-manifolds. 
For these  classes the adapted connection with skew-torsion takes a simpler form, see for example Propositions  \ref{normal_skewtorsion} and \ref{K1}.  The main result of this section  shows  that the class of $\mc{S}$-manifolds of CR-codimension $s\geq 2$ gives rise to a broad family of geometries with parallel skew-torsion in both odd and even dimensions, see Theorem \ref{S-mnfds1}.  This result provides a higher-dimensional analogue of the well-known  fact that on a  Sasakian manifold $(M^{2n+1},  \phi, \xi, \eta, g)$ the unique adapted connection with  skew-torsion $T$ satisfies the condition $\nabla T=0$ (\cite{FrIv}).  Another notable distinction between the odd-dimensional   Sasakian case $(s=1)$  and a general $\mc{S}$-manifold of higher CR-codimension $(s\geq 2)$  is that, in the latter case,  the torsion 3-form 
$T$  is degenerate,  see Proposition \ref{degeneracy}.  Note  that geometries with skew-torsion which have non-trivial  but degenerate  torsion 3-form $T$ do exist. For instance, in dimension 4 any such geometry is of this type, see also \cite{AFF15}.

Finally, drawing on our constructions, in Section \ref{Section4} we describe  new examples of geometries with skew-torsion in both odd and even dimensions, all    having CR-codimension $s>1$.
We fist study the Lie group $\U(2)$ equipped with a left-invariant metric 
$g$ and viewed as an $\mc{S}$-manifold of CR-codimension 2,  see \cite{TK07}.  We  compute the curvature tensor $R^{\nabla}$ of the adapted connection $\nabla$ with skew-torsion and show  that the corresponding holonomy algebra $\hol(\nabla)$ is 1-dimensional and abelian. In particular, we deduce that $R^{\nabla}$ is $\nabla$-parallel,  $\nabla R^{\nabla}=0$, thus the triple $(\U(2), g, \nabla)$ is an Ambrose-Singer manifold. 
 This  implies that the universal covering of $\U(2)$ is a 4-dimensional naturally reductive space, a result known by Kowalski and  Vanhecke \cite{KV83}. %We further show that $\nabla$ does {\it not} satisfy the $\nabla$-Einstein condition. 
Next we focus on the product $M^6=H_3\times\Tg^3$, where $H_3$ is the 3-dimensional Heisenberg group and $\Tg^3$  is the 3-dimensional torus. By  a construction in \cite{DL05}, our product $M^6$ serves  as an example of an $\mc{S}$-manifold of CR-codimension 4. In this example  the curvature of $\nabla$ is also parallel, $\nabla R^{\nabla}=0$,   hence $M^6$  is an Ambrose-Singer manifold as well. For both these examples we further describe the relation between the Ricci tensors corresponding to the connections  $\nabla$ and  $\nabla^g$ and compute the associated  scalar curvatures.
We finally treat  the 9-dimensional Lie group $\U(3)$, where we show that it can be seen as a normal metric $f$-manifold $(\U(3), \phi, \xi_i, \eta_i, g)$ of CR-codimension 3.  In this case we  prove  that the three Reeb vector fields are   Killing, but  the fundamental 2-form $F$ is not closed. Hence,  $(\U(3), \phi, \xi_1, \xi_2, \xi_3, \eta_1, \eta_2, \eta_3, g)$  does not provide an example of a $\mc{K}$-contact metric $f$-manifolds, but we are able to describe  the adapted connection with skew-torsion in  full detail.

To set the stage, we begin in Section \ref{Preliminaries} by  collecting essential facts about globally framed  $f$-structures and metric $f$-manifolds.
 In the same section % (see Section \ref{Use_lemmas}) 
 we present useful   results regarding Lie derivatives, covariant derivatives and the Nijenhuis tensor. % all useful for the description given in the subsequent sections.

\smallskip
\noindent\textbf{Acknowledgements.}  The first author  has been supported by the Polish National Science Center, project number 2022/47/D/ST1/02197.
The second author acknowledges the support of the Czech Grant Agency (project GA24-10887S), and the Horizon 2020 MSCA project CaLIGOLA,   ID 101086123.  He also thanks   the  Institute of Mathematics at Jagiellonian University in Krak\'ow, for its kind hospitality during a research visit in May 2025, as well as the   Institute of Mathematics at the University of  Z\"urich, where part of this paper was written during his academic stay in the  Winter/2025 semester. The authors would like to thank Anton Galaev and Jan Gregorovi\v{c} for  many  helpful discussions regarding this work.

\section{Globally framed $f$-structures}\label{Preliminaries}

\subsection{Preliminaries}\label{Definitions} 
To start, we fix the conventions that will be used throughout this paper.
%\br\label{conventions} 
In the sequel  we will work with connected smooth manifolds $M$.   
%For a vector bundle $E\to M$ we will denote by $\Gamma(E)$  the set of its smooth sections. For the case of the tangent bundle, $E=TM$, we will mainly use the notation $\fr{X}(M)(=\Gamma(TM))$. 
For the differential of a 1-form $\eta$ on $M$ we fix the convention  $\dd\eta (X, Y)=X\eta(Y)-Y\eta(X)-\eta([X, Y])$ for all $X, Y\in\fr{X}(M)$. 
Also, for a 2-form $\Phi$ we have
\[
\dd\Phi(X, Y, Z)=X\Phi(Y, Z)+Y\Phi(Z, X)+Z\Phi(X, Y)-\Phi([X, Y], Z)-\Phi([Z, X], Y)-\Phi([Y, Z], X)\,.
\]
These conventions are the same as those used in  \cite{FrIv},  but differ from those in  \cite{Blair}, for example.
For a tensor field $\phi\in\Gamma(\Ed(TM))$ we will  often  write $\phi X$ instead of $\phi(X)$.
%\er

Let  $(M^{2n+s}, \phi, \xi_i, \eta_j)$ be a  globally framed $f$-manifold, as defined in Section \ref{metric-f-structure}.  
  In the sequel, following  the work of Goldberg \cite[Section~2]{G72} (see also \cite[p.~94]{DIP01}), it will be convenient to include the trivial case $s=0$ in our discussion.
 As  mentioned in \cite[p.~345]{G72}, strictly speaking, this means that the indices $i, j$ should run through $0, 1, \ldots, s$ and set $\xi_{0}=0=\eta_{0}$.  However,  as is common in the literature and in order to simplify our exposition, we will instead follow the  convention $i, j\in\{1, \ldots, s\}$.  For  $s=0$ the relevant results follow by simply omitting terms depending on $\xi_j$ or their dual 1-forms, where appropriate.  
Let   $TM=\Imm(\phi)\oplus\Ker(\phi)$ be the  splitting  of $TM$  induced by $\phi$, with $\Ker(\phi)=\langle\xi_1, \ldots, \xi_s\rangle$. 
  Since $\phi|_{\Imm(\phi)}$ is an almost complex structure on the subbundle  $\mc{D}:=\Imm(\phi)$,
one can associate in a natural manner 
an almost CR-structure   $(\mc{D}, J=\phi|_{\mc{D}})$ on $M^{2n+s}$  of CR-dimension $n$ and CR-codimension $s=\dim M-2n$. 
 When $s=1$ we obtain the standard almost CR-structure $(\mc{D}, J)$ of hypersurface type  
associated to any almost contact manifold $(M^{2n+1}, \phi, \xi, \eta)$; in this case $\mc{D}=\Ker(\eta)=\Imm(\phi)$ is the distribution of rank $2n$ transversal to the Reeb vector field $\xi$, see for example \cite{T89, Miz93}.
   
The  Nijenhuis tensor  of a globally framed $f$-manifold $(M^{2n+s}, \phi, \xi_i, \eta_j)$ is defined by   
\[
N^{(1)}:=N_{\phi}+\sum_{i=1}^{s}\dd\eta_i\otimes\xi_i
\]
where $
N_{\phi}(X, Y):=[\phi X, \phi Y]+\phi^{2}[X, Y]-\phi[X, \phi Y]-\phi[\phi X, Y]$, 
for any $X, Y\in\fr{X}(M)$ (see \cite[p.~364]{GY70}). %Obviously, for $s=0$, we have $N^{(1)}=N_{\phi}$, and in this case $N^{(1)}$ 
%is the  Nijenhuis tensor associated to the almost complex structure $J=\phi$ on $M^{2n}$, i.e., $N^{(1)}=N_{J}$.  
A  globally framed $f$-structure $(\phi, \xi_i, \eta_j)$  for which the  Nijenhuis tensor  vanishes   identically is called \textsf{normal}. 
%Then  $(M^{2n+s}, \phi, \xi_i, \eta_j)$ is said to be a \textsf{normal globally framed $f$-manifold}.  For $s=0$, the vanishing of $N^{(1)}=N_{J}$  implies that  the almost complex structure $J=\phi$ is integrable, and hence in this case  $(M^{2n}, J)$ is a complex manifold.
%We refer to \cite{HOA86} for a geometric interpretation of the normality on a   globally framed $f$-manifold with $s\geq 1$ (see also \cite{Blair} for $s=1$). 
For a normal globally framed $f$-manifold  $(M^{2n+s}, \phi, \xi_i, \eta_j)$ it is known  that
\[
[\xi_i, \xi_j]=0\,,\quad \mc{L}_{\xi_i}\eta_{j}=0=\mc{L}_{\xi_i}\phi\,,\quad \dd\eta_{i}(\phi X, Y)+\dd\eta_{i}(X, \phi Y)=0\,,
\]
 for all vector fields $X, Y\in\fr{X}(M)$ and $i, j\in\{1, \ldots, s\}$, see \cite{GY70}. 
%\el
%5\bl\label{normal_M} \textnormal{(\cite{GY70})}
 
Let $g$ be an adapted metric on a globally framed $f$-manifold $(M^{2n+s}, \phi, \xi_i, \eta_j)$, that is, a Riemannian metric   satisfying (\ref{gXgY}).
In this case  the decomposition $TM=\Imm(\phi)\oplus\Ker(\phi)$ is orthogonal with respect to $g$.  For simplicity,  we will denote these distributions  by
\[
\mc{D}:=\Imm(\phi)=\bigcap_{i=1}^{s}\Ker(\eta_i)\,,\quad \mc{D}^{\perp}:=\Ker(\phi)=\langle \xi_1, \ldots, \xi_s\rangle\,.
\]
Therefore,  for $X\in\Gamma(\mc{D})$ we have $\eta_{i}(X)=0$  for all $i\in\{1, \ldots, s\}$
and for $X\in\Gamma(\mc{D}^{\perp})$  we get $\phi(X)=0$.
Next, we will refer to vectors belonging to  $\mc{D}$ as {\it horizontal} and vectors belonging to $\mc{D}^{\perp}$ as {\it vertical}. % A general vector field $X\in\fr{X}(M)$  can be expressed  as $X=X^{\mc{D}}+\sum_{k=1}^{s} f_k\xi_k$ for some smooth functions $f_k$ on $M$, with $X^\mc{D}\in\Gamma(\mc{D})$. 
The bundle endomorphism $\phi : TM\to TM$ is   skew-symmetric  with respect to $g$, i.e., 
\begin{equation}\label{gphi}
g(\phi X, Y)+g(X, \phi Y)=0\,,\quad   X, Y\in\fr{X}(M)\,.
\end{equation}
Thus, we can define a 2-form  $F$  on $M^{2n+s}$   by $F(X, Y):=g(X, \phi Y)$, for all $X, Y\in\fr{X}(M)$,  
 called the \textsf{fundamental 2-form} of $(M^{2n+s}, \phi, \xi_i, \eta_j, g)$. For $s=0$ and an almost Hermitian manifold $(M^{2n}, J=\phi, g)$, the fundamental 2-form is referred to as the \textsf{K\"ahler form},   denoted by $\Omega$, i.e., $\Omega(X, Y)=g(X, JY)$.
We have  the relations %$\eta_{i}(X)=g(X, \xi_i)$ for all $i\in\{1, \ldots, s\}$ and moreover
\[
g(\xi_i, \xi_j)=\delta_{ij}, \ (i, j\in\{1, \ldots, s\})\,,\quad \xi_{i}\lrcorner F=0, \ (i=1, \ldots, s)\,, \quad \eta_1\wedge\ldots\wedge\eta_s\wedge F^n\neq 0\,,
\]
%It follows that any metric $f$-manifold  $(M^{2n+s}, \phi, \xi_i, \eta_j, g)$  is orientable. 
and moreover %obviously the fundamental 2-form $F$ satisfies the  identities 
\begin{equation}\label{phi_invariance}
F(\phi X, \phi Y)=F(X, Y)\,,\quad 
 F(\phi X, Y)+F(X, \phi Y)=0\,, \quad X, Y\in\fr{X}(M)\,.
\end{equation}

\bd
(1)  An almost CR-structure   is said to be  \textsf{integrable} (or CR-integrable), 
if  $[JX, Y]+[Y, JX]\in\Gamma(\mc{D})$ and $N_{J}(X, Y)=[JX, JY]-[X, Y]-J[JX, Y]-J[X, JY]=0$ for any $X, Y\in\Gamma(\mc{D})$.\\
(2)   A metric $f$-manifold  $(M^{2n+s}, \phi, \xi_i, \eta_j, g)$ which is normal and has closed fundamental 2-form, i.e., $N^{(1)}=0$  and $\dd F=0$,  is called a   \textsf{$\mc{K}$-manifold}.\\
(3)  A metric $f$-manifold  $(M^{2n+s}, \phi, \xi_i, \eta_j, g)$ whose fundamental 2-form $F$ satisfies $2 F(X, Y)=\dd\eta_{i}(X, Y)$ for all $i=1, \ldots, s$ and  $X, Y\in\fr{X}(M)$
is called a \textsf{contact metric $f$-manifold} (or an  \textsf{almost $\mc{S}$-manifold}, see for example \cite[Definition 1.3]{DIP01} or \cite[Definition 3.1]{BP08}). 
\\
(4)  A  contact metric $f$-manifold which is normal  is called a (Riemannian)  \textsf{$\mc{S}$-manifold}, i.e.,
\[
2F=\dd\eta_1=\cdots=\dd\eta_s\,,\quad  N_{\phi}=-\sum_{i=1}^{s}\dd\eta_{i}\otimes\xi_i\,.
\]
\ed
The class of $\mc{K}$-manifolds was introduced in \cite{Blair70}, where it was demonstrated that for any $\mc{K}$-manifold  $(M^{2n+s},  \phi, \xi_i, \eta_j, g)$ the characteristic vector fields $\xi_1, \ldots, \xi_s$ are  Killing.  Note that for $s=1$ a $\mc{K}$-manifold is a   \textsf{quasi-Sasakian manifold}, i.e.,  a normal almost contact
metric manifold with closed fundamental 2-form.  %$N^{(1)} = 0$ and $\dd F = 0$.  
Moreover, for $s=1$ an $\mc{S}$-manifold is a \textsf{Sasakian manifold}, i.e., an almost contact metric manifold $(M^{2n+1}, \phi, \xi, \eta, g)$ satisfying $N^{(1)}=0$ and $2F=\dd\eta$.
Beyond $\mc{S}$-manifolds, another special class of $\mc{K}$-manifolds called \textsf{$\mc{C}$-manifolds} are those metric $f$-manifolds   $(M^{2n+s}, \phi, \xi_i, \eta_j, g)$ which are normal and satisfy $\dd F=0$ and $\dd\eta_j=0$ for all $j=1, \ldots, s$. %Note that for $s=1$, a $\mc{C}$-manifold is a so-called   \textsf{co-symplectic manifold}.

%%%%%%%%%%%%%%%%%%%%%%%%%%%%%%%%%%
%%%%%%%%useful_results%%%%%%%%%%%%%%%%%%%
%%%%%%%%%%%%%%%%%%%%%%%%%%%%%%%%%%%

\subsection{Results on metric $f$-manifolds}\label{Use_lemmas}
Let us now present several results concerning the geometry of metric $f$-manifolds, many of which will play a crucial role in the subsequent analysis. %We stress that   this section   focuses   on the general setting, without imposing additional conditions, such as normality ($N^{(1)} = 0$) or restrictions on the fundamental 2-form $F$.  Some results deal with the special case where the characteristic vector fields commute and are Killing, a  scenario   that will play a significant role later in this work.

%%%%%%%%%%%%%%%%%%%%%%%%%%%%%%%%%%
%%%%%%%%%Lie derivatives%%%%%%%%%%%%%%%%%%
%%%%%%%%%%%%%%%%%%%%%%%%%%%%%%%%%%%

\subsubsection{On Lie derivatives}
Let  $(M^{2n+s},  \phi, \xi_i, \eta_j, g)$ be  a metric $f$-manifold, and let $\nabla^{g}$  be  the Levi-Civita connection associated to $g$.  To begin with,
 observe that a direct  combination of  the relation $\eta_i(Y)= g(Y, \xi_i)$ and  the metric property of $\nabla^g$ yields the important relation % that %we mention the useful relation
\begin{equation}\label{gen_s}
(\nabla^{g}_{X}\eta_i)Y=g(\nabla^{g}_{X}\xi_{i}, Y)\,,\  i=1, \ldots, s\,,  \quad X, Y\in\fr{X}(M)\,.
\end{equation}
%This follows as a direct  combination of  the relation $\eta_i(Y)= g(Y, \xi_i)$ and  the metric property of $\nabla^g$:
%\begin{eqnarray*}
%(\nabla^{g}_{X}\eta_i)Y&=&\nabla^{g}_{X}\eta_i(Y)-\eta_i(\nabla^{g}_{X}Y)= \nabla^{g}_{X}g(Y, \xi_i)- g(\nabla^{g}_{X}Y, \xi_i)\\
%&=&  g(\nabla^{g}_{X}Y, \xi_i)+ g(Y, \nabla^{g}_{X}\xi_i)- g(\nabla^{g}_{X}Y, \xi_i)=  g(Y, \nabla^{g}_{X}\xi_i)\,,
%\end{eqnarray*}
%for all $X, Y\in\fr{X}(M)$ and   $i=1, \ldots, s$.   
%Let us now analyze  results regarding Lie derivatives.

 \bl\label{Lie_derivatives}
Let $(M^{2n+s},  \phi, \xi_i, \eta_j, g)$ be a  metric $f$-manifold.  Then
for any vector field $X\in\fr{X}(M)$ and any  $i, j\in\{1, \ldots, s\}$ we have 
\begin{equation}\label{Lie_derivative}
(\mc{L}_{\xi_i}\eta_{j})(X)=\dd\eta_j(\xi_i, X)=g(X, \nabla^g_{\xi_i}\xi_j)+g(\nabla^{g}_{X}\xi_i, \xi_j)\,.
\end{equation}
For $i=j$ this reduces to $(\mc{L}_{\xi_i}\eta_{i})(X)=g(X, \nabla^g_{\xi_i}\xi_i)$.  
Moreover $(\mc{L}_{\xi_i}\eta_{j})(\xi_k)=-\eta_{j}([\xi_i, \xi_k])=-g([\xi_i, \xi_k], \xi_j)$, 
while for $X\in\Gamma(\mc{D})$ we obtain $(\mc{L}_{\xi_i}\eta_j)(X)=-\eta_{j}([\xi_i, X])=-g([\xi_i, X], \xi_j)$.
\el

\bc 
Let $(M^{2n+s}, \phi, \xi_i, \eta_j, g)$ be a metric $f$-manifold satisfying $\mc{L}_{\xi_i}\eta_{j}=0$ for all $i, j\in\{1, \ldots, s\}$.
Then $[\xi_i, \xi_j]\in\Gamma(\mc{D})$ and $[X, \xi_i]\in\Gamma(\mc{D})$ for all $X\in\Gamma(\mc{D})$ and $i, j\in\{1, \ldots, s\}$.
\ec

The following proposition deals with the case where the vertical distribution $\mc{D}^{\perp}=\ker(\phi)=\langle\xi_1, \ldots, \xi_s\rangle$ is generated by commuting  Killing vector fields. Its proof is based on the Killing condition, the relations (\ref{gen_s}), (\ref{Lie_derivative}) and \cite[Lemma 3.1]{T04}. 
\bp\label{prop_derivatives}
Let  $(M^{2n+s},  \phi, \xi_i, \eta_j, g)$  be a metric $f$-manifold  
 whose characteristic
vector fields  $\xi_i$ are Killing and commute, i.e., $[\xi_i, \xi_j]=0$ for all $i, j\in\{1, \ldots, s\}$. 
Then  the following hold: \\  
(1)  $\dd\eta_j=2\nabla^g\eta_j$ for all $j\in\{1, \ldots, s\}$.\\
(2) $\mc{L}_{\xi_{i}}\eta_{j}=0$ and $\nabla^g_{\xi_i}\xi_j=0$ for all $i, j\in\{1, \ldots, s\}$.\\
(3) $\xi_j\lrcorner\dd\eta_i=0$ for all $i, j\in\{1, \ldots, s\}$, and in particular  $0=\dd\eta_i(\xi_j, X)=2(\nabla^{g}_{\xi_j}\eta_i)X$ for any $X\in\fr{X}(M)$ and $i, j\in\{1, \ldots, s\}$.\\
(4) $g(\nabla^{g}_{X}\xi_i, \xi_j)=0$ for all $X\in\fr{X}(M)$ and $i, j\in\{1, \ldots, s\}$.
\ep

\br\label{remark_derivatives}
(1)  Let $(M, g)$ be a $(2n+s)$-dimensional Riemannian manifold. Suppose that the  set $\{\xi_1, \ldots, \xi_s\}$ consists of $g$-orthonormal Killing vector fields satisfying $\mc{L}_{\xi_i}\eta_j=0$, for all $i, j\in\{1, \ldots, s\}$, where $\eta_1, \ldots, \eta_s$ are the corresponding dual 1-forms. 
In the proof of \cite[Lemma 3.1]{T04} it was derived that these conditions imply that $[\xi_i, \xi_j]=0$  for all   $i, j\in\{1, \ldots, s\}$.
On the other hand,  part   (2) of Proposition \ref{prop_derivatives} gives  a converse of this statement in our context. 
Therefore,  under the assumption that  the vertical subbundle  $\mc{D}^{\perp}=\langle \xi_1, \ldots, \xi_s\rangle$ is generated by Killing vector fields,  the conditions $[\xi_i, \xi_j]=0$ and   $\mc{L}_{\xi_{i}}\eta_j=0$ for all  $i, j\in\{1, \ldots, s\}$ are   equivalent. Moreover, note that the assumptions of Proposition \ref{prop_derivatives} imply the additional relations $g(\nabla^{g}_{\xi_{i}}\xi_{j}+\nabla^g_{\xi_j}\xi_i, X)=0$ and $\nabla^{g}_{X}\xi_i\in\Gamma(\mc{D})$ for all $X\in\Gamma(\mc{D})$, see also   \cite{T04}.\\
 (2)   Let $(M^{2n+s},  \phi, \xi_i, \eta_j, g)$ be a metric $f$-manifold  satisfying the  assumptions of Proposition \ref{prop_derivatives}. Then $M^{2n+s}$ admits a totally geodesic foliation whose leaves are tangent to $\mc{D}^{\perp}=\langle\xi_1, \ldots, \xi_s\rangle$.
 This is  because the vertical distribution is obviously integrable  and we have $\nabla^{g}_{X}Y\in\Gamma(\mc{D}^{\perp})$ for any $X, Y\in\Gamma(\mc{D}^{\perp})$ (since  $\nabla^{g}_{\xi_{i}}\xi_j=0$
 for all $i, j\in\{1, \ldots, s\}$).  
\er

%%%%%%%%%%%%%%%%%%%%%%%%%%%%%%%%%%
%%%%%%%%%Covariant derivatives%%%%%%%%%%%%%%
%%%%%%%%%%%%%%%%%%%%%%%%%%%%%%%%%%%

\subsubsection{On covariant derivatives} 
We now turn our attention to results concerning covariant derivatives. We begin by recording  the following identities, which are obtained by straightforward computations using Proposition  \ref{prop_derivatives}  and  Remark \ref{remark_derivatives} (see also \cite[Prop.~2.5]{BP08}).
\bl\label{main_lemma}
Let $(M^{2n+s},  \phi, \xi_i, \eta_j, g)$ be a  metric $f$-manifold.  Then  the following hold:
\begin{eqnarray}
(\nabla^g_{X}\phi)\phi Y&=&-\nabla^g_{X}Y+\sum_{i=1}^{s}\eta_{i}(Y)\nabla^{g}_{X}\xi_i+\sum_{i=1}^{s}X(\eta_{i}(Y))\xi_i-\phi(\nabla^g_{X}\phi Y)\,, \label{covariant_phi}\\
 (\nabla^g_{X}F)(Y, Z)&=&g((\nabla^g_{X}\phi)Z, Y)=-g((\nabla^{g}_{X}\phi)Y, Z)=-(\nabla^{g}_{X}F)(Z, Y)\,,  \label{cova_F}\\
\label{gen_s2} (\nabla^{g}_{X}F)(\xi_i, \phi Y)&=&
 \left\{
\begin{tabular}{lr}
$  (\nabla^{g}_{X}\eta_{i})Y+ \sum_{j=1}^{s}\eta_{j}(Y)g(\nabla^{g}_{X}\xi_j, \xi_i)\,,$ & $s>1$,\\
$ (\nabla^{g}_{X}\eta_{i})Y\,,$ & $s=1$,
\end{tabular}\right.
\end{eqnarray}
for any $X, Y, Z\in\fr{X}(M)$.  
If each $\xi_i$ is Killing satisfying $\mc{L}_{\xi_i}\eta_j=0$ for all $i, j\in\{1, \ldots, s\}$, then  
$(\nabla^{g}_{X}F)(\xi_i, \phi Y)=(\nabla^{g}_{X}\eta_{i})Y$  for all $i=1, \ldots, s$.
 \el

Based on this lemma one can now prove that

 \bp\label{symmetries_FN}
 On a metric $f$-manifold  $(M^{2n+s},  \phi, \xi_i, \eta_j, g)$  the covariant derivative $\nabla^gF$ of the fundamental form $F$ satisfies the identity 
\begin{equation}
 (\nabla^{g}_{X}F)(Y, Z)+(\nabla^{g}_{X}F)(\phi Y, \phi Z)=\sum_{i=1}^{s}\eta_{i}(Z)(\nabla^g_{X}\eta_i)\phi Y-\sum_{i=1}^{s}\eta_{i}(Y)(\nabla^{g}_{X}\eta_i)\phi Z\,.
\end{equation} 
\ep

%%%%%%%%%%%%%%%%%%%%%%%%%%%%%%%%
%%%%%%Nijenhuis_tensor_N^(1)%%%%%%%%%%%%%%%
%%%%%%%%%%%%%%%%%%%%%%%%%%%%%%%%%

\subsubsection{On the  Nijenhuis tensor}
%We conclude this preliminary section  by presenting  necessary results regarding  the  Nijenhuis  tensor  $N^{(1)}$. 
Let  $(M^{2n+s},  \phi, \xi_i, \eta_j, g)$ be a metric  $f$-manifold. Set
\[
N^{(1)}(X, Y, Z):=g(N^{(1)}(X, Y), Z)\,,\quad  N^{(2)}_{i}(X, Y):=\dd\eta_{i}(\phi X, Y)+\dd\eta_{i}(X, \phi Y)\,,
\]
for all vector fields $X, Y, Z\in\fr{X}(M)$ and $ i=1, \ldots, s$. 
Since we may write  $\phi^2[X, Y]=-[X, Y]+\sum_{k=1}^{s}\eta_{k}([X, Y])\xi_{k}$, we   see that
\begin{eqnarray}
N^{(1)}(X, Y, Z)
&=&g([\phi X, \phi Y], Z)-g(\phi[\phi X, Y],Z)-g(\phi[X, \phi Y], Z)-g([X, Y], Z)\nonumber\\
&&+\sum_{k=1}^{s}\big\{\eta_{k}([X, Y])+\dd\eta_{k}(X, Y)\big\}\eta_{k}(Z)\,.\label{general_N}  
\end{eqnarray}
This formula will play an important role in our subsequent arguments.  In addition,  for a  metric $f$-manifold $(M^{2n+s}, \phi, \xi_i, \eta_j,  g)$   a direct computation  yields the following identity 
\begin{equation}\label{tensor_N_T} 
  N^{(1)}(X, Y)=(\nabla^{g}_{\phi X}\phi)Y-(\nabla^{g}_{\phi Y}\phi)X+(\nabla^g_{X}\phi)\phi Y-(\nabla^{g}_{Y}\phi)\phi X
+\sum_{j=1}^{s}\big(\eta_{j}(X)\nabla^{g}_{Y}\xi_{j}-\eta_{j}(Y)\nabla^{g}_{X}\xi_{j}\big).
\end{equation}
We will make use of this later, in particular at the end of the proof of Lemma \ref{formulas_N12}, which plays a crucial role in establishing one of our main results, namely Theorem \ref{mainThm2}.   Another result contributing to the proof of Theorem \ref{mainThm2}   is the following:  
\bl\label{basic_1} \textnormal{(\cite[p.~93]{DIP01},  \cite[Prop.~2.4]{BP08})}
Let $(M^{2n+s}, \xi_i, \eta_j, \phi, g)$ be a  metric $f$-manifold.  Then  
\begin{eqnarray*}
2g((\nabla^{g}_{X}\phi)Y, Z)&=& \dd F(X, \phi Y, \phi Z)-\dd F(X, Y, Z)+N^{(1)}(Y, Z, \phi X)\\
&&+\sum_{i=1}^{s}\Big\{\eta_{i}(X)N^{(2)}_{i}(Y, Z)+\eta_{i}(Z)\dd\eta_{i}(\phi Y, X)+\eta_{i}(Y)\dd\eta_{i}(X, \phi Z)\Big\}\,.
\end{eqnarray*}
\el
% We now  present a result which is a generalization of the identity $s=1$  in \cite[p.~327]{FrIv}). 

\bl\label{identity_N}
The  Nijenhuis tensor  of a metric $f$-manifold  $(M^{2n+s},  \phi, \xi_i, \eta_j, g)$   obeys the following relation:  
\begin{eqnarray}
 N^{(1)}(X, Y, Z)&=& -N^{(1)}(\phi X, \phi Y, Z)+\sum_{i=1}^{s}\eta_{i}(X)N^{(1)}(\xi_i, Y, Z)+\sum_{i=1}^{s}\eta_{i}(Y)N^{(1)}(X, \xi_{i}, Z)\nonumber\\
&&+\sum_{i, j=1}^{s}\eta_{i}(X)\eta_{j}(Y)g([\xi_i, \xi_j], Z)\,, \label{Nij_tensor} 
\end{eqnarray}
for all $X, Y, Z\in\fr{X}(M)$. 
%For the special cases $s=0, 1$, or whenever   $[\xi_i, \xi_j]=0$ for all $i, j\in\{1, \ldots, s\}$,    the  second line  of {\rm(\ref{Nij_tensor})} vanishes. 
 \el
  \pr
An application of (\ref{general_N}) for $X, Y$ replaced by $\phi X$ and  $\phi Y$, respectively, yields the  expression
\begin{eqnarray*}
N^{(1)}(\phi X, \phi Y, Z)&=&g([\phi^2X, \phi^2Y], Z)-g([\phi X, \phi Y], Z)-g(\phi[\phi X, \phi^2Y], Z)-g(\phi[\phi^2X, \phi Y], Z)\\
&&+\sum_{k}\big\{\eta_{k}([\phi X, \phi Y])+\dd\eta_{k}(\phi X, \phi Y)\big\}\eta_{k}(Z)\,.
\end{eqnarray*}
Since $\eta_{k}\circ\phi=0$ we have $\dd\eta_{k}(\phi X, \phi Y)=-\eta_{k}([\phi X, \phi Y])$ for all $k=1, \ldots, s$, thus the sum in the second line vanishes.
For the rest terms we compute and get
\begin{eqnarray*}
g([\phi^2X, \phi^2Y], Z)=%&=&g([-X+\sum_{i}\eta_{i}(X)\xi_i, -Y+\sum_{j}\eta_{j}(Y)\xi_j], Z)\\
&=&g([X, Y], Z)-\sum_{j}\eta_{j}(Y)g([X, \xi_j], Z)-\sum_{j}X(\eta_{j}(Y))g(\xi_j, Z)\\
&&-\sum_{i}\eta_{i}(X)g([\xi_i, Y], Z)+\sum_{i}Y(\eta_{i}(X))g(\xi_i, Z)+\sum_{i, j}\eta_{i}(X)\eta_{j}(Y)g([\xi_i, \xi_j], Z)\\
&&+\sum_{i, j}\eta_{i}(X)\xi_{i}(\eta_j(Y))g(\xi_j, Z)-\sum_{j, i}\eta_{j}(Y)\xi_j(\eta_i(X))g(\xi_i, Z)\,,\\
g(\phi[\phi X, \phi^2 Y], Z)&=&-g(\phi[\phi X, Y], Z)+\sum_{j}\eta_{j}(Y)g(\phi[\phi X, \xi_j], Z)\,, \\
g(\phi[\phi^2X, \phi Y], Z)&=&-g(\phi[X, \phi Y], Z)+\sum_{i}\eta_{i}(X)g(\phi[\xi_i, \phi Y], Z)\,,
\end{eqnarray*}
%where for the last two relations we didn't  write terms that vanish by  the relation $\phi(\xi_i)=0$. This gives
%\begin{eqnarray*}
%N^{(1)}(\phi X, \phi Y, Z)&=&g([X, Y], Z)-g([\phi X, \phi Y], Z)+g(\phi[\phi X, Y], Z)+g(\phi[X, \phi Y], Z)\\
%&&-\sum_{j}\eta_{j}(Y)g([X, \xi_j], Z)-\sum_{j}X(\eta_{j}(Y))g(\xi_j, Z)-\sum_{j}\eta_{j}(Y)g(\phi[\phi X, \xi_j], Z)\\
%&&-\sum_{i}\eta_{i}(X)g([\xi_i, Y], Z)+\sum_{i}Y(\eta_{i}(X))g(\xi_i, Z)-\sum_{i}\eta_{i}(X)g(\phi[\xi_i, \phi Y], Z)\\
%&&+\sum_{i, j}\eta_{i}(X)\xi_{i}(\eta_j(Y))\eta_{j}(Z)-\sum_{j, i}\eta_{j}(Y)\xi_j(\eta_i(X))\eta_{i}(Z)\\
%&&+\sum_{i, j}\eta_{i}(X)\eta_{j}(Y)g([\xi_i, \xi_j], Z)\,.
%\end{eqnarray*}
Putting this together and subsituting
\[
g([X, Y], Z)=-g(\phi^2[X, Y], Z)+\sum_{k}\eta_{k}([X, Y])g(\xi_k, Z)=-g(\phi^2[X, Y], Z)+\sum_{k}\eta_{k}([X, Y])\eta_{k}(Z)\,.
\]
we can   reconstruct the part $N_{\phi}$ of   $N^{(1)}$. In particular,  since the relation $N^{(1)}(X, Y, Z)=g(N_{\phi}(X, Y), Z)+\sum_{k}\dd\eta_{k}(X, Y)\eta_{k}(Z)$ is valid for any $X, Y, Z\in\fr{X}(M)$,    
we obtain that
\begin{eqnarray*}
N^{(1)}(\phi X, \phi Y, Z)&=&-N^{(1)}(X, Y, Z)+\sum_{k}\dd\eta_{k}(X, Y)\eta_{k}(Z)+\sum_{k}\eta_{k}([X, Y])\eta_{k}(Z)\\
&&-\sum_{j}\eta_{j}(Y)g([X, \xi_j], Z)-\sum_{j}X(\eta_{j}(Y))g(\xi_j, Z)-\sum_{j}\eta_{j}(Y)g(\phi[\phi X, \xi_j], Z)\\
&&-\sum_{i}\eta_{i}(X)g([\xi_i, Y], Z)+\sum_{i}Y(\eta_{i}(X))g(\xi_i, Z)-\sum_{i}\eta_{i}(X)g(\phi[\xi_i, \phi Y], Z)\\
&&+\sum_{i, j}\eta_{i}(X)\xi_{i}(\eta_j(Y))\eta_{j}(Z)-\sum_{j, i}\eta_{j}(Y)\xi_j(\eta_i(X))\eta_{i}(Z)\\
&&+\sum_{i, j}\eta_{i}(X)\eta_{j}(Y)g([\xi_i, \xi_j], Z)\,.
\end{eqnarray*}
It is easy to see that in this expression the sum of the second and third term cancels by the sum of the fifth and eighth term, as we can consider these sums with respect to a common index. It remains to use the expressions of $N^{(1)}(\xi_i, Y, Z)$ and $N^{(1)}(X, \xi_j,  Z)$, which are respectively given by
\begin{eqnarray*}
N^{(1)}(\xi_i, Y, Z)%&=&-g(\phi[\xi_i, \phi Y], Z)-g([\xi_i, Y], Z)+\sum_{k}\big\{\eta_{k}([\xi_i,  Y])+\dd\eta_{k}(\xi_i, Y)\big\}\eta_{k}(Z)\\
&=&-g(\phi[\xi_i, \phi Y], Z)-g([\xi_i, Y], Z)+\sum_{k}\xi_{i}(\eta_{k}(Y))\eta_{k}(Z)\,,\\
N^{(1)}(X, \xi_j,  Z)%&=&-g(\phi[\phi X, \xi_j], Z)-g([X, \xi_j], Z)+\sum_{k}\big\{\eta_{k}([X, \xi_j])+\dd\eta_{k}(X, \xi_j)\big\}\eta_{k}(Z)\\
&=&-g(\phi[\phi X, \xi_j], Z)-g([X, \xi_j], Z)-\sum_{k}\xi_{j}(\eta_{k}(X))\eta_{k}(Z)\,,
\end{eqnarray*}
for all $i, j\in\{1, \ldots, s\}$. The stated formula is now readily obtained.
\pro

%%%%%%%%%%%%%%%%%%%%%%%%%%%%%%%%%
%%%%%%%%%%%%%%%%%%%%%%%%%%%%%%%%%
%%%%%%Connections_Skew_torsion%%%%%%%%%%%%%
%%%%%%%%%%%%%%%%%%%%%%%%%%%%%%%%%

\section{Adapted connections with skew-torsion}\label{Section2}
\subsection{Connections with skew-torsion}\label{subparagraph_Skew_Torsion}
Let $(M, g)$ be a Riemannian manifold and let $\nabla^g$ be the Levi-Civita connection associated to  the Riemannian metric $g$. 
 Recall that the torsion $T$ of an affine connection $\nabla$ on $M$  is the vector-valued 2-form  defined by
 $T(X, Y)=\nabla_{X}Y-\nabla_{Y}X-[X, Y]$, for all $X, Y\in\fr{X}(M)$.
The connection $\nabla$     is said to have {\it totally skew-symmetric torsion}, or  {\it skew-torsion} in short,
if  the induced tensor field
\[
T(X, Y, Z):=g(T(X, Y), Z)\,,\quad X, Y, Z\in\fr{X}(M)
\]
is a 3-form on $M$. 
Then we may express $\nabla$ as 
\[
g(\nabla_{X}Y, Z)=g(\nabla^{g}_{X}Y, Z)+\frac{1}{2}g(T(X, Y), Z)=g(\nabla^{g}_{X}Y, Z)+\frac{1}{2}T(X, Y, Z)\,.
\]
 Such a connection  is necessarily metric, $\nabla g=0$.
 Moreover, the $\nabla$-Killing vector fields coincide with the Riemannian Killing vector fields, see \cite{AF04, Srni}.   
Note that for any 3-form $T$ on $M$
we can consider the 4-form
\[
\sigma_{T}:=\frac{1}{2}\sum_{i=1}^{\dim M}(e_i\lrcorner T)\wedge(e_i\lrcorner T), 
\]
where $\{e_i\}$ is a  local orthonormal frame of  $(M, g)$. This 4-form is also expressed as 
\[
\sigma_{T}(X, Y, Z, W)= g(T(X, Y), T(Z, W))+g(T(Y, Z), T(X, W))+g(T(Z, X), T(Y, W))
\]
for all $X, Y, Z, W\in\fr{X}(M)$,  and plays an important role in the theory of metric connections with skew-torsion, see for example \cite{AF04, Srni, AFF15}. We will meet this 4-form later. We finally recall
that for the norm of the 3-form $T$ one  can adopt the normalization  
\begin{equation}\label{normT}
\|T\|^2=\frac{1}{6}\sum_{i, j=1}^{\dim M}\|T(e_i, e_j)\|^2\,. 
\end{equation}

\subsection{Connections with skew-torsion on metric $f$-manifolds}\label{metric-phi-connections}
Next we will study metric $f$-manifolds $(M^{2n+s},  \phi, \xi_i, \eta_j, g)$  
which admit an adapted connection $\nabla : \Gamma(TM)\to\Gamma(T^*M\otimes TM)$ whose torsion $T$ is totally skew-symmetric. Here, by the term ``adapted''  we mean that $\nabla$   preserves the structure, i.e., 
\begin{equation}\label{parallel_objects}
\nabla g=0\,,\quad \nabla\phi=0\,,\quad \nabla\xi_i=0\,, \quad \nabla\eta_i=0\,,\quad \forall \ i=1,\ldots, s\,.
\end{equation}
Of course, for $s=0$ and an almost Hermitian manifold $(M^{2n}, J=\phi, g)$  one  only requires $\nabla g=0=\nabla J$. 
Next we will adopt the following definition.
\bd
On a  metric $f$-manifold $(M^{2n+s},  \phi, \xi_i, \eta_j, g)$  an affine connection $\nabla$ satisfying  the conditions in (\ref{parallel_objects}) 
is called  an  \textsf{almost metric $\phi$-connection}.
\ed
Since we are interested in the case where the torsion $T$ of $\nabla$  is a 3-form,  
we will often write   $\nabla=\nabla^g+\frac{1}{2}T$ and refer to $\nabla$ as an  \textsf{almost metric $\phi$-connection with skew-torsion}. 
A connection with totally skew-symmetric torsion that preserves a geometric structure is also referred to as a   \textsf{characteristic connection}, see \cite{Srni}. 
 Friedrich and Ivanov  \cite{FrIv} presented the  (unique) characteristic connections for the cases $s=0, 1$. Let us recall the statement for the case $s=1$, i.e., for an almost contact metric manifold $(M^{2n+1}, \phi, \xi, \eta, g)$. For $s=0$, i.e., an almost Hermitian manifold $(M^{2n}, J=\phi, g)$ the result  occurs by Theorem \ref{FrIvTHM} by dropping the redundant condition on  $\xi$, while the formula of the torsion 3-form $T$ is  obtained by the given formula by replacing $N^{(1)}$ by $N_{J}$, $\dd^{\phi}F$ by $\dd^{J}\Omega$, and   deleting   the two non-relevant  terms including the 1-form $\eta$, see \cite[Theorem 10.1]{FrIv}. 
\bt\label{FrIvTHM} \textnormal{(\cite[Theorem 8.2]{FrIv})} Let $(M^{2n+1}, \phi, \xi, \eta, g)$ be an almost contact metric manifold with  fundamental 2-form denoted by $F$. Then $M$ admits a metric connection $\nabla$ with skew-torsion such that
$\nabla g=\nabla\phi=\nabla\eta=\nabla\xi=0$ if and only if the Nijenhuis tensor $N^{(1)}$ is totally skew-symmetric and $\xi$ is Killing.  In this case the connection $\nabla$ is uniquely determined and its torsion
is given by
\[
T=\eta\wedge\dd\eta+\dd^{\phi}F+N^{(1)}-\eta\wedge(\xi\lrcorner N^{(1)})\,,
\]
where   $\dd^{\phi} F(X, Y, Z):=-\dd F(\phi X, \phi Y, \phi Z)$ for any  $X, Y, Z\in\fr{X}(M)$. 
\et

Below our goal is to generalize the above theorem for metric $f$-manifolds $(M^{2n+s},  \phi, \xi_i, \eta_j, g)$ with   $s>1$.
  We mention that the $\nabla$-parallelism of   $\phi$ and of $\eta$ (or equivalently of $\xi$)    implies that the characteristic connection $\nabla$ in Theorem \ref{FrIvTHM}
preserves  both the distributions     $\mc{D}= \langle\xi\rangle^{\perp}$ and $\mc{D}^{\perp}= \langle\xi\rangle$.  Later on we would like to maintain this scenario and hence we  may need to  consider additional assumptions, as for example $\mc{L}_{\xi_i}\xi_j=[\xi_i, \xi_j]=0$ for all $i, j\in\{1, \ldots, s\}$, see below.
 
 We begin with the following proposition, which  should be understood as a  higher-dimensional generalization    
 of  statements that appear in the proof of the forward implication of Theorem \ref{FrIvTHM}, given in 
\cite{FrIv} for the almost contact metric case $(s=1)$.

%%%%%%%%%%%%%%%%%%%%%%%%%%%%%%%%
%%%%%%%%Torsion_Characterization%%%%%%%%%%%%
%%%%%%%%%%%%%%%%%%%%%%%%%%%%%%%%%

\bp\label{characterize_torsion}
Let  $(M^{2n+s},  \phi, \xi_i, \eta_j, g)$  be a metric $f$-manifold    and let $\nabla=\nabla^g+\frac{1}{2}T$
be a metric connection with skew-torsion on $M$ preserving the structure, that is, an almost metric $\phi$-connection with skew-torsion. 
Then the following hold:
\begin{enumerate}
\item  The conditions $\nabla\xi_i=0=\nabla\eta_i$  are equivalent for all $i=1, \ldots, s$.
\item   Each $\xi_i$ is a Killing vector field.  
\item  For all $i=1, \ldots, s$ we have $\dd\eta_i=\xi_i\lrcorner T=2\nabla^g\eta_i$, that is, 
\[
\dd\eta_{i}(X, Y)=T(\xi_i, X, Y)=2(\nabla^g_{X}\eta_i)Y
\]
for all $X, Y\in\fr{X}(M)$. Thus also $\xi_i\lrcorner\dd\eta_i=0$ for all $i=1, \ldots, s$.
\item The Nijenhuis tensor satisfies the identity 
\[
N^{(1)}(X, Y, Z)=-T^{-}(X, Y, Z)\,,\quad X, Y, Z\in\fr{X}(M)\,,
\]
where  $T^{-}(X, Y, Z):=T(X, \phi Y, \phi Z)+T(\phi X, Y, \phi Z)+T(\phi X, \phi Y, Z)-T(X, Y, Z)$. 
Hence $N^{(1)}$ is totally skew-symmetric. 
\item  The covariant derivative of the fundamental form $F$ with respect to $\nabla^g$ satisfies the identity 
\[
2(\nabla^{g}_{X}F)(Y, Z)= -2g((\nabla^g_{X}\phi)Y, Z)= T(X, Y, \phi Z)+T(X, \phi Y, Z)\,, \quad X, Y, Z\in\fr{X}(M)\,.
\]
\item The differential of the fundamental 2-form $F$ is the 3-form given by  the following cyclic sum:
\[
\dd F(X, Y, Z)=\fr{S}_{X, Y, Z}T(X, Y, \phi Z)\,, \quad X, Y, Z\in\fr{X}(M)\,.
\]
\item For any $X, Y, Z\in\fr{X}(M)$ we have 
\begin{equation}\label{TxTyTz}
T(\phi X, \phi Y, \phi Z)=\dd F(X, Y, Z)-N^{(1)}(X, Y, \phi(Z))-\sum_{i=1}^{s}\eta_{i}(Z)N_{i}^{(2)}(X, Y)\,.
\end{equation}
\end{enumerate}
\ep
\pr
(1)  By the identity $\eta_i(X)=g(X, \xi)$, for all $X, Y\in\fr{X}(M)$  and  $i\in\{1, \ldots, s\}$,  we have
\[
(\nabla_{X}\eta_i)Y=(\nabla_{X}^{g}\eta_i)Y-\frac{1}{2}\eta_{i}(T(X, Y))=(\nabla_{X}^{g}\eta_i)Y-\frac{1}{2}g(T(X, Y), \xi_i)\,.
\]
Thus,
\begin{equation}\label{nabla_eta}
(\nabla_{X}\eta_i)Y=0 \quad \Longleftrightarrow \quad (\nabla_{X}^{g}\eta_i)Y=\frac{1}{2}g(T(X, Y), \xi_i)=\frac{1}{2}T(X, Y, \xi_i)\,,
\end{equation}
for all $i$ and $X, Y\in\fr{X}(M)$.   Based on (\ref{gen_s})  and  on the fact that $T$  is a 3-form,   one can equivalently write
\[
(\nabla_{X}\eta_i)Y=0 \quad \Longleftrightarrow  \quad g(\nabla^{g}_{X}\xi_i, Y)=\frac{1}{2}T(X, Y, \xi_i)=-\frac{1}{2}T(X, \xi_i, Y)=-\frac{1}{2}g(T(X, \xi_i), Y)\,,
\]
that is, $\nabla\eta_i=0$ for all $i$ if and only if $\nabla^{g}_{X}\xi_i=-\frac{1}{2}T(X, \xi_i)$  for all $X\in\fr{X}(M)$ and all $i$, or equivalently $\nabla_{X}\xi_i=0$ for all $X\in\fr{X}(M)$, i.e., $\nabla\xi_i=0$. \\
(2)
We saw above that the condition $\nabla\xi_i=0$ for all $i\in\{1, \ldots, s\}$ is equivalent to the condition
\[
0=g(\nabla^{g}_{X}\xi_i, Z)+\frac{1}{2}T(X, \xi_{i}, Z)\,, \quad X, Z\in\fr{X}(M)\,.\qquad (\ast)
\]
As we require $T(X, Y, Z)=g(T(X, Y), Z)$ to be a 3-form for all $X, Y, Z\in\fr{X}(M)$, we should have 
\[
T(X, \xi_i, Z)+T(Z, \xi_i, X)=0
\]
 for all $X, Z\in\fr{X}(M)$ and $i=1, \ldots, s$. By $(\ast)$ this is equivalent to
the condition 
\[
g(\nabla^{g}_{X}\xi_i, Z)+g(\nabla^{g}_{Z}\xi_i, X)=0
\]
 for all $X, Z\in\fr{X}(M)$ and $i=1, \ldots, s$.  This  shows that each $\xi_i$ is a Killing vector field. \\
(3)  By (\ref{nabla_eta}) the condition $\nabla\eta_i=0$ for all $i$ is equivalent to  $(\nabla^{g}_{X}\eta_i)Y=\frac{1}{2}T(X, Y, \xi_i)$ for all $i$. 
Thus
\begin{align}
\label{dd_eta1}
\dd\eta_{i}(X, Y)=(\nabla^{g}_{X}\eta_i)Y-(\nabla^g_{Y}\eta_i)X&=\frac{1}{2}T(X, Y, \xi_i)-\frac{1}{2}T(Y, X, \xi_i)\\
&=T(X, Y, \xi_{i})=T(\xi_i, X, Y)\nonumber\,,
\end{align}
for all $i=1, \ldots, s$ and $X, Y\in\fr{X}(M)$. This proves that $\dd\eta_i=\xi_{i}\lrcorner T$  
and the relation $\dd\eta_i=2\nabla^g\eta_i$ follows   by $(\ast)$.\\
(4) This part  relies on the following standard identity
\begin{eqnarray*}\label{N_T}
N^{(1)}(X, Y)&=&
(\nabla_{\phi X}\phi)Y-(\nabla_{\phi Y}\phi)X+\phi((\nabla_{Y}\phi)X-(\nabla_{X}\phi)Y)\nonumber\\
&&+T(X, Y)-T(\phi X, \phi Y)+\phi\big(T(\phi X, Y)+T(X, \phi Y)\big)\,.\nonumber
\end{eqnarray*}
%whose proof  is direct and h avoid to present the proof. %whose proof is presented in the appendix. 
Taking the contraction via $g$ and using  (\ref{gphi}) this gives
\begin{align}
&N^{(1)}(X, Y, Z) =g(N^{(1)}(X, Y), Z)=\nonumber\\
&\quad=g((\nabla_{\phi X}\phi)Y, Z)-g((\nabla_{\phi Y}\phi)X, Z)-g((\nabla_{Y}\phi)X, \phi Z)+g((\nabla_{X}\phi)Y,  \phi Z)\nonumber\\
&\quad+T(X, Y, Z)-T(\phi X, \phi Y, Z)-T(\phi X, Y, \phi Z)-T(X, \phi Y, \phi Z)\nonumber\\
&\quad= g((\nabla_{\phi X}\phi)Y, Z)-g((\nabla_{\phi Y}\phi)X, Z)-g((\nabla_{Y}\phi)X, \phi Z)+g((\nabla_{X}\phi)Y,  \phi Z)-T^{-}(X, Y, Z)\,.\nonumber  
\end{align}
Since  we assume that $\nabla$ preserves $\phi$, i.e., $\nabla\phi=0$, each of the first four terms in this relation vanish, and this gives the result.\\
(5)  
This claim follows from the condition $\nabla\phi=0$, which can be equivalently written
as
\[
(\nabla^{g}_{X}\phi)Y=\frac{1}{2}\{-T(X, \phi Y)+\phi(T(X, Y))\}\,,\quad X, Y\in\fr{X}(M)\,. 
\]
(6) The  Levi-Civita connection $\nabla^g$ is torsion free, hence   
\[
\dd F(X, Y, Z)=(\nabla^{g}_{X}F)(Y, Z)-(\nabla^{g}_{Y}F)(X, Z)+(\nabla^g_{Z}F)(X, Y)\,,\quad X, Y, Z\in\fr{X}(M)\,.\quad (\sharp)
\]
Moreover,  by (\ref{cova_F}) we have
\[
(\nabla^g_{X}F)(Y, Z)=g((\nabla^g_{X}\phi)Z, Y)=-g((\nabla^{g}_{X}\phi)Y, Z)\,, \quad X, Y, Z\in\fr{X}(M)\,.
\]
A  combination of this relation with the result derived  above in part (5), yields the condition
\[
(\nabla^{g}_{X}F)(Y, Z)=-g((\nabla^{g}_{X}\phi)Y, Z)=\frac{1}{2}\{T(X, Y, \phi Z)+T(X, \phi Y, Z)\}\,,\quad X, Y, Z\in\fr{X}(M)\,.
\]
This identity  allows us to replace  any of the three components   in the right-hand-side of the relation in $(\sharp)$ by expressions depending on $T$, and this gives 
\begin{eqnarray*}
\dd F(X, Y, Z)&=&\frac{1}{2}\{T(X, Y, \phi Z)+T(X, \phi Y, Z)-T(Y, X, \phi Z)-T(Y, \phi X, Z)\\
&&+T(Z, X, \phi Y)+T(Z, \phi X, Y)\}\,.
\end{eqnarray*}
It is now  straightforward to see that  in this expression  the right-hand side simplifies to the cyclic sum  
\[
\fr{S}_{X, Y, Z}T(X, Y, \phi Z)=T(X, Y, \phi Z)+T(Y, Z, \phi X)+T(Z, X, \phi Y)\,.
\]
(7) By the result in (4) we know that $N^{(1)}(X, Y, Z)=-T^{-}(X, Y, Z)$ for all $X, Y, Z\in\fr{X}(M)$.  In this relation replace  $Z$ by $\phi Z$ to get
\[
-N^{(1)}(X, Y, \phi Z)=T(X, \phi Y, \phi^2Z)+T(\phi X, Y, \phi^2 Z)+T(\phi X, \phi Y, \phi Z)-T(X, Y, \phi Z)\,.\qquad (\flat)
\]
Next,  based on the  relation $\dd\eta_i=\xi_i\lrcorner T$ obtained above, we see that
\begin{eqnarray*}
T(X, \phi Y, \phi^2Z)&=&-T(X, \phi Y, Z)+\sum_{i=1}^{s}\eta_{i}(Z)T(X, \phi Y, \xi_i)=-T(X, \phi Y, Z)+\sum_{i=1}^{s}\eta_{i}(Z)\dd\eta_{i}(X, \phi Y)\,,\\
T(\phi X, Y, \phi^2 Z)&=&-T(\phi X, Y,  Z)+\sum_{i=1}^{s}\eta_{i}(Z)T(\phi X, Y, \xi_i)=-T(\phi X, Y,  Z)+\sum_{i=1}^{s}\eta_{i}(Z)\dd\eta_{i}(\phi X, Y)\,.
\end{eqnarray*}
Returning these expressions in $(\flat)$ and by the relation $N^{(2)}_{i}(X, Y):=\dd\eta_{i}(\phi X, Y)+\dd\eta_{i}(X, \phi Y)$,  we obtain
\[
 -N^{(1)}(X, Y, \phi Z)=-\fr{S}_{X, Y, Z}T(X, Y, \phi Z)+T(\phi X, \phi Y, \phi Z)+\sum_{i=1}^{s}\eta_{i}(Z)  N_{i}^{(2)}(X, Y)\,.
 \]
Our claim then follows by the relation $\dd F(X, Y, Z)=\fr{S}_{X, Y, Z}T(X, Y, \phi Z)$, proved in (6).
\pro

\subsubsection{The direct statement}
 Next we derive the explicit formula of the torsion 3-form $T$ of an adapted  connection $\nabla$, as presented above.  As we will see  below,  our   formula for $T$   reduces for $s=1$ to the torsion 3-form  appearing  in Friedrich-Ivanov's Theorem \ref{FrIvTHM}.
 Consequently, the result presented here provides a higher-dimensional generalization of the results in  \cite{FrIv} for almost contact metric   and almost Hermitian manifolds.
 Our approach is based on the relation (\ref{TxTyTz}), in which  we replace  $X, Y, Z$ by $\phi X, \phi  Y$ and $\phi Z$, respectively.
  Compared to the cases  $s=0, 1$ in \cite{FrIv}, this process is  more  involved and  demands numerous computations. We begin with the following crucial observation:
   \bl\label{T3phi}
 Let  $(M^{2n+s},  \phi, \xi_i, \eta_j, g)$  be a metric \negthinspace$f$-manifold endowed with the connection $\nabla=\nabla^g+\frac{T}{2}$ satisfying
 the relations in {\rm(\ref{parallel_objects})}, where $T$ is its torsion 3-form. 
 Then
  \begin{align*}
 &T(\phi^2X, \phi^2Y, \phi^2Z)=-T(X, Y, Z)+\sum_{i=1}^{s}\eta_{i}(X)\dd\eta_{i}(Y, Z)-\sum_{i=1}^{s}\eta_{i}(Y)\dd\eta_{i}(X, Z)+\sum_{i=1}^{s}\eta_{i}(Z)\dd\eta_{i}(X, Y)\\
&\quad\quad\quad-\sum_{i, j=1}^{s}\eta_{i}(X)\eta_{j}(Y)\dd\eta_{j}(Z, \xi_i)-\sum_{j, k=1}^{s}\eta_{j}(Y)\eta_{k}(Z)\dd\eta_{k}(X, \xi_j)+\sum_{i, k=1}^{s}\eta_{i}(X)\eta_{k}(Z)\dd\eta_{k}(Y, \xi_i)\\
&\quad\quad\quad +\sum_{1\leq i, j, k \leq s}\eta_{i}(X)\eta_{j}(Y)\eta_{k}(Z)\dd\eta_{k}(\xi_i, \xi_j)\,,
 \end{align*}
 for any $X, Y, Z\in\fr{X}(M)$.  For $s=0$ this gives the  obvious relation $T(\phi^2X, \phi^2Y, \phi^2Z)=-T(X, Y, Z)$, while for  $s=1$  it reduces to the identity $T(\phi^2X, \phi^2Y, \phi^2Z)=-T(X, Y, Z)+(\eta\wedge\dd\eta)(X, Y, Z)$. More in general, if $[\xi_i, \xi_j]=0$ for all $i, j\in\{1, \ldots, s\}$ then the following holds: 
 \[
 T(\phi^2X, \phi^2Y, \phi^2Z)=-T(X, Y, Z)+\sum_{i=1}^{s}(\eta_i\wedge d\eta_i)(X, Y, Z)\,,\quad X, Y, Z\in\fr{X}(M).
 \]
 \el
\pr
We simply write $\phi^2X=-X+\sum_{i}\eta_{i}(X)\xi_i$, $\phi^2Y=-Y+\sum_{j}\eta_{j}(X)\xi_j$  and $\phi^2Z=-Z+\sum_{k}\eta_{k}(Z)\xi_k$, and rely on the fact that $T$ is a 3-form. This gives
\begin{align*}
&T(\phi^2X, \phi^2Y, \phi^2Z)=T(-X+\sum_{i}\eta_{i}(X)\xi_i, -Y+\sum_{j}\eta_{j}(Y)\xi_j, -Z+\sum_{k}\eta_{k}(Z)\xi_k)\\
&=-T(X, Y, Z)+\sum_{k}\eta_{k}(Z)T(X, Y, \xi_k)-\sum_{j}\eta_{j}(Y)T(X, Z, \xi_j)+\sum_{i}\eta_{i}(X)T(Y, Z, \xi_i)\\
&-\sum_{i, j}\eta_{i}(X)\eta_{j}(Y)T(\xi_i, \xi_j, Z)-\sum_{i, k}\eta_{i}(X)\eta_{k}(Z)T(\xi_i, Y, \xi_k)-\sum_{j, k}\eta_{j}(Y)\eta_{k}(Z)T(X, \xi_j, \xi_k)\\
&+\sum_{i, j, k}\eta_{i}(X)\eta_{j}(Y)\eta_{k}(Z)T(\xi_i, \xi_j, \xi_k)\,.
\end{align*}
By repeatedly applying the relation (\ref{dd_eta1}),  namely  $\dd\eta_{i}(X, Y)=T(\xi_i, X, Y)=T(X, Y, \xi_i)$ for all $X, Y\in\fr{X}(M)$ and $i=1, \ldots, s$, we obtain the general formula given in the statement.
When $[\xi_i, \xi_j]=0$ for all $i, j\in\{1, \ldots, s\}$ by {\rm(\ref{parallel_objects})} we get $T(\xi_i, \xi_j)=\nabla_{\xi_i}\xi_j-\nabla_{\xi_j}\xi_i=0$.   Thus  $T(\xi_i, \xi_j, X)=g(T(\xi_i, \xi_j), X)=0$ for all $X\in\fr{X}(M)$  and hence,  under the assumption $[\xi_i, \xi_j]=0$ for all $i, j\in\{1, \ldots, s\}$,  we deduce that   in the stated expression the sums in the second and third line   vanish.  
\pro

\bt\label{Main_Thm_1}
  Let  $(M^{2n+s},  \phi, \xi_i, \eta_j, g)$  be a metric \negthinspace$f$-manifold  
 whose characteristic
vector fields  $\xi_i$ commute, i.e., $[\xi_i, \xi_j]=0$ for all $i, j\in\{1, \ldots, s\}$.  Suppose that  $\nabla=\nabla^{g}+\frac{1}{2}T $ is a metric connection with skew-torsion $T$  satisfying the conditions in
{\rm(\ref{parallel_objects})}.  Then, the torsion 3-form $T$ 
is given by
\begin{equation}\label{s_Torsion}
T=\sum_{i=1}^{s}\eta_{i}\wedge\dd\eta_i+\dd^{\phi}F+N^{(1)}-\sum_{i=1}^{s}(\eta_{i}\wedge(\xi_i\lrcorner N^{(1)}))\,.
\end{equation}
\et
\pr
By part (7) in Proposition \ref{characterize_torsion},     the torsion $T$ should satisfy the identity
\[
T(\phi X, \phi Y, \phi Z)=\dd F(X, Y, Z)-N^{(1)}(X, Y, \phi(Z))-\sum_{i=1}^{s}\eta_{i}(Z)N_{i}^{(2)}(X, Y)\,.
\]
Thus, recalling that $\eta_i\circ\phi=0$ for all $i\in\{1, \ldots, s\}$,  we obtain
\begin{eqnarray}
 T(\phi^2X, \phi^2Y, \phi^2Z)&=&\dd F(\phi X, \phi Y, \phi Z)-N^{(1)}(\phi X, \phi Y, \phi^2Z)-\sum_{i}\eta_{i}(\phi Z)N_{i}^{2}(\phi X, \phi Y)\nonumber\\
&=&-\dd^{\phi}F(X, Y, Z)+N^{(1)}(\phi X, \phi Y, Z)-\sum_{j}\eta_{j}(Z)N^{(1)}(\phi X, \phi Y, \xi_j)\,,\label{TT3phi2}
\end{eqnarray}
where  we replaced $\phi^2Z$ by $-Z+\sum_{j}\eta_{j}(Z)\xi_j$. Let us compute $N^{(1)}(\phi X, \phi Y, \xi_j)$.  Under our main assumption $[\xi_i, \xi_j]=0$ we will demonstrate that

\smallskip
{\bf Claim:}  $N^{(1)}(\phi X, \phi Y, \xi_j)=-N^{(1)}(X, Y, \xi_j)$ for all $j\in\{1, \ldots, s\}$.\\
{\it Proof:} When the condition $[\xi_i, \xi_j]=0$ holds for all $i, j\in\{1, \ldots, s\}$, by   (\ref{Nij_tensor}) we obtain that
\[
N^{(1)}(\phi X, \phi Y, Z)= -N^{(1)}(X, Y, Z) +\sum_{i=1}^{s}\eta_{i}(X)N^{(1)}(\xi_i, Y, Z)+\sum_{i=1}^{s}\eta_{i}(Y)N^{(1)}(X, \xi_{i}, Z)\,,\qquad (\star)
\]
for all $X, Y, Z\in\fr{X}(M)$.
Thus for any $j=1, \ldots, s$ the following holds:
\[
N^{(1)}(\phi X, \phi Y, \xi_j)=-N^{(1)}(X, Y, \xi_j)+\sum_{i=1}^{s}\eta_{i}(X)N^{(1)}(\xi_i, Y, \xi_j)+\sum_{i=1}^{s}\eta_{i}(Y)N^{(1)}(X, \xi_{i}, \xi_j)\,.
\]
Now, a  direct computation of $N^{(1)}(\xi_i, Y, \xi_j)$ and $N^{(1)}(X, \xi_{i}, \xi_j)$ shows  that these terms vanish 
  for all $i, j\in\{1, \ldots, s\}$.   For example, a combination of (\ref{general_N}) and  (\ref{gphi}) gives
\begin{eqnarray*}
N^{(1)}(\xi_i, Y, \xi_j)&=& -g(\phi[\xi_i, \phi Y], \xi_{k})-g([\xi_i, Y], \xi_j)+\sum_{k}\{\eta_{k}([\xi_i, Y])+\dd\eta_{k}(\xi_i, Y)\}\eta_{k}(\xi_j)\\
&=&-\eta_{j}([\xi_i, Y])+\eta_{j}([\xi_i, Y])+\dd\eta_{j}(\xi_i, Y)=\dd\eta_{j}(\xi_i, Y)=T(\xi_j, \xi_i, Y)\,,
\end{eqnarray*}
for all $Y\in\fr{X}(M)$, 
as $g(\phi[\xi_i, \phi Y], \xi_{k})=0$ and $\eta_{k}(\xi_j)=\delta_{kj}$ (hence in the sum over $k$ only the term with $k=j$ survives).
Above, we also applied the identity  $\dd\eta_{j}(X, Y)=T(\xi_j, X, Y)$, obtained in Proposition \ref{characterize_torsion}. 
This expression vanishes since $T(\xi_j, \xi_i, Y)=g(T(\xi_j, \xi_i), Y)=0$, for any $i, j\in\{1, \ldots, s\}$; this  relies on the condition  $[\xi_i, \xi_j]=0$ and the $\nabla$-parallelism of $\xi_i$,  for all $i$, as in the proof of Lemma \ref{T3phi} above. As an alternative, note that the relation $N^{(1)}(\xi_i, Y, \xi_j)=-N^{(1)}(Y, \xi_{i}, \xi_j)=0$ follows directly by  the identity in (\ref{general_N}), where  for the first equality above we use the fact that $N^{(1)}$ is totally skew-symmetric  (see (4) in Proposition \ref{characterize_torsion}). This proves our claim. \\
Therefore, we can now rewrite  (\ref{TT3phi2})
as follows:
\begin{eqnarray*}
 T(\phi^2X, \phi^2Y, \phi^2Z)
&=& -\dd^{\phi}F(X, Y, Z)+N^{(1)}(\phi X, \phi Y, Z)+\sum_{i=1}^{s}\eta_{i}(Z)N^{(1)}( X,  Y, \xi_i)\,.
 \end{eqnarray*}
Next,  we use $(\star)$ (or (\ref{Nij_tensor}))  to replace the term $N^{(1)}(\phi X, \phi Y, Z)$ inside this expression; this    gives 
 \begin{eqnarray*}
 T(\phi^2X, \phi^2Y, \phi^2Z)&=&-\dd^{\phi}F(X, Y, Z)-N^{(1)}(X, Y, Z)\\
 &&+\sum_{i=1}^{s}\eta_{i}(X)N^{(1)}(\xi_i, Y, Z)-\sum_{i=1}^{s}\eta_{i}(Y)N^{(1)}(\xi_{i} X, Z)+\sum_{i=1}^{s}\eta_{i}(Z)N^{(1)}( X,  Y, \xi_i)\\
&=& -\dd^{\phi}F(X, Y, Z)-N^{(1)}(X, Y, Z)+\sum_{i=1}^{s}(\eta_{i}\wedge(\xi_i\lrcorner N^{(1)}))(X, Y, Z)
 \end{eqnarray*}
   for all $X, Y, Z\in\fr{X}(M)$. Here we employed the fact that $N^{(1)}$ is totally skew-symmetric,  and that   for a 1-form $\eta$,  its dual vector field $\xi$ and a 3-form $N$ on a smooth manifold $M$, the 3-form $(\eta\wedge(\xi\lrcorner N))$ 
 has the expansion
 \[
 (\eta\wedge(\xi\lrcorner N))(X, Y, Z)=\eta(X)N(\xi, Y, Z)-\eta(Y)N(\xi, X, Z)+\eta(Z)N(\xi, X, Y)\,,\quad X, Y, Z\in\fr{X}(M)\,.
 \]
 On the other hand,  in  Lemma \ref{T3phi} we proved that 
  \[
  T(\phi^2X, \phi^2Y, \phi^2Z)=-T(X, Y, Z)+\sum_{i=1}^{s}(\eta_i\wedge d\eta_i)(X, Y, Z)\,.
  \]
  Comparing these two expressions of $T(\phi^2X, \phi^2Y, \phi^2Z)$  the stated formula for $T$ follows.   
\pro

\subsubsection{The converse statement}
Let us now   derive the ``converse'' of Theorem \ref{Main_Thm_1}, which will  finally enable  us to  establish a  higher-dimensional analogue  of the correspondence presented in \cite{FrIv}. 
In order to do so, we  first introduce the following form, which will be useful in the proof:
\[
 \dd F^{-}(X, Y, Z):=\dd F(X, \phi Y, \phi Z)+\dd F(\phi X, Y, \phi Z)+\dd F(\phi X, \phi Y, Z)-\dd F(X, Y, Z)\,,\  X, Y, Z\in\fr{X}(M)\,.
\]

\bl\label{ddF_minus} \textnormal{(see also \cite[p.~327]{FrIv} for $s=1$)}
Let $(M^{2n+s},  \phi, \xi_i, \eta_j, g)$ be a metric $f$-manifold. Then  
\[
\dd F^{-}(X, Y, Z)=-N^{(1)}(	X, Y, \phi Z)-N^{(1)}(Y, Z, \phi X)-N^{(1)}(Z, X, \phi Y)\,, \quad X, Y, Z\in\fr{X}(M)\,.
\]
\el
\pr
This is a  lengthy but direct computation based  on the formula 
\[
\dd F(X, Y, Z)=X F(Y, Z)+Y F(Z, X)+Z F(X, Y)-F([X, Y], Z)-F([Z, X], Y)-F([Y, Z], X)\,.
\]
This approach relies on the identities  in (\ref{phi_invariance}), that is, the $\phi$-invariance of $F$ and the skew-symmetry of $\phi$ with respect to $F$. 
The identities $\xi_i\lrcorner F=0$ $(i=1, \ldots, s)$ and 
\[
F(N^{(1)}(X, Y), Z)=g(N^{(1)}(X, Y), \phi Z)=N^{(1)}(X, Y, \phi Z)\,.
\]
are also necessary.  We leave the remaining details to the reader.
\pro

We also need the following lemma:
\bl\label{formulas_N12}
Let  $(M^{2n+s},  \phi, \xi_i, \eta_j, g)$  be a metric $f$-manifold  
 whose characteristic
vector fields  $\xi_i$ commute, i.e., $[\xi_i, \xi_j]=0$ for all $i, j\in\{1, \ldots, s\}$. Suppose that each $\xi_i$  is a Killing  vector field and 
that $N^{(1)}$ is totally skew-symmetric.  Then  the following identities hold:    $N^{(1)}(X, \xi_i,  \xi_j)=0$ for all $i, j\in\{1, \ldots, s\}$ and moreover
\begin{equation}\label{final_N12}
N^{(1)}(\phi X, Y, \xi_i)=N^{(1)}(X, \phi Y, \xi_i)=N^{(2)}_{i}(X, Y)=\dd F(X, Y, \xi_i)=-\dd F(\phi X, \phi Y, \xi_i)
\end{equation}
for any $X, Y\in\fr{X}(M)$ and $i=1, \ldots, s$.
\el

\pr
Note that for $s=1$ these identities  are already known from \cite{FrIv}; see also Remark \ref{differ/} at the end of the proof. 
We have already mentioned that the first relation $N^{(1)}(X, \xi_i,  \xi_j)=0$ for all $i, j\in\{1, \ldots, s\}$ and $X\in\fr{X}(M)$,  follows by a direct computation based on (\ref{general_N}). Relying on this fact,   by Lemma \ref{identity_N}  
we obtain  $N^{(1)}(X, Y, \xi_k)= -N^{(1)}(\phi X, \phi Y, \xi_k)$ for any $k=1, \ldots, s$ and $X, Y\in\fr{X}(M)$.  
As a consequence, we see that
\begin{eqnarray*}
N^{(1)}(\phi X, Y, \xi_k)&=&-N^{(1)}(\phi^2X, \phi Y, \xi_k)=N^{(1)}(X, \phi Y, \xi_k)-\sum_{i}\eta_{i}(X)N^{(1)}(\xi_i, \phi Y, \xi_k)\\
&=&N^{(1)}(X, \phi Y, \xi_k)\,.
\end{eqnarray*}
This proves the first  equality in (\ref{final_N12}). \\
\noindent Next we will show  that $\dd F(X, Y, \xi_i)=-\dd F(\phi X, \phi Y, \xi_i)$ for all $i=1, \ldots, s$ and  $X, Y\in\fr{X}(M)$.  Let us first mention that using  (\ref{cova_F}) we see that
\begin{eqnarray*}
\dd F(X, Y, Z)&=&(\nabla^g_{X}F)(Y, Z)-(\nabla^{g}_{Y}F)(X, Z)+(\nabla^g_{Z}F)(X, Y)\\
&=&-g((\nabla^g_{X}\phi)Y, Z)+g((\nabla^g_{Y}\phi)X, Z)-g((\nabla^g_{Z}\phi)X, Y)\,,\quad X, Y, Z\in\fr{X}(M)
\end{eqnarray*}
which implies that
\begin{equation}\label{difFxij}
\dd F(\xi_j, Y, Z) =-g((\nabla^{g}_{\xi_j}\phi)Y, Z)+g((\nabla^g_{Y}\phi)\xi_j, Z)-g((\nabla^g_{Z}\phi)\xi_j, Y)\,,
\end{equation}
for any $j=1, \ldots, s$. 
Next, a  combination of (\ref{cova_F}) and Proposition \ref{symmetries_FN} gives
\[
g((\nabla^{g}_{X}\phi)Y, Z)+g((\nabla^g_{X}\phi)\phi Y, \phi Z)=\sum_{i=1}^{s}\eta_{i}(Y)(\nabla^{g}_{X}\eta_i)\phi Z-\sum_{i=1}^{s}\eta_{i}(Z)(\nabla^g_{X}\eta_i)\phi Y\,, \quad X, Y, Z\in\fr{X}(M)\,.
\]
Let us   take the cyclic sum of this relation and put $X=\xi_j$. This gives
\begin{eqnarray}
&&g((\nabla^g_{\xi_j}\phi)Y, Z)+g((\nabla^g_{\xi_j}\phi)\phi Y, \phi Z)+g((\nabla^g_{Y}\phi)Z, \xi_j)+g((\nabla^{g}_{Z}\phi)\xi_j), Y)=\nonumber\\
&&\sum_{i}\eta_{i}(Y)(\nabla^g_{\xi_j}\eta_i)\phi Z-\sum_{i}\eta_{i}(Z)(\nabla^g_{\xi_{j}}\eta_i)\phi Y
-\sum_{i}\eta_{i}(\xi_j)(\nabla^{g}_{Y}\eta_i)\phi Z+\sum_{i}\eta_{i}(\xi_j)(\nabla^g_{Z}\eta_i)\phi Y\,.\nonumber\\ \label{update_rel}
\end{eqnarray}
However,  by part (3) in Proposition \ref{prop_derivatives}  we know that 
$0=\dd\eta_i(\xi_j, Y)=2(\nabla^{g}_{\xi_j}\eta_i)Y$ for any $Y\in\fr{X}(M)$, so in (\ref{update_rel})   the first two sums vanish.
 Using the middle equality  in (\ref{cova_F}), the identity $\eta_{i}(\xi_j)=\delta_{ij}$ and the relation (\ref{difFxij})  mentioned above, we can  thus rewrite (\ref{update_rel}) as 
\[
-\dd F(\xi_j, Y, Z)+g((\nabla^g_{\xi_j}\phi)\phi Y, \phi Z)=-(\nabla^g_{Y}\eta_j)\phi Z+(\nabla^{g}_{Z}\eta_j)\phi(Y)\,,\quad (\dag)
\]
 for all $j=1, \ldots, s$ and $Y, Z\in\fr{X}(M)$. On the other hand, an application of (\ref{difFxij}) gives
 \[
 \dd F(\xi_j, \phi Y, \phi Z)=-g((\nabla^g_{\xi_j}\phi)\phi Y, \phi Z)+g((\nabla^g_{\phi Y}\phi)\xi_j, \phi Z)-g((\nabla^g_{\phi Z}\phi)\xi_j, \phi Y)
 \]
 which we can rewrite as
 \begin{eqnarray*}
 -\dd F(\xi_j, \phi Y, \phi Z)-g((\nabla^g_{\xi_j}\phi)\phi Y, \phi Z)&=&-g((\nabla^g_{\phi Y}\phi)\xi_j, \phi Z)+g((\nabla^g_{\phi Z}\phi)\xi_j, \phi Y)\\
 &\overset{(\ref{cova_F})}{=}&g((\nabla^{g}_{\phi Y}\phi)\phi Z, \xi_j)-g((\nabla^g_{\phi Z}\phi)\phi Y, \xi_j)\,.
 \end{eqnarray*}
 We will now use the identity 
 \begin{equation}\label{g_nabla_phi}
 g((\nabla^{g}_{X}\phi)\phi Y, \xi_j)=(\nabla^g_{X}\eta_j)Y\,,\quad j=1, \ldots, s
 \end{equation}
which we mentioned after the relation (\ref{gen_s2}),  as $(\nabla^g_{X}F)(\xi_j, \phi Y)=(\nabla^g_{X}\eta_j)Y$; this holds since the set $\{\xi_1, \ldots, \xi_s\}$ consists of commuting   Killing vector fields and %  as we saw in Proposition \ref{prop_derivatives},
 in this scenario   we have
$(\mc{L}_{\xi_j}\eta_{i})(X)=\dd\eta_i(\xi_j, X)=0$, for all $X\in\fr{X}(M)$. 
Thus, for all $j=1, \ldots s$ and $Y, Z\in\fr{X}(M)$  the following relation holds:
\[
-\dd F(\xi_j, \phi Y, \phi Z)-g((\nabla^g_{\xi_j}\phi)\phi Y, \phi Z)=(\nabla^{g}_{\phi Y}\eta_j)Z-(\nabla^g_{\phi Z}\eta_j)Y\,. \quad (\ddag)
\]
Adding $(\dag)$ and $(\ddag)$ we finally obtain
 \begin{eqnarray*}
-\dd F(\xi_j, Y, Z) -\dd F(\xi_j, \phi Y, \phi Z)&=&-(\nabla^g_{Y}\eta_j)\phi Z+(\nabla^{g}_{Z}\eta_j)\phi(Y)
+(\nabla^{g}_{\phi Y}\eta_j)Z-(\nabla^g_{\phi Z}\eta_j)Y\\
&=&-\frac{1}{2}\dd\eta_j(Y, \phi Z)+\frac{1}{2}\dd\eta_j(Z, \phi Y)+\frac{1}{2}\dd\eta_{j}(\phi Y, Z)-\frac{1}{2}\dd\eta_{j}(\phi Z, Y)\\
&=&0\,,
 \end{eqnarray*}
 for all $j\in\{1, \ldots, s\}$. This proves the claim. \\
 Next, the relation $N^{(1)}(\xi_j, \phi X, Y)=N^{(2)}_{j}(X, Y)$ for all $X, Y\in\fr{X}(M)$ and $j=1, \ldots, s$ follows from the identity 
\begin{equation}\label{NXYxi}
N^{(1)}(X, Y, \xi_j)=N^{(1)}(\xi_j, X, Y)=\dd\eta_{j}(X, Y)-\dd\eta_{j}(\phi X, \phi Y)\,.
\end{equation}
Indeed, using this  we get
\begin{eqnarray*}
N^{(1)}(\xi_j, \phi X, Y)&=&\dd\eta_j(\phi X, Y)-\dd\eta_j(\phi^2X, \phi Y)\\
&=&\dd\eta_j(\phi X, Y)+\dd\eta_{j}(X, \phi Y)-\sum_{k}\eta_k(X)\dd\eta_j(\xi_k, \phi Y)=N^{(2)}_{j}(X, Y)
\end{eqnarray*}
 for any $j=1, \ldots, s$, as $\xi_k\lrcorner \dd\eta_j=0$ for all $k, j\in\{1, \ldots, s\}$. We will prove  (\ref{NXYxi}) a few below. \\
Let us first show that $\dd F(\xi_j, Y, Z)=N^{(2)}_{j}(Y, Z)=\dd\eta_j(\phi Y, Z)+\dd\eta_j(Y, \phi Z)$, as this relies on the  identity 
 $N^{(1)}(\xi_j, \phi Y, Z)=N^{(2)}_{j}(Y, Z)$. Indeed,   by Lemma \ref{ddF_minus}
 one can write
 \begin{eqnarray*}
 \dd F(X, Y, Z)&=& \dd F(X, \phi Y, \phi Z)+\dd F(\phi X, Y, \phi Z)+\dd F(\phi X, \phi Y, Z)\\
 &&+N^{(1)}(X, Y, \phi Z)+N^{(1)}(Y, Z, \phi X)+N^{(1)}(Z, X, \phi Y)\,, \quad X, Y, Z\in\fr{X}(M)\,.
 \end{eqnarray*}
 Since $N^{(1)}$ is assumed to be a 3-form, for $X=\xi_j$ the above formula gives
 \begin{eqnarray*}
 \dd F(\xi_j, Y, Z)&=&\dd F(\xi_j, \phi Y, \phi Z)+N^{(1)}(\xi_j, Y, \phi Z)+N^{(1)}(Z, \xi_j, \phi Y)\\
 &=&\dd F(\xi_j, \phi Y, \phi Z)+N^{(1)}(Y, \phi Z, \xi_j)+N^{(1)}(\phi Y, Z, \xi_j)\\
 &=&\dd F(\xi_j, \phi Y, \phi Z)+2N^{(1)}(\xi_j, \phi Y, Z)=\dd F(\xi_j, \phi Y, \phi Z)+2N^{(2)}_{j}(Y, Z)\,,
 \end{eqnarray*}
 where it was used  the identity $N^{(1)}(\phi Y, Z, \xi_j)=N^{(1)}(Y, \phi Z, \xi_j)$  for all $j$.
The claim follows by replacing $\dd F(\xi_j, \phi Y, \phi Z)$ by $-\dd F(\xi_j, Y, Z)$, as it was proved above.

\noindent It remains to prove the identity (\ref{NXYxi}).
Recall that by part (4) of Propositon \ref{prop_derivatives}  
we have $g(\nabla^{g}_{X}\xi_i, \xi_j)=0$
for all $X\in\fr{X}(M)$ and $i, j\in\{1, \ldots, s\}$. Thus,  an application of the identity in (\ref{tensor_N_T}) gives
\begin{eqnarray*}
 N^{(1)}(X, Y, \xi_j)&=&g((\nabla^{g}_{\phi X}\phi)Y, \xi_j)-g((\nabla^{g}_{\phi Y}\phi)X, \xi_j)+g((\nabla^g_{X}\phi)\phi Y, \xi_j)-g(\nabla^{g}_{Y}\phi)\phi X, \xi_j)\\
&&+\sum_{i=1}^{s}\eta_{i}(X)g(\nabla^{g}_{Y}\xi_{i}, \xi_j)-\sum_{i=1}^{s}\eta_{i}(Y)g(\nabla^{g}_{X}\xi_{i}, \xi_j) \\
&=&g((\nabla^{g}_{\phi X}\phi)Y, \xi_j)-g((\nabla^{g}_{\phi Y}\phi)X, \xi_j)+g((\nabla^g_{X}\phi)\phi Y, \xi_j)-g(\nabla^{g}_{Y}\phi)\phi X, \xi_j)\\
&\overset{(\ref{g_nabla_phi})}{=}&g((\nabla^{g}_{\phi X}\phi)Y, \xi_j)-g((\nabla^{g}_{\phi Y}\phi)X, \xi_j)+(\nabla^g_{X}\eta_j) Y -  (\nabla^{g}_{Y}\eta_j) X\\
&\overset{(\ref{cova_F})}{=}&-g((\nabla^{g}_{\phi X}\phi)\xi_j, Y)+g((\nabla^{g}_{\phi Y}\phi)\xi_j, X)+(\nabla^g_{X}\eta_j) Y -  (\nabla^{g}_{Y}\eta_j) X\\
&=&g(\phi(\nabla^{g}_{\phi X}\xi_j), Y)-g(\phi(\nabla^{g}_{\phi Y}\xi_j), X))+(\nabla^g_{X}\eta_j) Y -  (\nabla^{g}_{Y}\eta_j) X\\
&\overset{(\ref{gphi})}{=}&-g(\nabla^{g}_{\phi X}\xi_j, \phi Y)+g(\nabla^{g}_{\phi Y}\xi_j, \phi X)+(\nabla^g_{X}\eta_j) Y -  (\nabla^{g}_{Y}\eta_j) X\\
&\overset{(\ref{gen_s})}{=}&-(\nabla^{g}_{\phi X}\eta_j)\phi Y+(\nabla^{g}_{\phi Y}\eta_{j})\phi X+(\nabla^g_{X}\eta_j) Y -  (\nabla^{g}_{Y}\eta_j) X\\
&=&-\dd\eta_j(\phi X, \phi Y)+\dd\eta_j(X, Y)\,,
\end{eqnarray*}
for all $X, Y\in\fr{X}(M)$ and $j=1, \ldots, s$. This completes the proof.
  \pro
  
\br\label{differ/}
The  identities  in (\ref{final_N12}) provide  a higher-dimensional analogue  of those given in  \cite[Lemma~8.3]{FrIv} for $s=1$.
We note, however,  that in that setting the condition $g(\nabla^{g}_{X}\xi, \xi)=0$  holds without assuming that $\xi$ is Killing. For $s\geq 2$,  by contrast, there are terms of the form $g(\nabla^{g}_{X}\xi_i, \xi_j)$ for $i\neq j$, hence in the proof there are some remarkable differences.  For instance, for $s=1$ is not necessary to   assume that $\xi$ is Killing to obtain the relation $ N^{(1)}(X, Y, \xi)=-\dd\eta(\phi X, \phi Y)+\dd\eta(X, Y)$. 
\er

\bt\label{mainThm2}
Let  $(M^{2n+s},  \phi, \xi_i, \eta_j, g)$  be a metric $f$-manifold  
 whose characteristic
vector fields  $\xi_i$ commute, i.e., $[\xi_i, \xi_j]=0$ for all $i, j\in\{1, \ldots, s\}$. Suppose that each $\xi_i$  is a Killing  vector field and 
that $N^{(1)}$ is totally skew-symmetric. Then $(M^{2n+s},  \phi, \xi_i, \eta_j, g)$  admits a unique metric connection $\nabla$ with skew-torsion
$T$ preserving the structure,  that is, satisfying the conditions in {\rm(\ref{parallel_objects})}. The torsion 3-form $T$ is given by {\rm(\ref{s_Torsion})} and satisfies 
$T(\xi_i, \xi_j)=0$ for all $i, j\in\{1, \ldots, s\}$  and  $\xi_j\lrcorner T=\dd\eta_j=2\nabla^g\eta_j$ for all $j=1, \ldots, s$.
\et
\pr
Let us define the connection $\nabla$ by the formula $\nabla=\nabla^g+\frac{1}{2}T$, where $T$ is given by {\rm(\ref{s_Torsion})}.
Under the assumption that $N^{(1)}$ is totally skew-symmetric, it is obvious that $T$ is a 3-form on $M^{2n+s}$. 
Regarding the relation $\xi_j\lrcorner T=\dd\eta_j=2\nabla^g\eta_j$ for all $j=1, \ldots, s$,  the second equality was obtained in Proposition \ref{prop_derivatives}. The first equality  follows by a direct computation,  using the explicit form   {\rm(\ref{s_Torsion})}  of $T$ and relying on the  fact   
that the vector fields $\xi_1, \ldots, \xi_s$ are commuting  Killing vector fields (hence, by Proposition \ref{prop_derivatives} we also have   $\xi_j\lrcorner\dd\eta_i=0$ for all $i, j\in\{1, \ldots, s\}$).

  Now, as $T$ is a 3-form  we  automatically  have $\nabla g=0$. It remains to verify the remaining conditions $\nabla\xi_j=\nabla\eta_j=0=\nabla\phi$.
 The $\nabla$-parallelism of $\xi_j$, expressed by  $g(\nabla^{g}_{X}\xi_j, Y)=\frac{1}{2}T(\xi_j, X, Y)$  for all $j=1, \ldots, s$,  is a direct consequence of the fact that each $\xi_j$ is Killing. 
 Since for all  $j=1, \ldots, s$,  the conditions $\nabla{\xi_j}=0$ and $\nabla\eta_j=0$ are equivalent (see  Proposition \ref{characterize_torsion}),  
  it remains to prove that $\nabla\phi=0$. 
 Set $T^{\phi}(X, Y, Z):=T(X, Y, \phi Z)+T(X, \phi Y, Z)$ 
 for any $X, Y, Z\in\fr{X}(M)$.   
It turns out that the condition $(\nabla_{X}\phi)Y=0$ for all $X, Y\in\fr{X}(M)$,  %or equivalently  $g((\nabla_{X}\phi)Y, Z)=0$ for all $X, Y, Z\in\fr{X}(M)$, 
is equivalent to the identity 
\begin{equation}\label{main_Tphi}
-2g((\nabla^g_{X}\phi)Y, Z)=T^{\phi}(X, Y, Z)\,,\quad X, Y, Z \in\fr{X}(M)\,,
\end{equation}
 see    (5) in Proposition \ref{characterize_torsion}.   Using the definition of $T$, see (\ref{s_Torsion}),  we will show that this identity holds.
  We will use Lemmas \ref{basic_1}, \ref{ddF_minus} and  \ref{formulas_N12}.  First we see that
\begin{eqnarray*}
T(X, Y, \phi Z)&=&\sum_{i}(\eta_i\wedge\dd\eta_i)(X, Y, \phi Z)+\dd^{\phi}F(X, Y, \phi Z)+N^{(1)}(X, Y, \phi Z)\\
&&-\sum_{i}(\eta_i\wedge(\xi_i\lrcorner N^{(1)}))(X, Y, \phi Z)\\
&=&\sum_{i}\left\{\eta_{i}(X)\dd\eta_{i}(Y, \phi Z)-\eta_{i}(Y)\dd\eta_{i}(X, \phi Z)\right\}-\dd F(\phi X, \phi Y, \phi^2Z)+N^{(1)}(X, Y, \phi Z)\\
&&-\sum_{i}\{\eta_{i}(X)N^{(1)}(\xi_i, Y, \phi Z)-\eta_i(Y)N^{(1)}(\xi_i, X, \phi Z)\}\\
&=&\sum_{i}\left\{\eta_{i}(X)\dd\eta_{i}(Y, \phi Z)-\eta_{i}(Y)\dd\eta_{i}(X, \phi Z)\right\}+\dd F(\phi X, \phi Y, Z)-\sum_{i}\eta_{i}(Z)\dd F(\phi X, \phi Y, \xi_i)\\
&&+N^{(1)}(X, Y, \phi Z)-\sum_{i}\eta_{i}(X)N^{(1)}(\xi_i, Y, \phi Z)+\sum_i\eta_i(Y)N^{(1)}(\xi_i, X, \phi Z)\,,
\end{eqnarray*}
where we used that $\eta_i\circ\phi=0$ for all $i=1, \ldots, s$ and that $\phi^2Z=-Z+\sum_{i}\eta_{i}(Z)\xi_i$.
%In a similar way  we compute %5one can compute $T(X, \phi Y, Z)$ %Similarly, we compute the
Similarly, 
\begin{eqnarray*}
T(X, \phi Y, Z)&=&\sum_{i}\{\eta_{i}(X)\dd\eta_{i}(\phi Y, Z)+\eta_{i}(Z)\dd\eta_{i}(X, \phi Y)\}+\dd F(\phi X, Y, \phi Z)-\sum_{i}\eta_{i}(Y)\dd F(\phi X, \xi_i, \phi Z)\\
&&+N^{(1)}(X, \phi Y, Z)-\sum_{i}\eta_{i}(X)N^{(1)}(\xi_i, \phi Y, Z)-\sum_{i}\eta_{i}(Z)N^{(1)}(\xi_i, X, \phi Y)\,.
\end{eqnarray*}
Thus, %we deduce that % for the quantity $T^{\phi}$  deduce that
\begin{eqnarray}
T^{\phi}(X, Y, Z) 
&=&\sum_{i}\eta_{i}(X)\{\dd\eta_{i}(Y, \phi Z)+\dd\eta_i(\phi Y, Z)\}-\sum_{i}\eta_{i}(Y)\dd\eta_{i}(X, \phi Z)+\sum_{i}\eta_{i}(Z)\dd\eta_{i}(X, \phi Y)\nonumber\\
&&+\dd F(\phi X, \phi Y, Z)+\dd F(\phi X, Y, \phi Z)+N^{(1)}(X, Y, \phi Z)+N^{(1)}(X, \phi Y, Z)\nonumber\\
&&-\sum_{i}\eta_{i}(X)\{N^{(1)}(\xi_i, Y, \phi Z)+N^{(1)}(\xi_i, \phi Y, Z)\}\nonumber\\
&& +\sum_i\eta_i(Y)\{N^{(1)}(\xi_i, X, \phi Z)-\dd F(\phi X, \xi_i, \phi Z)\}\nonumber\\
&&-\sum_{i}\eta_{i}(Z)\{N^{(1)}(\xi_i, X, \phi Y)+\dd F(\phi X, \phi Y, \xi_i)\}\,,\label{Tphi}
\end{eqnarray}
for any $X, Y, Z\in\fr{X}(M)$.  Now, by Lemma \ref{ddF_minus} we see that
\begin{eqnarray*}
\dd F(\phi X, \phi Y, Z)+\dd F(\phi X, Y, \phi Z)&=&-N^{(1)}(	X, Y, \phi Z)-N^{(1)}(Y, Z, \phi X)-N^{(1)}(Z, X, \phi Y)\\
&&+\dd F(X, Y, Z)-\dd F(X, \phi Y, \phi Z)\,,
\end{eqnarray*}
thus the second line in  (\ref{Tphi}) equals  to $-N^{(1)}(Y, Z, \phi X)+\dd F(X, Y, Z)-\dd F(X, \phi Y, \phi Z)$. 
This gives 
\begin{align}
&T^{\phi}(X, Y, Z)=\sum_{i}\eta_{i}(X)N^{(2)}_{i}(Y, Z) -\sum_{i}\eta_{i}(Y)\dd\eta_{i}(X, \phi Z)+\sum_{i}\eta_{i}(Z)\dd\eta_{i}(X, \phi Y)\nonumber\\
&-N^{(1)}(Y, Z, \phi X)+\dd F(X, Y, Z)-\dd F(X, \phi Y, \phi Z)
-\sum_{i}\eta_{i}(X)\{N^{(1)}(\xi_i, Y, \phi Z)+N^{(1)}(\xi_i, \phi Y, Z)\}\nonumber\\
& +\sum_i\eta_i(Y)\{N^{(1)}(\xi_i, X, \phi Z)-\dd F(\phi X, \xi_i, \phi Z)\}
-\sum_{i}\eta_{i}(Z)\{N^{(1)}(\xi_i, X, \phi Y)+\dd F(\phi X, \phi Y, \xi_i)\}\,,\label{Tphi2}
\end{align}
where we have applied the definition of $N_{i}^{(2)}$. The condition in (\ref{main_Tphi}) now follows directly from Lemma \ref{basic_1} and   the
identities (\ref{final_N12})  in Lemma \ref{formulas_N12}. % the first    three lines  of (\ref{Tphi2})  combine to give $-2g((\nabla^g_{X}\phi)Y, Z)$, 
%and the remaining two lines vanish.
\pro
We are now in a position to encode the correspondence established through Proposition \ref{characterize_torsion} 
and the Theorems \ref{Main_Thm_1} and \ref{mainThm2}, as follows  %(recall that for $s=0, 1$ this correspondence was  derived  in \cite{FrIv}).
\bc\label{mainCorol}
Let  $(M^{2n+s},  \phi, \xi_i, \eta_j, g)$  be a metric $f$-manifold  
 whose characteristic
vector fields  $\xi_i$ commute, i.e., $[\xi_i, \xi_j]=0$ for all $i, j\in\{1, \ldots, s\}$. Then the following conditions are equivalent: \\
1) $N^{(1)}$ is totally skew-symmetric and  $\xi_1, \ldots, \xi_s$ are Killing   vector fields.\\ 
2)    $(M^{2n+s},  \phi, \xi_i, \eta_j, g)$  admits an almost metric $\phi$-connection with totally skew-symmetric torsion tensor 
$T$.  In particular, this connection $\nabla$ is unique and  determined by the relation $\nabla=\nabla^g+\frac{T}{2}$, where the torsion 3-form $T$   is defined by {\rm(\ref{s_Torsion})}. 
\ec

%%%%%%%%%%%%%%%%%%%%%%%%%%%%%%%%%
%%%%%%%%%Special_cases%%%%%%%%%%%%
%%%%%%%%%%%%%%%%%%%%%%%%%%%%%%%%%%%

\section{Special cases}\label{Section3}
In the sequel we discuss the previous results for some special types of  metric $f$-manifolds.

\subsection{Normal metric $f$-manifolds with skew-torsion}
Let  $(M^{2n+s},  \phi, \xi_i, \eta_j, g)$  be a metric $f$-manifold which is normal, i.e., $N^{(1)}=0$. Recall that in this case we have $[\xi_i, \xi_j]=0=\mc{L}_{\xi_i}\eta_j$    for all $1\leq i, j\leq s$ and  $N^{(2)}_{i}=0$ for all $i=1, \ldots, s$.  By the previous corollary it is thus direct that 

\bp\label{normal_skewtorsion}
A  metric $f$-manifold  $(M^{2n+s},  \phi, \xi_i, \eta_j, g)$  which is normal admits a unique almost metric $\phi$-connection $\nabla$ with skew-torsion
if and only if each $\xi_i$ is a Killing vector field. In this case the torsion 3-form $T$ of $\nabla$ is given by
\[
T=\sum_{i=1}^{s}\eta_{i}\wedge\dd\eta_i+\dd^{\phi} F\,.
\]
\ep
Note that for $s=1$ this proposition reduces to  part 2) of Theorem 8.4 in  \cite{FrIv}.

\subsection{$\mc{K}$-manifolds with skew-torsion}
Recall that a metric $f$-manifold $(M^{2n+s},  \phi, \xi_i, \eta_j, g)$   which is normal and has closed fundamental 2-form $F$ is called a $\mc{K}$-manifold.   For $\mc{K}$-manifolds with $s\geq 1$, by  Corollary \ref{mainCorol} we deduce that 
\bp\label{K1}
Any $\mc{K}$-manifold $(M^{2n+s}, \phi, \xi_i, \eta_j, g)$ with $s\geq 1$ admits a unique almost  metric $\phi$-connection $\nabla$ with skew-torsion $T$. This connection is   given by  $\nabla=\nabla^{g}+\frac{1}{2}\sum_{j=1}^{s}\eta_{j}\wedge\dd\eta_j$.
\ep
For $s=1$ Proposition \ref{K1} reduces to the classical result that   any quasi-Sasakian manifold $(M^{2n+1}, \phi, \xi, \eta, g)$ admits a metric connection with skew-torsion given by $T=\eta\wedge\dd\eta$, see \cite{FrIv}. Note also that the assumption $s\geq 1$ is necessary to  exclude the integrable K\"ahler case.

Next we will show that    an almost $\mc{S}$-manifold $(M^{2n+s}, \phi, \xi_j, \eta_i, g)$  admits  an almost metric $\phi$-connection with skew-torsion  if and only if $M$ is normal, i.e., an $\mc{S}$-manifold. Moreover, in this case we obtain the remarkable relation $\nabla T=0$.  This  result for $s> 1$  furnishes us with a higher-dimensional analogue of   part 1) in \cite[Theorem 8.4]{FrIv}.

\bt\label{S-mnfds1}
Let  $(M^{2n+s},  \phi, \xi_i, \eta_j, g)$ be a   contact  metric $f$-manifold. Then $M^{2n+s}$ admits an almost  metric $\phi$-connection $\nabla$ with skew-torsion $T$ if and only if $M^{2n+s}$ is an $\mc{S}$-manifold. If this is the case, the uniquely determined   connection $\nabla$  has $\nabla$-parallel torsion 3-form $T$  given by 
\[
T=\sum_{i=1}^{s}\eta_{i}\wedge\dd\eta_i=2\bar{\eta}\wedge F\,,\quad\text{where}\quad \bar{\eta}=:\sum_{i=1}^{s}\eta_i\,.
\]
\et
\pr
Recall that a contact metric $f$-manifold is a metric $f$-manifold  $(M^{2n+s},  \phi, \xi_i, \eta_j, g)$  with 
$2 F=\dd\eta_1=\dd\eta_2=\ldots=\dd\eta_s$. It is known that  in this case we have $[\xi_i, \xi_j]=0$ for all $i, j$, see \cite[Corollary 2.4]{CFF90}. 
Suppose that $M^{2n+s}$ admits such a connection $\nabla$. Then, by Corollary \ref{mainCorol} the vector fields $\xi_i$ for $i\in\{1, \ldots, s\}$  are all  Killing,  and $N^{(1)}$ is totally skew-symmetric.  
We need to show that $N^{(1)}$ vanishes identically. 
Since the assumptions of Lemma  \ref{formulas_N12} are  satisfied, we may apply it and use   the relation   $\xi_j\lrcorner \dd F=N^{(2)}_{j}$ for all $j\in\{1, \ldots, s\}$. 
Since $F$ is closed, we thus obtain  $N^{(2)}_{j}=0$ for all $j$. This also implies that $\xi_j\lrcorner N^{(1)}=0$ for all $j\in\{1, \ldots, s\}$.   Indeed, recall by  (\ref{final_N12})  that 
 $N^{(1)}(\phi X, Y, \xi_j)=N^{2}_{j}(X, Y)$  for any $X, Y\in\fr{X}(M)$, which gives $0=N^{(1)}(\phi X, Y, \xi_j)$.  Substituting 
  $X$ by  $\phi X$ in this identity,  yields
\[
0=-N^{(1)}(X, Y, \xi_j)+\sum_{i}\eta_{i}(X)N^{(1)}(\xi_i, Y, \xi_j)=-N^{(1)}(X, Y, \xi_j)
\]
 for any $X, Y\in\fr{X}(M)$ and $j\in\{1, \ldots, s\}$. 
This proves that  $\xi_j\lrcorner N^{(1)}=0$ for all $j\in\{1, \ldots, s\}$. 

Now, since  $\dd F=0$ by Lemma \ref{ddF_minus} we also get
$N^{(1)}(X, Y, \phi Z)+N^{(1)}(Y, Z, \phi X)+N^{(1)}(Z, X, \phi Y)=0$.  Relying then on Lemma \ref{identity_N} it is easy to see that 
\[
N^{(1)}(X, Y, \phi Z)=N^{(1)}(Y, Z, \phi X)=N^{(1)}(Z, X, \phi Y)=-N^{(1)}(\phi X, \phi Y, \phi Z)\,.
\]
It follows that $N^{(1)}(\phi X, Y, Z)=0$ for any $X, Y, Z\in\fr{X}(M)$ and this implies that  $N^{(1)}(X, Y, Z)=0$ for any $X, Y, Z\in\fr{X}(M)$ (notice that for this step the identity $\xi_j\lrcorner N^{(1)}=0$ is necessary). 
This shows that  $(M^{2n+s},  \phi, \xi_i, \eta_j, g)$ is an $\mc{S}$-manifold. The converse relies on the fact that for a $\mc{K}$-manifold $(M^{2n+s},  \phi, \xi_i, \eta_j, g)$ the characteristic vector fields $\xi_1, \ldots, \xi_s$ are all Killing, see \cite{Blair70}, and $\mc{S}$-manifolds are special examples of $\mc{K}$-manifolds. Finally, let us mention that the condition $\nabla T=0$ is an immediate consequence of the $\nabla$-parallelism of  the fundamental 2-form $F$ and of $\eta_i$ for all $i\in\{1, \ldots, s\}$.
\pro

\bp\label{Torsion_Smnfd}
Let $(M^{2n+s},  \phi, \xi_i, \eta_j, g)$ be an $\mc{S}$-manifold of CR-codimension $s$ endowed with the unique almost metric $\phi$-connection $\nabla$
with skew-torsion $T$ presented in  Theorem \ref{S-mnfds1}. Then for the corresponding vector-valued 2-form $T$ the following hold:  \begin{enumerate}
\item $T(X, Y)=2\sum_{j=1}^{s}\big\{F(X, Y)\xi_j-\eta_{j}(X)\phi(Y)+\eta_j(Y)\phi(X)\big\}$, for all $X, Y\in\fr{X}(M)$.
\item $T(X, \xi_i)=2\phi(X)$ for all $X\in\mc{D}$ and $i=1, \ldots, s$, and  $T(\xi_i, \xi_j)=0$ for all $i, j\in\{1, \ldots, s\}$.
\item $T(X, Y)=2\sum_{i=1}^{s}F(X, Y)\xi_{i}\in\mc{D}^{\perp}$ for all $X, Y\in\mc{D}$ and thus $T(X, Y, Z)=0$ for all $X, Y, Z\in\mc{D}$.
\end{enumerate}
Moreover,   $\dd T=2\sigma_{T}$, where
 $\sigma_T$ is the 4-form  introduced in Section \ref{subparagraph_Skew_Torsion}.
\ep

% For Sasakian manifolds $M^{2n+1}$ it is well known that the kernel $\Ker(T)$ of the $\nabla$-parallel torsion 3-form is trivial. 
%For higher CR-codimensions $s>1$ it follows that all the geometries with parallel skew-torsion  induced by an $\mc{S}$-manifold $M^{2n+s}$ are  \textsf{degenerate}, in the   sense that the 3-form $T$ has non-trivial kernel. % $\Ker(T_x)\neq 0$ at any $x\in M^{2n+s}$.   
%This is the content of the following proposition:
Let us now describe a result regarding the kernel of the torsion 3-form $T$ on a $\mc{S}$-manifold.
\bp\label{degeneracy}
Let $(M^{2n+s}, \phi, \xi_i, \eta_j, g)$ be an $\mc{S}$-manifold of CR-codimension $s\geq 1$, endowed with the unique almost metric $\phi$-connection $\nabla$
with (parallel) skew-torsion $T$, given  in  Theorem \ref{S-mnfds1}, that is,  $T=2\bar{\eta}\wedge F$ where $\bar{\eta}=\sum_{i=1}^{s}\eta_j$. Then,    no non-zero horizontal vector   $U\in \mc{D}_x$ lies in the  kernel $\Ker(T_x)$ of $T_x$, at any $x\in M$.  Moreover,  $\Ker(T_x)$ coincides with  the orthogonal complement
of the vertical vector $\bar{\xi}_x=\sum_{j=1}^{s}\xi_j(x)$ in $\mc{D}^{\perp}_{x}=\langle\xi_1(x), \ldots, \xi_s(x)\rangle$, and in particular  there is an isomorphism $\Ker(T_x)\cong\R^{s-1}$ at any $x\in M$.  Therefore the corresponding distribution $\Ker(T)$ is a smooth subbundle of $\mc{D}^{\perp}$ of rank $s-1$.
\ep

\pr
%At any point $x\in M$ we have the decomposition $T_xM=\mc{D}_x\oplus\mc{D}_x^{\perp}$ with
%$\mc{D}_x=\bigcap_{i=1}^{s}\Ker(\eta_i)_x$ and $\mc{D}^{\perp}_{x}=\langle \xi_1(x), \ldots, \xi_s(x)\rangle$, respectively.
At any point $x\in M$, it is easy to see that a tangent vector $U\in T_xM$ satisfies $U\in \mc{D}^{\perp}_{x}$ if and only if $U\lrcorner F_x=0$.  Now, the torsion 3-form $T$ at  $x\in M$  is an exterior 3-form $T_x \in\Lambda^{3}T_x^*M$ and
by its expression one deduces that %the expression of $T$ we deduce that 
% its kernel 
%is the subspace $\Ker(T_{x}):=\{U\in T_xM : U\lrcorner T_x=0\}\subset T_xM$. By the expression of $T$ we deduce that
\[
\Ker(T_x)=\{U\in T_xM : U\lrcorner T_x=0\}=\{U\in T_xM : 2\big(\bar{\eta}_x(U)F_x-\bar{\eta}_x\wedge(U\lrcorner F_x)\big)=0\}\subset T_xM\,.
\]
Thus, one may reduce the analysis to the following two cases:\\
\noindent a) $U\in\mc{D}_x$. Then $\bar{\eta}_{x}(U)=0$  since $(\eta_i)_{x}(U)=0$ for all $i\in\{1, \ldots, s\}$.
However, $U\lrcorner F_x\neq 0$ and thus in this case we get  $U\lrcorner T_x=-2\bar{\eta}_x\wedge(U\lrcorner F_x)$. 
 Then  an evaluation on a basis shows that this exterior 2-form is non-zero and  our first claim follows. \\ %at any $x\in M$ no non-zero horizontal vector $U\in D_x$ lies in the  kernel $\Ker(T_x)$.\\
 \noindent b)  $U\in \mc{D}_x^{\perp}$.  Then 
  $U\lrcorner F_x=0$, which means that  for vertical vectors we get the exterior 2-form $ U\lrcorner T_x=2\bar{\eta}_x(U)F_x$. 
Therefore, 
 $U\in\Ker(T_x)$ if and only if $\bar{\eta}_x(U)=0$. Let us express $U\in \mc{D}_x^{\perp}$ as  $U=\sum_{j=1}^{s}c_j\cdot\xi_j(x)$ for some scalars $c_1, \ldots, c_s$. Then   
 $\bar{\eta}_x(U)=\sum_{i=1}^{s}(\eta_{i})_{x}\big(\sum_{j=1}^{s}c_j\cdot\xi_j(x)\big)=\sum_{i, j=1}^{s}c_{j}\cdot(\eta_{i})_{x}\big(\xi_j(x)\big)=\sum_{i=1}^{s}c_{i}$, hence $
 \Ker(T_x)=\{U=\sum_{i=1}^{s}c_i\xi_i(x)\in \mc{D}_{x}^{\perp}: \sum_{i=1}^{s}c_{i}=0\}\cong\R^{s-1}$. 
This also shows that at any $x\in M$ the vertical vector  $\bar{\xi}_x=\sum_{j=1}^{s}\xi_j(x)$ does {\it not} belong to the kernel of $T_x$.  Since  $\bar{\eta}_x$ is   the dual  form of $\bar{\xi}_x$ with respect to $g_x$, we thus deduce that
   \[
 \Ker(T_x)=\{U\in \mc{D}_x^{\perp} : \bar\eta_x(U)=0\}=\{U\in \mc{D}_x^{\perp} : g_x(U, \bar{\xi}_x)=0\}=\bar{\xi}_x^{\perp}\cap\mc{D}_x^{\perp}\,.
 \]  
Varying $x\in M$ smoothly we finally obtain a smooth subbundle of  $\mc{D}^{\perp}$ of rank $(s-1)$,  $\Ker(T)=\bigsqcup_{x\in M}\Ker(T_x)\subset\mc{D}^{\perp}$.   This completes the proof.
\pro

%\br
%Obviously, the analysis  given above applies also for   $s=1$  and  in this case we get $\Ker(T)=0$, which is a standard fact for Sasakian manifolds. 
For $s=1$ the above statement gives $\Ker(T)=0$, which is a standard fact for Sasakian manifolds. 
This is a   central difference between the geometries
with parallel torsion on Sasakian manifolds   and on   $\mc{S}$-manifolds of   higher  CR-codimension  $s\geq 2$.
%In fact, any $\mc{S}$-manifold $(M^{2n+s}, \phi, \xi_i, \eta_j, g)$ is locally a Riemannian product of a  $\sqrt{s}$-Sasakian manifold and an $(s-1)$-dimensional flat manifold, see \cite[Theorem 2.1]{TP13}. Thus only the case $s=1$ produces irreducible examples. 
 Another notable  difference is   that $\mc{S}$-structures  of CR-codimension $s\geq 2$ with parallel skew-torsion occur in both odd and even dimensions.   Below,
in Sections \ref{U2_example!} and \ref{HeisenbergxT},  we  present examples illustrating these facts for the cases  $(n=1, s=2)$ and $(n=1, s=4)$, respectively. %while a study
%of the holonomy features of $\nabla$ on  general $\mc{S}$-manifolds with $s\geq 2$ will be presented in a separate work.
%\er

We end this section by comparing the characteristic connection $\nabla$ on  an $\mc{S}$-manifold with  the  so-called \textsf{Tanaka-Webster connection} (see \cite{T89, Miz93, LP04} for  more details on this connection).

\br\label{comparison}
In \cite[Theorem.~2.2]{LP04}  Lotta and Pastore gave a geometric interpretation of   CR-integrable contact metric $f$-manifolds (or CR-integrable almost $\mc{S}$-manifolds) 
in terms of the Tanaka-Webster connection $\tilde\nabla$. This connection $\tilde\nabla$ preserves the defining tensors $\phi$, $\xi_i$, $\eta_i$ and $g$, and  provides a  higher-dimensional analogue of the  Tanaka-Webster connection that admits a contact metric manifold in the case $s=1$ (see \cite[Proposition.~3.1]{T89}).
 Lotta and Pastore  also proved that a CR-integrable  almost $\mc{S}$-manifold is normal, i.e., an $\mc{S}$-manifold, if and only if the torsion form $\tilde{T}$ of $\tilde{\nabla}$
satisfies $\tilde T(\xi_i, X)=0$ for all $X\in\Gamma(\mc{D})$ and $i\in\{1, \ldots, s\}$, see \cite[Corollary 2.5]{LP04}. If this is the case, then the connection $\tilde\nabla$ is expressed by
$\tilde\nabla = \nabla^g+\tilde A$, where the $(1, 2)$-tensor field $\tilde A$  is given by $\tilde A(X, Y)=\sum_{j=1}^{s}\big\{F(X, Y)\xi_j+\eta_{j}(X)\phi(Y)+\eta_j(Y)\phi(X)\big\}$ for any $X, Y\in\fr{X}(M)$. 

On the other hand, in this work we proved    that an almost $\mc{S}$-manifold $(M^{2n+s}, \phi, \xi_i, \eta_j, g)$ admits an almost metric $\phi$-connection with skew-torsion $T$ if and only if  $M^{2n+s}$  is an $\mc{S}$-manifold (see Theorem \ref{S-mnfds1}), and then the characteristic connection $\nabla$ has torsion satisfying
$T(X, \xi_i)=2\phi(X)$ for all $X\in\Gamma(\mc{D})$ and all $i\in\{1, \ldots, s\}$ (see Proposition \ref{Torsion_Smnfd}).  Moreover,   $\nabla$ can be expressed as $\nabla=\nabla^g+\frac{1}{2}T=\nabla^g+A$, where in this case the $(1, 2)$-tensor field $A$ has the form $A(X, Y)=\sum_{j=1}^{s}\big\{F(X, Y)\xi_j-\eta_{j}(X)\phi(Y)+\eta_j(Y)\phi(X)\big\}$, for any  $X, Y\in\fr{X}(M)$. 
  Thus, the adapted connection $\nabla$  introduced in this article 
 is distinct from the Tanaka-Webster connection $\tilde{\nabla}$  presented in \cite{LP04}, as metric connections are uniquely determined by their torsion tensor.  
In fact, at least for $s=1$ it is known that the Tanaka-Webster connection has torsion of mixed type.
\er

%%%%%%%%%%%%%%%%%%%%%%%%%%%%%%%%
%%%%%%%%%%Examples%%%%%%%%%%%%
%%%%%%%%%%%%%%%%%%%%%%%%%%%%%%%%%

\section{Examples}\label{Section4}
In this  final section we describe    examples of metric $f$-structures admitting a characteristic connection with skew-torsion.
These examples arise either from the use of the compact Lie group $\U(n)$ for appropriate $n$, or through a construction of $\mc{S}$-manifolds outlined in \cite{DL05}, as described in the following remark:
\br\label{DL_rem}
Let $(N^{2n+1}, \phi, \xi, \eta, g)$ be a Sasakian manifold and let  $(\mc{D}=\Ker(\eta), \phi|_{\mc{D}})$ be the associated 
integrable CR-structure. Let $g_1$ be the (Hermitian) metric on $N$ with
\[
g_1(X, Y)=g(X, Y)\,,\ \text{for}\  X, Y\in\Gamma(\mc{D})\,, \ \ \text{and} \ \ g_1(X, \xi)=0\,,\quad g_1(\xi, \xi)=s\,.
\]
Consider the product manifold $M^{2n+s}:=N\times \Gg^{s-1}$, of CR-codimension $s$, where 
$\Gg$ is any commutative Lie group of dimension $s-1$.  Dileo and Lotta proved in \cite[Prop.~2.4]{DL05}  that $M^{2n+s}$   endowed with the metric $h=g_1\oplus g_2$, where $g_2$ is any  left-invariant metric of $\Gg$, 
is an $\mc{S}$-manifold. For instance, any product of the form $\Ss^{2n+1}\times\Tg^{s-1}$ is a compact $\mc{S}$-manifold of CR-codimension $s$.  
\er

% Remark \ref{DL_rem}.   
% Another class of examples based on a Riemannian submersion
%that admits any $\mc{S}$-manifold over a K\"ahler manifold (see \cite{Blair70, BLY73}), will be described in a forthcoming work.

\subsection{The Lie group $\U(2)$ as an Ambrose-Singer manifold}\label{U2_example!}
\subsubsection{The Lie group $\U(2)$ as an $\mc{S}$-manifold}
Consider the 4-dimensional Lie group $M^4=\U(2)$ and let $\fr{u}(2)=\Lie(\U(2))$ be its Lie algebra.
  \br In the sequel, all computations will be carried out at the identity element $e\in\U(2)$,  although we will use this fact mostly implicitly and usually avoid to mention that our tensor fields are evaluated at $e$. Since finally we will work only with Lie groups, we will adopt this approach for all the examples presented below.
 \er
Let $E_{ij}$ be the $2\times 2$ matrix with $1$ in the $(i, j)$ entry and zero everywhere else.
As a basis of $\fr{u}(2)$   consider the  left-invariant vector fields $\{X_{12}, Y_{12}, \xi_1, \xi_2\}$
with
\begin{eqnarray*}
&& (X_{12})_{e}=E_{12}-E_{21}=\begin{pmatrix}
0 & 1 \\ 
-1 & 0 \end{pmatrix}\,,\quad \  (\xi_1)_{e}=\sqrt{-1}E_{11}=\begin{pmatrix}
i & 0 \\ 
0 & 0 \end{pmatrix}\,,\\
&& (Y_{12})_{e}=\sqrt{-1}(E_{12}+E_{21})=\begin{pmatrix}
0 & i \\ 
i & 0 \end{pmatrix}\,,\quad  (\xi_2)_{e}=-\sqrt{-1}E_{22}=\begin{pmatrix}
0 & 0 \\ 
0 & -i \end{pmatrix}\,.
\end{eqnarray*}
The Lie brackets are given by 
\[
[X_{12}, Y_{12}]=2(\xi_1+\xi_2)\,,\quad [X_{12}, \xi_i]=-Y_{12}\,,\quad  [Y_{12}, \xi_i]=X_{12}\,,\quad [\xi_1, \xi_2]=0\,,
\]
where $i\in\{1, 2\}$.
Let $g$ be a left-invariant Riemannian metric on $\U(2)$ such that $\{X_{12}, Y_{12}, \xi_1, \xi_2\}$  is an orthonormal basis,
\[
g(X_{12}, X_{12})=1=g(Y_{12}, Y_{12})=g(\xi_i, \xi_i)\,\quad g(X_{12}, Y_{12})=0=g(\xi_1, \xi_2)=g(X_{12}, \xi_i)=g(Y_{12}, \xi_i)
\]
for all $i\in\{1, 2\}$.
Denote by $X_{12}^{\flat}, Y_{12}^{\flat}$ the dual 1-forms corresponding to $X_{12}, Y_{12}$ and by $\eta_1:=\xi_1^{\flat}, \eta_2:=\xi_2^{\flat}$ the dual 1-forms
corresponding to $\xi_1, \xi_2$, respectively.
Then we obtain a metric $f$-structure $(\phi, \eta_1, \eta_2, \xi_1, \xi_2, g)$, where the endomorphism $\phi$ is defined by  $\phi(X_{12})=Y_{12}$, $\phi(Y_{12})=-X_{12}$, $\phi(\xi_1)=\phi(\xi_2)=0$, such that
\[
TM=\mc{D}\oplus\mc{D}^{\perp}\,,\quad \mc{D}=\Imm(\phi)=\langle X_{12}, Y_{12}\rangle\,,\quad \mc{D}^{\perp}=\Ker(\phi)=\langle\xi_1, \xi_2\rangle\,.
\]
Of course, the fundamental 2-form has only one horizontal non-trivial component, $F(X_{12}, Y_{12})=-1$ and it is easy to see that $F=-X_{12}^{\flat}\wedge Y_{12}^{\flat}$ (note that   our convention is $u\wedge v=u\otimes v-v\otimes u$). Now,  by \cite{TK07} it is known that $(\U(2), \phi, \eta_1, \eta_2, \xi_1, \xi_2, g)$ is an even-dimensional $\mc{S}$-manifold (though for topological reasons it cannot be a K\"ahler manifold).  To realize the $\mc{S}$-structure, observe  that since  $\eta_i$ are all left-invariant,  
one computes
\[
\dd\eta_{1}(X_{12}, Y_{12})=-\eta_1([X_{12}, Y_{12}])=-\eta_{1}(2(\xi_1+\xi_2))=-2
\]
 and similarly for $\eta_2$. It follows that $2F=\dd\eta_1=\dd\eta_2$, while the vanishing of $N^{(1)}$ is direct. 
 \subsubsection{The adapted metric $\phi$-connection with skew-torsion}  Let us now describe the characteristic connection $\nabla$ on  the $\mc{S}$-manifold   $(\U(2), \phi, \xi_i, \eta_i, g)$.  
 \bp\label{Torsion_U2_expr}
Consider $M^{4}=\U(2)$ endowed with the  $\mc{S}$-structure  $(\phi, \xi_i, \eta_i, g)$ of CR-codimension 2 described above, and set   $\bar{\eta}:=\eta_1+\eta_2$. Then the following hold:
\begin{enumerate}
\item $(\U(2), \phi, \xi_i, \eta_i, g)$  admits a unique almost metric $\phi$-connection $\nabla$ with skew-torsion given by $\nabla=\nabla^g+\frac{T}{2}$, where  the torsion 3-form $T$ is expressed by
\[
T=2\bar{\eta}\wedge F=-2(\eta_1+\eta_2)\wedge X_{12}^{\flat}\wedge Y_{12}^{\flat}\,.
\]
\item  The corresponding vector-valued 2-form $T$  satisfies the relations
\[
 T(X_{12}, \xi_i)=2\phi(X_{12})=2Y_{12}\,,\quad T(Y_{12}, \xi_i)=2\phi(Y_{12})=-2X_{12}\,,\quad T(X_{12}, Y_{12})=-2(\xi_1+\xi_2)\,.
\]
\item The 3-form $T$ is $\nabla$-parallel, $\Ker(T)\subset\mc{D}^{\perp}$ is  a subbundle of rank 1, and  the (squared) norm $\|T\|^2$  is given by $\|T\|^2=8$ {\rm (}see the relation {\rm(\ref{normT})} for the definition of $\|T\|^2${\rm)}.
\item  The associated 4-form $\sigma_T$ vanishes and hence the torsion 3-form $T$ is  closed, $\dd T=\sigma_T=0$. 
\end{enumerate}
\ep

\subsubsection{The curvature tensor and the Ambrose-Singer structure}
Next we want to compute the curvature tensor of $\nabla$, which is defined by  $R^{\nabla}(X, Y)Z=\nabla_{X}\nabla_{Y}Z-\nabla_{Y}\nabla_{X}Z-\nabla_{[X, Y]}Z$.  This will allow us to show that $\nabla R^{\nabla}=0$. Again all computations are carried out   at the identity element. We first compute the expression of $\nabla$ and $\nabla^g$ on the $g$-orthonormal basis $\{X_{12}, Y_{12}, \xi_1, \xi_2\}$ (for $\nabla^g$ we are based on the  Koszul formula). 
 \bl\label{LCandT_u2}
 The Levi-Civita connection on the $\mc{S}$-manifold $(\U(2), \phi, \xi_1, \xi_2, \eta_1, \eta_2, g)$ satisfies the relations
\[
\begin{array}{l | l | l}
 \nabla^{g}_{X_{12}}Y_{12}=\xi_1+\xi_2=-\nabla^g_{Y_{12}}X_{12}   &  \nabla^g_{X_{12}}\xi_i=-Y_{12}  &   \nabla^{g}_{Y_{12}}\xi_i=X_{12}  \\
 \nabla^{g}_{X_{12}}X_{12}=0=\nabla^{g}_{Y_{12}}Y_{12} & \nabla^{g}_{\xi_i}X_{12}=0=\nabla^g_{\xi_i}Y_{12}  & \nabla^{g}_{\xi_i}\xi_j=0
\end{array}
 \]
 where  $i, j\in\{1, 2\}$. 
Moreover, the almost metric $\phi$-connection $\nabla=\nabla^g+\frac{T}{2}$ is such that
 \[
\begin{array}{l | l | l}
\nabla_{X_{12}}Y_{12}=0=\nabla_{Y_{12}}X_{12} & \nabla_{X_{12}}\xi_i=0=\nabla_{Y_{12}}\xi_i\ & \nabla_{X_{12}}X_{12}=0=\nabla_{Y_{12}}Y_{12}\\
 \nabla_{\xi_i}X_{12}=-Y_{12} & \nabla_{\xi_i}Y_{12}=X_{12} & \nabla_{\xi_i}\xi_j=0
\end{array}
 \]
 where $i, j\in\{1, 2\}$.
 \el
   %Using this lemma a direct calculation shows that
  \bp\label{U2_curvature}
  Consider the $\mc{S}$-manifold $(\U(2), \phi, \xi_i, \eta_i, g)$ $(i=1, 2)$ with the orthonormal basis $\{e_1=X_{12},\ e_2=Y_{12},\ e_3=\xi_1,\  e_4=\xi_2\}$ and the almost metric $\phi$-connection $\nabla=\nabla^g+\frac{T}{2}$  of Proposition \ref{Torsion_U2_expr}. Let $X=\sum_{i=1}^4 f_i e_i$, $Y=\sum_{j=1}^4 g_j e_j$ and $Z=\sum_{k=1}^4 h_k e_k$ be arbitrary left-invariant vector fields of $\U(2)$, for some constants $f_i, g_j, h_k$, with $i, j, k\in\{1, \ldots, 4\}$. Then the following hold:
  \begin{enumerate}
  \item The only nonzero curvature values of $R^{\nabla}$  are
\[
R^\nabla(e_1,e_2)e_1=4e_2=4\phi(e_1),\qquad R^\nabla(e_1,e_2)e_2=-4e_1=4\phi(e_2)
\]
and the corresponding sign  variants obtained by the antisymmetry in the first two arguments. In particular,  we have $R^{\nabla}(e_1, e_2)=4\phi$ and $R^{\nabla}(e_3, e_4)=0$. 
\item More in general,   $R^{\nabla}(X,Y)Z \;=\; 4\,(f_1 g_2 - f_2 g_1)\,(h_1 e_2 - h_2 e_1)\in\mc{D}$, and hence 
  $R^{\nabla}(X, Y)Z$ depends only on the horizontal components of $X, Y, Z$.
In particular, we can express $R^{\nabla}$ as $R^{\nabla}(X,Y)Z = -4\,F(X,Y)\,\phi(Z)$, or in other words  $R^{\nabla} = -4\,F\otimes\phi$. 
\end{enumerate}
\ep
\pr
(1) This part is   direct  and relies  on the expressions $\nabla_{e_{i}}e_j$ given in Lemma \ref{LCandT_u2}.\\
(2) Obviously,  $R^{\nabla}(X, Y)=\sum_{i, j}f_{i}g_{j}R^{\nabla}(e_i, e_j)=\sum_{i<j}(f_{i}g_{j}-f_jg_{i})R^{\nabla}(e_i, e_j)$, and thus 
\[
R^{\nabla}(X, Y)Z=\sum_{i<j}(f_{i}g_{j}-f_jg_{i})R^{\nabla}(e_i, e_j)Z\,.\quad (\ast)
\]
However we have $R^{\nabla}(e_i, e_j)\neq 0$ if and only if $(i,  j)=(1, 2)$ with $R^{\nabla}(e_1, e_2)=4\phi$, which means that
$R^{\nabla}(e_1, e_2)Z=4\phi(\sum_{k=1}^{4}h_k e_k)=4(h_1e_2-h_2e_1)$. 
Combining this with $(\ast)$ we obtain the given expression for $R^{\nabla}(X,Y)Z$. The final expression $R^{\nabla}=-4F\otimes\phi$
is based on the observation that $F(X,Y)=f_2 g_1 - f_1 g_2$ and  $\phi(Z)=h_1 e_2 - h_2 e_1$.
\pro

The $\nabla$-parallelism of $R^{\nabla}$ now follows as an important, yet straightforward, consequence of the formulas derived above.  In particular, 
 \bt\label{hol1}
 The holonomy algebra $\fr{hol}(\nabla)$ is isomorphic to $\fr{so}(2)\cong\fr{u}(1)$ and   the curvature tensor $R^{\nabla}$ is $\nabla$-parallel. 
  \et
 \pr
The $\nabla$-parallelism is an immediate consequence of the relation $R^{\nabla}=-4F\otimes\phi$ and the $\nabla$-parallelism of both $F$ and $\phi$. It also occurs by a more sophisticated way,  based on the holonomy algebra $\fr{hol}(\nabla)$ of $\nabla$.
Recall that  according to the Ambrose-Singer theorem (\cite{AS53}) the holonomy algebra $\fr{hol}(\nabla)$ is generated by the curvature endomorphisms $R^{\nabla}(X, Y)$ and their   conjugates under parallel transport.  Above we saw that any curvature endomorphism $R^{\nabla}(X, Y)$ is expressed as 
\[
R^{\nabla}(X, Y)=\sum_{i<j}(f_{i}g_{j}-f_jg_{i})R^{\nabla}(e_i, e_j)=(f_1g_2-f_2g_1)R^{\nabla}(e_1, e_2)=4(f_1g_2-f_2g_1)\phi
\]
and hence is a real scalar multiple of the endomorphism $\phi$.  
Since $\phi$ is $\nabla$-parallel,  parallel transport along any path  preserves $\phi$ and thus conjugation by parallel transport preserves the 1-dimensional subspace $\R\phi$. Hence we have $\fr{hol}(\nabla)\cong\langle\phi\rangle:=\R\phi$. Due to the expression of $\phi$ we see that $\langle\phi\rangle\subset\fr{so}(4)$ is isomorphic to $\fr{so}(2)\cong\fr{u}(1)$ and hence deduce that $\fr{hol}(\nabla)$ is 1-dimensional and abelian. The $\nabla$-parallelism of $R^{\nabla}$   then  follows from   part (3) of Proposition 2.2 in \cite{AFF15}.
  \pro

Based on this result and the $\nabla$-parallelism of $T$,  we deduce that  the connection $\nabla$ on the $\mc{S}$-manifold $(\U(2), \phi, \xi_i, \eta_i, g)$ $(i=1, 2)$ is an \textsf{Ambrose-Singer connection}, in particular, $(\U(2), g, \nabla)$  is an \textsf{Ambrose-Singer manifold} (we refer to \cite{CS04} for more details on Ambrose-Singer manifolds). 
Hence 
\bt
The Lie group $\U(2)$ endowed with the $\mc{S}$-structure $(\phi, \xi_i, \eta_i, g)$ $(i=1, 2)$ described above, and  the almost  metric $\phi$-connection $\nabla$ with parallel skew-torsion $T=\sum_{i=1}^{2}\eta_i\wedge\dd\eta_i$, is locally isometric to a naturally reductive homogeneous space. 
\et
Moreover,   $g$ is a left-invariant Riemannian metric and hence complete. Thus, it also  follows that the universal covering $\widetilde{\U}(2)=\SU(2)\times\R$ of $\U(2)$ is a naturally reductive space, a result known form 
 Kowalski and Vanhecke \cite{KV83}, see also \cite{AFF15}.

\subsubsection{The Ricci tensors} In this part we compute the Ricci tensor of $\nabla$ and compare it with that of the Levi-Civita connection $\nabla^g$.  
For the Riemannian curvature tensor and  the Ricci tensor $\Ric^g$ on $\mc{S}$-manifolds we refer to \cite{KT72, HOA86, CFF90} (see also   the Remark \ref{nabla_Einstein} below).

\br\label{curvature_symmetries}
Since $\nabla$ is a metric connection on $M^4=(\U(2), \phi, \xi_i, \eta_i, g)$, we have
\begin{equation}\label{symS}
 R^{\nabla}(X, Y, Z, V)=-R^{\nabla}(Y, X, Z, V)=-R^{\nabla}(X, Y, V, Z)\,,
 \end{equation}
 where   $R^{\nabla}(X, Y,Z, V):=g(R^{\nabla}(X, Y)Z, V)$. 
On the other hand, we saw that $\nabla T=0=\sigma_T$  and hence the following identity holds (see \cite{Srni})
\[
R^{\nabla}(X, Y, Z, V)=R^{g}(X, Y, Z, V)+\frac{1}{4}g(T(X, Y), T(Z, V))
\]
 for any $X, Y, Z, V\in\fr{u}(2)$. Moreover, the first Bianci identity for $R^{\nabla}$ has the classical form, i.e.,   $\fr{S}_{X, Y, Z}R^{\nabla}(X, Y, Z, V)=0$. 
 A combination of the first Bianci identity with the relations in (\ref{symS})  yields that $R^{\nabla}\in S^{2}\Lambda^2(M)$, which means that the curvature endomorphism $R^{\nabla} : \Lambda^{2}(M)\to\fr{hol}^{\nabla}\subset\Lambda^{2}(M)$  is actually a  symmetric endomorphism,  $g(R^{\nabla}(X, Y)Z, V)=g(R^{\nabla}(Z, V)X, Y)$ or in other words $R^{\nabla}(X, Y, Z, V)=R^{\nabla}(Z, V, X, Y)$. 
 
 Consider now the symmetric  tensor $S$  defined by $S(U,V):=\sum_{p=1}^4 g\big(T(e_p,U),\,T(e_p,V)\big)$, 
where $\{e_1:=X_{12},e_2:=Y_{12},e_3:=\xi_1, e_4:=\xi_2\}$ is the $g$-orthonormal frame.
Since the torsion $T$ is $\nabla$-parallel, the co-differential of $T$   vanishes and thus  the  Ricci tensor $\Ric^{\nabla}$ of $\nabla$ is symmetric, i.e.,
$\Ric^{\nabla}=\Ric^g-\frac{1}{4}S$,  see \cite{FrIv, AF04}.  
 \er

\bp\label{Ricci_tensors_U2}
Consider the $\mc{S}$-manifold $(\U(2), \phi, \xi_1, \xi_2, \eta_1, \eta_2, g)$. 
Then, relative to the $g$-orthonormal frame $\{e_1:=X_{12},e_2:=Y_{12},e_3:=\xi_1, e_4:=\xi_2\}$, the Ricci tensors $\Ric^g$, 
 $\Ric^{\nabla}$ and the symmetric   tensor  $S$ are given by the matrices
\[
\Ric^\nabla=\begin{pmatrix}
-4 & 0 & 0 & 0\\
0 & -4 & 0 & 0\\
0 & 0 & 0 & 0\\%[4pt]
0 & 0 & 0 & 0
\end{pmatrix},\quad 
\Ric^g=\begin{pmatrix}
0 & 0 & 0 & 0\\
0 & 0 & 0 & 0\\
0 & 0 & 2 & 2\\
0 & 0 & 2 & 2
\end{pmatrix},\qquad
S=\begin{pmatrix}
16 & 0 & 0 & 0\\
0 & 16 & 0 & 0\\
0 & 0 & 8 & 8\\
0 & 0 & 8 & 8
\end{pmatrix}.
\]
The scalar curvatures are given by  $\Scal^{\nabla}=-8$ and $\Scal^g=4$, respectively. 
\ep

\pr
Both the Ricci tensors are symmetric and  by Remark \ref{curvature_symmetries}  it holds that $\Ric^\nabla(U,V)=\sum_{p=1}^4 g\big(R^{\nabla}(e_p,U)V,e_p\big)=\sum_{p=1}^{4}g\big(R^{\nabla}(U, e_p)e_p, V\big)$. 
Similarly for $\Ric^{g}$. The  given matrix expression $\Ric^{\nabla}=(\Ric^{\nabla}_{ij})_{1\leq i, j\leq 4}$
   follows by a direct   computation
based on Proposition \ref{U2_curvature}, where $\Ric^{\nabla}_{ij}:=\Ric^\nabla(e_i, e_j)=\sum_{p=1}^{4}g\big(R^{\nabla}(e_i, e_p)e_p, e_j\big)$. 
For example,  we compute $\Ric^{\nabla}_{11}=\Ric^{\nabla}_{22}=-4$. In full details 
 \begin{eqnarray*}
 \Ric^{\nabla}(e_1, e_1)&=&\Ric^{\nabla}(X_{12}, X_{12})=g(R^{\nabla}(X_{12}, X_{12})X_{12}, X_{12})+g(R^{\nabla}(X_{12}, Y_{12})Y_{12}, X_{12})\\
 &&+g(R^{\nabla}(X_{12}, \xi_1)\xi_1, X_{12})+g(R^{\nabla}(X_{12}, \xi_2)\xi_2, X_{12})\\
 &=&-4g(X_{12}, X_{12})=-4\,,\\
 \Ric^{\nabla}(e_2, e_2)&=&\Ric^{\nabla}(Y_{12}, Y_{12})=g(R^{\nabla}(Y_{12}, X_{12})X_{12}, Y_{12})+g(R^{\nabla}(Y_{12}, Y_{12})Y_{12}, Y_{12})\\
 &&+g(R^{\nabla}(Y_{12}, \xi_1)\xi_1, Y_{12})+g(R^{\nabla}(Y_{12}, \xi_2)\xi_2, Y_{12})\\
 &=&-4g(Y_{12}, Y_{12})=-4\,,
 \end{eqnarray*}
  while all the other entries of $\Ric^{\nabla}$ vanish. 
Similarly we have   $\Ric^{g}=(\Ric^{g}_{ij})_{1\leq i, j\leq 4}$,  with $\Ric^{g}_{ij}:=\Ric^{g}(e_i, e_j)=\sum_{p=1}^{4}g\big(R^{g}(e_i, e_p)e_p, e_j\big)$. 
For this part we mention  that the only  non-zero components of the  Riemannian curvature tensor $R^{g}$ are given by
  \[
\begin{tabular}{l | l | l}
$R^g(X_{12},Y_{12})X_{12}= 2Y_{12}=2\phi(X_{12})\,,$ &  $R^g(X_{12}, \xi_i)X_{12}=-(\xi_1+\xi_2)\,,$ & $R^g(X_{12},\xi_i)\xi_j=X_{12}\,,$ \\
$R^g(X_{12},Y_{12})Y_{12}= -2X_{12}=2\phi(Y_{12})\,,$ & $R^g(Y_{12}, \xi_i)Y_{12}=-(\xi_1+\xi_2)\,,$ &  $R^g(Y_{12},\xi_i)\xi_j=Y_{12}\,,$ 
\end{tabular}
\]
 where $i, j\in\{1, 2\}$, together with the corresponding sign  variants obtained by the antisymmetry in the first two arguments. To obtain these expressions we employed   
Lemma  \ref{LCandT_u2}.\\
\noindent The verification of the relation  $\Ric^{\nabla}=\Ric^g-\frac{1}{4}S$  relies now on the matrix    $S=(S_{ij})_{1\leq i, j\leq 4}$,  with  $S_{ij}:=S(e_i, e_j)=\sum_{p=1}^{4}g(T(e_i, e_p), T(e_j, e_p))$. 
 To obtain these entries  we use the relations in part (2) of  Proposition \ref{Torsion_U2_expr}. 
 Finally, taking traces gives the scalar curvatures.   
\pro

\br\label{nabla_Einstein}
(1) The  values  or $\Scal^{\nabla}$ and $\Scal^g$ obtained above  in combination with the relation  $\|T\|^2=8$ from Proposition \ref{Torsion_U2_expr}, 
 confirm the standard identity $\Scal^{\nabla}=\Scal^g-\frac{3}{2}\|T\|^2$  (see   \cite{AF04}).\\
(2) By \cite{KT72}  (see also \cite[Prop.~3.5]{CFF90}) it is known that  an $\mc{S}$-manifold $(M^{2n+s}, \phi, \xi_i, \eta_j, g)$ with $s\geq 2$ is never Einstein. This is because the Ricci tensor $\Ric^g$ satisfies $\Ric^{g}(X, \xi_i)=2 n\sum_{j=1}^{s}\eta_{j}(X)$ for all $X\in\Gamma(TM)$ and $i=1, \ldots, s$.  For the case of $\U(2)$,  using this identity one can confirm   all the given entries $\Ric^{g}_{ij}$ of $\Ric^g$, 
where at least one of $i, j$ is such that $i, j\in\{3, 4\}$. In particular, 
\[
\Ric^{g}(X, \xi_i)=\begin{cases} 
0, & \text{for} \  X\in\Gamma(\mc{D}), \\
2, & \text{for} \  X=e_3=\xi_1 \ \text{or} \ X=e_4=\xi_2.
\end{cases}
\]
(3)  By Proposition \ref{Ricci_tensors_U2}, the restriction of the Ricci tensor $\Ric^{\nabla}$ to the horizontal distribution $\mc{D}$ is a multiple of the restricted metric $g|_{\mc{D}}$, while its restriction  
to the vertical part $\mc{D}^{\perp}$ vanishes. Therefore, the $\mc{S}$-manifold $(\U(2), \phi, \xi_1, \xi_2, \eta_1, \eta_2, g)$ endowed with the connection $\nabla$  does not furnish an example   of a $\nabla$-Einstein manifold with (parallel) skew-torsion. 
\er

\subsection{$\mc{S}$-manifolds of the form $H_{2n+1}\times\Tg^{s-1}$}\label{HeisenbergxT}
Let us now describe examples based on  Remark \ref{DL_rem}. We will  discuss   an example of CR-codimension 4 (and of CR-dimension $1$), which naturally extends to general CR-dimension $n$ and general CR-codimension $s$. 

\subsubsection{The 6-dimensional product $M^{2\cdot 1+4}=H_3\times\Tg^3$} 
Let us consider the 3-dimensional Heisenberg group $H_{3}$. Its Lie algebra $\fr{h}_3=\Lie(H_3)$ is generated by three left-invariant vector fields $\{X, Y, Z\}$, with $[X, Y]=Z$ being the only non-zero Lie bracket.  A  Sasakian structure is defined by $\hat{\phi}(X)=Y$, $\hat{\phi}(Y)=-X$, and $\hat{\phi}(Z)=0$ with    Reeb vector field $\xi:=Z$ and horizontal distribution $\hat{\mc{D}}=\Imm(\hat{\phi})$.  The corresponding left-invariant  Sasakian metric $g$ on $H_{3}$ is associated to  an inner product in $\fr{h}_3$, which we still denote by $g$, and  can be obtained   by declaring the relations
\[
g(X, Y)=0=g(X, Z)=g(Y, Z)\,, \quad  g(X, X)=\frac{1}{2}=g(Y, Y)\,, \quad g(Z, Z)=1\,,
\] 
such that $\dd\hat{\eta}=2\hat{F}$.  
Consider also the 3-dimensional torus $\Tg^3$ endowed with a left-invariant metric $g_2$, and let $\{\zeta_1, \zeta_2, \zeta_3\}$ be a  $g_2$-orthonormal frame of $\Tg^3$.
Obviously,  $\bb{V}:=\langle Z, \zeta_1, \zeta_2, \zeta_3\rangle$ is a 4-dimensional commutative Lie algebra.

 Consider now the 6-dimensional product $M=H_3\times\Tg^3$ endowed with the metric $h=g_1\oplus g_2$, where $g_1$ is as in  Remark \ref{DL_rem}, i.e.,
 \[
 g_1(U, V)=g(U, V)\,,\ \text{for}\  U, V\in\Gamma(\hat{\mc{D}})\,, \ \ \text{and} \ \ g_1(U, Z)=0\,,\quad g_1(Z, Z)=4\,.
 \]
The metric $h$ satisfies $h(Z, Z)=4$,  hence there exist $h$-orthonormal  vector fields $\xi_1, \xi_2, \xi_3, \xi_4\in\bb{V}$ such that $Z=\sum_{i=1}^{4}\xi_i$.  
We  define an $f$-structure $\phi : TM\to TM$ on $M^{2+4}$ by extending the Sasakian structure $\hat\phi$ as follows: $\phi(X):=\hat\phi(X)=Y$, $\phi(Y):=\hat\phi(Y)=-X$ and $\phi(\xi_i):=0$ for all $i=1, \ldots, 4$. Then  we get the $h$-orthogonal decomposition $TM=\mc{D}\oplus\mc{D}^{\perp}$ with 
\[
\mc{D}=\Imm(\phi)=\Ker(\eta_1)\cap\cdots\cap\Ker(\eta_4)\,,\quad \mc{D}^{\perp}=\Ker(\phi)=\langle \xi_1, \xi_2, \xi_3, \xi_4\rangle\,,
\] 
where $\eta_i$ are the dual 1-forms of $\xi_i$. Let  $X^{\flat}$, $Y^{\flat}$ be the dual 1-forms of $X, Y$, as well.  The fundamental 2-form satisfies $F(X, Y)=h(X, \phi Y)=-g_1(X, X)=-g(X, X)=-\frac{1}{2}$.
Since $X^{\flat}(X)=\frac{1}{2}=Y^{\flat}(Y)$, we have $(X^{\flat}\wedge Y^{\flat})(X, Y)=\frac{1}{4}$, and   hence $F=-2(X^{\flat}\wedge Y^{\flat})$. Since 
\[
\dd\eta_i(X, Y)=-\eta_i([X, Y])=-\eta_i(Z)=-\eta_i(\sum_{j=1}^{4}\xi_j)=-\sum_{j=1}^{4}\eta_{i}(\xi_j)=-\sum_{j=1}^{4}h(\xi_i, \xi_j)=-1
\]
for all $i=1,\ldots, 4$, we get $\dd\eta_i=-4(X^{\flat}\wedge Y^{\flat})$.  Thus $2F=\dd\eta_i$ for all $i=1,\ldots, 4$ and $M^6$ is a contact metric $f$-manifold of CR-codimension 4.
Moreover, $N^{(1)}$ vanishes, hence $M$ is an $\mc{S}$-manifold.  
\bp\label{H3T3_torsion}
 Consider the product $M=H_3\times\Tg^3$ endowed with the $\mc{S}$-structure $(\phi, \xi_i, \eta_j, h)$ of CR-codimension 4 described above, and set  $\bar{\eta}=\eta_1+\cdots+\eta_4$. Then the following hold:
 \begin{enumerate}
 \item $(M=H_3\times\Tg^3, \phi, \xi_i, \eta_j, h)$ admits a unique almost metric $\phi$-connection  $\nabla$ with skew-torsion given by $\nabla=\nabla^g+\frac{T}{2}$, where the torsion 3-form $T$   is expressed by
\[
T=\sum_{i=1}^{4}(\eta_i\wedge\dd\eta_i)=2(\eta_1+\cdots+\eta_4)\wedge F=2\bar{\eta}\wedge F=-4\bar{\eta}\wedge X^{\flat}\wedge Y^{\flat}\,.
\]
\item The  corresponding vector-valued 2-form  $T$ satisfies
\[
T(X, \xi_i)=2\phi(X)=2Y\,,\quad T(Y, \xi_i)=2\phi(Y)=-2X\,, \quad T(X, Y)=-\sum_{i=1}^{4}\xi_i=-Z\,.
\]
\item The torsion 3-form $T$  is $\nabla$-parallel, $\nabla T=0$, and its (squared) norm $\|T\|^2$ is given by $\|T\|^2=16$.
\item  The kernel of the torsion 3-form $T$ is a  smooth subbundle of  $\mc{D}^{\perp}$ of rank 3. 
\item The 4-form $\sigma_T$ vanishes, $\sigma_T=0$, and hence  $\dd T=0$ as well. 
\end{enumerate}
\ep
 
% \br 
%1) An analogous construction clearly applies to the products $H_3\times\Tg^1$ and $H_3\times\Tg^2$: the first yields a 4-dimensional $\mc{S}$-manifold of CR-codimension 2, while the second  gives a 5-dimensional $\mc{S}$-manifold of CR-codimension 3 (hence  distinct  from a five-dimensional quasi-Sasaki manifold, which by definition has CR-codimension 1). On these spaces the left-invariant (parallel) torsion 3-form $T$  still takes  the form  $T=2\bar{\eta}\wedge F$, where $\bar{\eta}$ is defined analogously and $F$ denotes the corresponding fundamental 2-form. This is still given by $F=-2(X^{\flat}\wedge Y^{\flat})$. Moreover, the kernel $\Ker(T)$ is non-trivial.  The universal spaces $H_3\times\R$ and $H_3\times\R^2$ can likewise be regarded as (complete)
% $\mc{S}$-manifolds with   parallel degenerate skew-torsion $T=2\bar{\eta}\wedge F$.   
%This pattern naturally extends  to higher dimensions: products of the form $H_3\times\Tg^{s}$ yield $(3+s)$-dimensional $\mc{S}$-manifolds of CR-codimension $s+1$, all admitting parallel skew-torsion with non-trivial kernel. \\ 
%2)  We can also consider products $H_{2n+1}\times\Tg^{s-1}$, where $H_{2n+1}$ is the $(2n+1)$-dimensional  Heisenberg group endowed with its standard Sasakian structure,  or their universal coverings, i.e., products 
%of the form $H_{2n+1}\times\R^{s-1}$. These are both $(2n+s)$-dimensional $\mc{S}$-manifolds of CR-codimension  $s$ with parallel degenerate skew-torsion.  
%\er

 \subsubsection{The curvature tensor and the Ambrose-Singer structure} Next we will describe the curvature of $\nabla$ on $H_3\times\Tg^3$. 
 We first present the expressions of
 the Levi-Civita connection $\nabla^h$ and the unique almost metric $\phi$-connection $\nabla$ of Proposition \ref{H3T3_torsion} on the $h$-orthonormal basis used above. Again all our computations  are carried out  at the identity element of $H_3\times\Tg^3$.

\bl\label{productLCn}
Consider the $\mc{S}$-manifold $(M=H_3\times\Tg^3, \phi, \xi_i, \eta_j, h)$ of CR-codirmension 4 as described above, endowed with the $h$-orthonormal frame  $\{e_1=\sqrt{2}\,X,\ e_2=\sqrt{2}\,Y,\ 
e_3=\xi_1,\  e_4=\xi_2,\  e_5=\xi_3,\ e_6=\xi_4\}$.   
Then the Levi--Civita connection $\nabla^h$ satisfies the relations  
\[
\begin{array}{l | l | l}
 \nabla^h_{e_1} e_1 = 0       &  \nabla^h_{e_2} e_2 = 0    & \nabla^h_{e_i} e_1 = -e_2   \\
 \nabla^h_{e_1} e_2 = Z       &   \nabla^h_{e_2} e_1 = -Z  & \nabla^h_{e_i} e_2 = e_1    \\
 \nabla^h_{e_1} e_i = -e_2   &    \nabla^h_{e_2} e_i = e_1 & \nabla^h_{e_i} e_j = 0  
\end{array}
\]
where $i,j \in \{3,4,5,6\}$.  
Moreover, for the almost metric $\phi$-connection $\nabla=\nabla^g+\frac{T}{2}$ we see that
\[
\begin{array}{l | l | l}
 \nabla_{e_1} e_1 = 0  &  \nabla_{e_2} e_2 = 0  & \nabla_{e_i} e_1 = -2e_2  \\
 \nabla_{e_1} e_2 = 0  &  \nabla_{e_2} e_1 = 0 &   \nabla_{e_i} e_2 = 2e_1 \\
 \nabla_{e_1} e_i = 0  &  \nabla_{e_2} e_i = 0  & \nabla_{e_i} e_j = 0
\end{array}
\]
\el
In a similar way with Proposition \ref{U2_curvature} we can now prove that
\bp
Consider the $\mc{S}$-manifold $(M=H_3\times\Tg^3, \phi, \xi_i, \eta_i, g)$ of CR-codimension 4  endowed with the orthonormal frame 
$\{e_1=\sqrt{2}\,X,\ e_2=\sqrt{2}\,Y,\ 
e_3=\xi_1,\  e_4=\xi_2,\  e_5=\xi_3,\ e_6=\xi_4\}$ and the almost metric $\phi$-connection $\nabla=\nabla^h+\frac{T}{2}$ of Proposition \ref{H3T3_torsion}.
 Then the following hold:
\begin{enumerate}
\item The only nonzero curvature values of $R^{\nabla}$  are
\[
R^\nabla(e_1,e_2)e_1=16e_2=16\phi(e_1),\qquad R^\nabla(e_1,e_2)e_2=-16e_1=16\phi(e_2)
\]
and the corresponding sign  variants obtained by the antisymmetry in the first two arguments.  In particular,  we have $R^{\nabla}(e_1, e_2)=16\thinspace\phi$ and $R^{\nabla}(\xi_i, \xi_j)=0$ for all $i, j\in\{1, \ldots, 4\}$.
\item The curvature $R^{\nabla}$ is a scalar multiple of 	$F\otimes\phi$, i.e., $R^{\nabla}=-16 F\otimes\phi$.
\end{enumerate}
\ep
\pr
The computation of the curvature tensor $R^{\nabla}$ is direct and relies on Lemma \ref{productLCn}. 
We also mention that in terms of the $h$-orthonormal frame $\{e_i\}$ we have $F(e_1, e_2)=-1$ and hence $F=-e_1^{\flat}\wedge e_2^{\flat}=-e^{1}\wedge e^{2}$. Since $(F\otimes\phi)(U, V)W=F(U, V)\phi(W)$ we thus get 
\begin{eqnarray*}
(-16F\otimes\phi)(e_1, e_2)e_1=-16F(e_1, e_2)\phi(e_1)=16\phi(e_1)=16e_2\,,\\
(-16F\otimes\phi)(e_1, e_2)e_2=-16F(e_1, e_2)\phi(e_2)=16\phi(e_2)=-16e_1\,.
\end{eqnarray*}
The relation $R^{\nabla}=-16 F\otimes\phi$ is now direct.
\pro
Therefore, as in Theorem \ref{hol1}, we deduce that
\bt
The holonomy algebra $\fr{hol}(\nabla)$ of $\nabla$ is abelian and 1-dimensional, $\fr{hol}(\nabla)\cong\fr{so}(2)$, and moreover the curvature tensor $R^{\nabla}$ is $\nabla$-parallel, $\nabla R^{\nabla}=0$.
\et 
 
 It follows that  the connection $\nabla$ on the $\mc{S}$-manifold $(H_3\times\Tg^3, \phi, \xi_i, \eta_i, h)$ $(i=1,\ldots, 4)$ is an  Ambrose-Singer connection, in particular, $(H_3\times\Tg^3, h, \nabla)$  is an  Ambrose-Singer manifold.
 Hence 
\bt
The product $H_3\times\Tg^3$ endowed with the $\mc{S}$-structure $(\phi, \xi_i, \eta_i, h)$  of CR-codimension 4, as described above, and the almost metric $\phi$-connection $\nabla$ with  parallel skew-torsion $T=\sum_{i=1}^{4}\eta_i\wedge\dd\eta_i$, is locally isometric to a naturally reductive homogeneous space. 
\et
Since $h$ is  complete, it also follows that the universal covering $\widetilde{M}^6=H_3\times \R^3$ of $M$ is a naturally reductive space whose torsion has non-trivial kernel.\footnote{By \cite[p.~71]{AFF15} it is known that a naturally reductive space in six dimensions whose torsion has non-trivial kernel is locally a product of lower-dimensional
naturally reductive spaces.}

\subsubsection{The Ricci tensors} Let us finally present the Ricci tensors of $\nabla$ and $\nabla^g$ (which are both symmetric).
\bp
Consider the $\mc{S}$-manifold $(M=H_3\times\Tg^3, \phi, \xi_i, \eta_i, g)$ of CR-codimension 4. 
With respect to the $h$-orthonormal frame
$\{e_1=\sqrt{2}X, e_2=\sqrt{2}Y, e_3=\xi_1, e_4=\xi_2, e_5=\xi_3, e_6=\xi_4\}$
the Ricci tensor $\Ric^{\nabla}$ of the almost metric $\phi$-connection 
$\nabla$, the Riemannian Ricci tensor $\Ric^{g}$ and the symmetric tensor  $S$ satisfying  $\Ric^{\nabla}=\Ric^g-\frac{1}{4}S$,  are given by the following $6\times6$ matrices:
\[
\Ric^{\nabla}=\begin{pmatrix}
16 & 0 & 0 & 0 & 0 & 0\\[4pt]
0 & 16 & 0 & 0 & 0 & 0\\[4pt]
0 & 0 & 0 & 0 & 0 & 0\\[4pt]
0 & 0 & 0 & 0 & 0 & 0\\[4pt]
0 & 0 & 0 & 0 & 0 & 0\\[4pt]
0 & 0 & 0 & 0 & 0 & 0
\end{pmatrix},\quad 
\Ric^{g}=\begin{pmatrix}
24 & 0 & 0 & 0 & 0 & 0\\[4pt]
0 & 24 & 0 & 0 & 0 & 0\\[4pt]
0 & 0 & 2 & 2 & 2 & 2\\[4pt]
0 & 0 & 2 & 2 & 2 & 2\\[4pt]
0 & 0 & 2 & 2 & 2 & 2\\[4pt]
0 & 0 & 2 & 2 & 2 & 2
\end{pmatrix}, \quad 
S=\begin{pmatrix}
32 & 0 & 0 & 0 & 0 & 0\\[4pt]
0 & 32 & 0 & 0 & 0 & 0\\[4pt]
0 & 0 & 8 & 8 & 8 & 8\\[4pt]
0 & 0 & 8 & 8 & 8 & 8\\[4pt]
0 & 0 & 8 & 8 & 8 & 8\\[4pt]
0 & 0 & 8 & 8 & 8 & 8
\end{pmatrix}.
\]
Moreover,  the corresponding scalar curvatures are given by $\Scal^{\nabla}=32$ and $\Scal^g=56$.
\ep
\pr
The proof follows the same approach as those presented in Proposition \ref{Ricci_tensors_U2} and we will avoid to describe the details.
We only mention that  using Lemma \ref{productLCn} we see that 
\[
\begin{array}{lll}
R^{h}(e_1,e_2)e_1 = 12\,e_2, & R^{h}(e_1,e_2)e_2 = -12\,e_1, & R^{h}(e_1,e_2)e_i = 0,\\[4pt]
R^{h}(e_1,e_i)e_1 = -Z, & R^{h}(e_1,e_i)e_j = e_1, & R^{h}(e_i,e_j) = 0,\\[4pt]
R^{h}(e_2,e_i)e_2 = -Z, & R^{h}(e_2,e_i)e_j = e_2, & R^{h}(e_1,e_2)=-R^{h}(e_2,e_1),
\end{array}
\]
for all $i,j\in\{3,4,5,6\}$, and all other components vanish. Note also that the values presented for $\Scal^{\nabla}$ and $\Scal^g$  in combination with the relation  $\|T\|^2=16$ by Proposition \ref{H3T3_torsion}, 
verify the    identity $\Scal^{\nabla}=\Scal^g-\frac{3}{2}\|T\|^2$.
\pro

\subsection{A normal metric $f$-structure on the Lie group $\U(3)$}
Let us finally discuss  an example of a metric $f$-structure of CR-codimension 3,  based on the Lie group $\U(3)$. As we will see below, this metric $f$-structure  differs from the one on $\U(2)$ in the sense that it does not provide   an example of an $\mc{S}$-structure.

Consider the   Lie group $\U(3)$ 
and  denote by  $E_{ij}$  the  $3\times 3$ matrices    with $1$ in the $(i, j)$ slot and $0$ elsewhere. Define left-invariant vector field $X_{ij}$, $Y_{ij}$  with $(ij)=(12), (13), (23)$ and $\xi_i$ for $i=1, 2, 3$, 
by the relations
\[
(X_{ij})_{e}=E_{ij}-E_{ji}\,, \quad (Y_{ij})_{e}=\sqrt{-1}(E_{ij}+E_{ji})\,, \quad (\xi_{i})_{e}=\sqrt{-1}E_{ii}\,,
\]
where $e$ is the neutral element in $\U(3)$. Then the set  $\mc{B}:=\{X_{12}, Y_{12}, X_{13}, Y_{13}, X_{23}, Y_{23}, \xi_{1}, \xi_2, \xi_3\}$ 
 form a basis of the corresponding Lie algebra $\fr{u}(3)=\Lie(\U(3))$.
 For the Lie brackets we have $[\xi_i,  \xi_j]=0$ for all $1\leq i, j\leq 3$,
hence $\fr{t}^{3}:=\langle \xi_1, \xi_2, \xi_3\rangle$ coincides with the Lie algebra of a maximal torus $\Tg^{3}\subset \U(3)$ in $\U(3)$. Set $\mc{D}:=\langle X_{12}, Y_{12}, X_{13}, Y_{13}, X_{23}, Y_{23}\rangle$. 
 We see that $[X_{ij}, Y_{ij}]=2(\xi_i-\xi_j)\in\fr{t}^{3}$ and the rest Lie brackets belong to $\mc{D}$. For instance, for all $(ij)\in\{(12), (13), (23)\}$ we have  
\begin{equation}\label{LieXxi}
\begin{tabular}{ l ll l ll l}
$[X_{ij}, \xi_i]=-Y_{ij}\,,$ & $[X_{ij}, \xi_j]=Y_{ij}\,,$ & $[X_{ij},\xi_k]=0\,, \ (k\neq i, j)\,,$ \\  
$[Y_{ij}, \xi_i]=X_{ij}\,,$ & $[Y_{ij}, \xi_j]=-X_{ij}\,,$ & $[Y_{ij},\xi_k]=0\,, \ \thinspace (k\neq i, j)\,.$  
\end{tabular}
\end{equation}
We include all the non-zero Lie brackets (along with a few important zero ones) in the table below.
\[
\begin{tabular}{ l || l || l }
\thickline
$[X_{12}, Y_{12}]=2 (\xi_1 - \xi_2)$ & $[X_{13}, Y_{13}]=2(\xi_1-\xi_3)$ & $[X_{23}, Y_{23}]=2(\xi_2 - \xi_3)$ \\
$[X_{12}, X_{13}] =- X_{23}$ & $[X_{13}, X_{23}]=-X_{12}$  & $[X_{23}, Y_{12}]=-Y_{13}$  \\
$[X_{12}, X_{23}]=X_{13}$ &     {$[X_{13}, Y_{12}]=- Y_{23}$} & $[X_{23}, Y_{13}]=Y_{12}$ \\
$[X_{12}, Y_{13}]=- Y_{23}$ & $[X_{13}, Y_{23}]=Y_{12}$ & {$[Y_{12}, Y_{13}]=-X_{23}$}  \\
$[X_{12}, Y_{23}]=Y_{13}$  & $[Y_{13}, Y_{23}]=-X_{12}$ &   $[Y_{12}, Y_{23}]=- X_{13}$      \\ 
\hline
$[X_{12}, \xi_1]=-Y_{12}$  & $ [X_{13}, \xi_1]=-Y_{13}$    & $[X_{23}, \xi_1]=0$ \\
$[X_{12}, \xi_2]= Y_{12}$ & $[X_{13}, \xi_2]=0$   & $[X_{23}, \xi_2]=-Y_{23}$ \\
$[X_{12}, \xi_3]=0$ & $[X_{13}, \xi_3]= Y_{13}$  & $[X_{23}, \xi_3]=Y_{23}$  \\ 
\hline
$[Y_{12}, \xi_1]= X_{12}$ &  $[Y_{13}, \xi_1]= X_{13}$  & $[Y_{23}, \xi_1]=0$  \\
 $[Y_{12}, \xi_2]=-X_{12}$  & $[Y_{13}, \xi_2]=0$ & $[Y_{23}, \xi_2]=X_{23}$  \\
$[Y_{12}, \xi_3]=0$ & $[Y_{13}, \xi_3]=-X_{13}$ &    $[Y_{23}, \xi_3]=-X_{23}$ \\
\thickline
\end{tabular}
\]
Let us now consider the left-invariant Riemannian metric $g$ on $\U(3)$ which  makes the given  basis $\mc{B}$ an orthonormal one.     We denote by the same notation the  corresponding inner product in  $\fr{u}(3)$.
\bp
With respect to $g$ each  (Reeb) vector field $\xi_i$ is Killing.
\ep
\pr
The metric $g$ is left-invariant, hence $g(X, Y)$ is constant and  the condition $(\mc{L}_{\xi_i}g)(X, Y)=0$ is equivalent to say that $\ad(\xi_i)$ is skew-symmetric with respect to  $g$, i.e.,  $g(\ad(\xi_i)X, Y)+g(X,\ad(\xi_i)Y)=0$ 
for any two left-invariant vector fields $X, Y\in\fr{u}(3)$.
Then one can easily see that this holds as an identity due to the relations mentioned in (\ref{LieXxi}).
\pro

\noindent We now  define an $f$-structure on $\U(3)$ by
\[
\phi(X_{ij})=Y_{ij}\,, \quad \phi(Y_{ij})=-X_{ij}\,, \ \text{for} \ (ij)=(12), (13), (23)\,,\quad \phi(\xi_i)=0\,,\quad i=1, 2, 3\,.
\]
Then we obtain a $g$-orthogonal decomposition  $T\U(3)=\Imm(\phi)\oplus\Ker(\phi)=\mc{D}\oplus\mc{D}^{\perp}$ with   
\[
\Imm(\phi)=\mc{D}=\langle X_{12}, Y_{12}, X_{13}, Y_{13}, X_{23}, Y_{23}\rangle\,,\quad  
 \Ker(\phi)=\mc{D}^{\perp}=\langle \xi_1, \xi_2, \xi_3\rangle=\fr{t}^{3}\,,
 \]
  respectively. 
We denote by $X_{ij}^{\flat}$, $Y_{ij}^{\flat}$ the dual 1-forms of $X_{ij}, Y_{ij}$ for $(ij)\in\{(12), (13), (23)\}$
and by $\eta_i:=\xi_{i}^{\flat}$ the dual of the Reeb vector fields $\xi_i$ for $i=1, 2, 3$. 
Then it is easy to see that $(\phi, \xi_i, \eta_i, g)$ defines a (left-invariant) metric $f$-structure on $\U(3)$, that is, $g$ satisfies (\ref{gXgY}). Moreover, the 2-form $F$ is given by
\[
F=-X_{12}^{\flat}\wedge Y_{12}^{\flat}-X_{13}^{\flat}\wedge Y_{13}^{\flat}-X_{23}^{\flat}\wedge Y_{23}^{\flat}\,,
\]
or in terms of the abbreviation $e_{ij}:=X_{ij}^{\flat}\wedge Y_{ij}^{\flat}$  as $F=-(e_{12}+e_{13}+e_{23})$. 
For the 2-forms $\dd\eta_1$, $\dd\eta_2$ and $\dd\eta_3$ we compute
\[
\begin{tabular}{ l | l | l}
\thickline
$\dd\eta_1(X_{12}, Y_{12})=-2$ & $\dd\eta_2(X_{12}, Y_{12})=2$ & $\dd\eta_3(X_{12}, Y_{12})=0$ \\
$\dd\eta_1(X_{13}, Y_{13})=-2$ &  $\dd\eta_2(X_{13}, Y_{13})=0$ & $\dd\eta_3(X_{13}, Y_{13})=2$ \\
$\dd\eta_1(X_{23}, Y_{23})=0$ &  $\dd\eta_2(X_{23}, Y_{23})=-2$  & $\dd\eta_3(X_{23}, Y_{23})=2$ \\
\thickline
\end{tabular}
\]
and all the other components vanish (since  all the other Lie brackets lie on the horizontal distribution $\mc{D}$ and each $\eta_i$ is left-invariant).
Thus we deduce that
  \begin{align*}
\dd\eta_1&=-2(X_{12}^{\flat}\wedge Y_{12}^{\flat}+X_{13}^{\flat}\wedge Y_{13}^{\flat})=-2(e_{12}+e_{13})\,,\\
\dd\eta_2&=2(X_{12}^{\flat}\wedge Y_{12}^{\flat}-X_{23}^{\flat}\wedge Y_{23}^{\flat})=2(e_{12}-e_{23})\,,\\
\dd\eta_3&=2(X_{13}^{\flat}\wedge Y_{13}^{\flat}+X_{23}^{\flat}\wedge Y_{23}^{\flat})=2(e_{13}+e_{23})\,,
\end{align*}
and it follows that $\sum_{i=1}^{3}\dd\eta_i=0$.  The condition $2F=\dd\eta_i$ for all $i=1, 2, 3$  is thus not satisfied, so $(\phi, \eta_i, \xi_i, g)$
  does  not provide a contact metric $f$-structure in 9 dimensions.
  Moreover,  $\dd F\neq 0$, i.e., $F$ is not closed. For instance $\dd F(X_{12}, X_{13}, Y_{23})=-1\neq 0$. 
    However,  we see that
 \bp
 The metric $f$-manifold $(\U(3), \phi, \xi_1, \xi_2, \xi_3, \eta_1, \eta_2, \eta_3, g)$    is normal. 
 \ep
 \pr For $X=X_{ij}$ and $Y=Y_{ij}$ with $(ij)\in\{(12), (13), (23)\}$ 
we   obtain the relations
\[
N_{\phi}(X_{ij}, Y_{ij})=2(\xi_i-\xi_j)\in\mc{D}^{\perp}\,,\quad \sum_{k=1}^{3}\dd\eta_{k}(X_{ij}, Y_{ij})\xi_k=2(\xi_j-\xi_i)\in\mc{D}^{\perp}\,,
\]
and thus it follows thats $N^{(1)}(X_{ij}, Y_{ij})=2(\xi_i-\xi_j)+2(\xi_j-\xi_i)=0$. Next, since $\phi(\xi_k)=0=\eta_{k}(\mc{D})$ for all $k=1, 2, 3$ we see that
 \[
 N^{(1)}(X_{ij}, \xi_k)=N_{\phi}(X_{ij}, \xi_k)=\phi^2[X_{ij}, \xi_k]-\phi[\phi X_{ij}, \xi_k]
 \]
  and similarly for $N^{(1)}(Y_{ij}, \xi_k)$.  Based now on the definition of $\phi$ and on the table with the Lie brackets,  it easy to see that    
  \[
  N^{(1)}(X_{ij}, \xi_k)=0=N^{(1)}(Y_{ij}, \xi_k)\,, \ \text{for all} \  (ij)\in\{(12), (13), (23)\} \ \text{and}  \ k=1, 2, 3\,.
  \]
   Finally a direct computation shows that $N^{(1)}$ vanishes for all the pairs of horizontal vectors of the form $(X_{ij}, Y_{rs})$, $(X_{ij}, X_{rs})$, $(Y_{ij}, X_{rs})$ and $(Y_{ij}, Y_{rs})$ with $(ij)\neq (rs)$ and $(ij), (rs)\in\{(12), (13), (23)\}$.
 Let us describe a few cases. Take $X=X_{12}, Y=Y_{13}$ with $[X_{12}, Y_{13}]=-Y_{23}\in\mc{D}$.  Then
 \[
\sum_{k=1}^{3}\dd\eta_{k}(X_{12}, Y_{13})\xi_k=-\sum_{k=1}^{3}\eta_{k}([X_{12}, Y_{13}])\xi_k=0\,.
\]
  Thus
 \begin{eqnarray*}
N^{(1)}(X_{12}, Y_{13})&=&N_{\phi}(X_{12}, Y_{13})=[\phi X_{12}, \phi Y_{13}]+\phi^2[X_{12}, Y_{13}]-\phi[X_{12}, \phi Y_{13}]-\phi[\phi X_{12}, Y_{13}]\\
&=&-[Y_{12}, X_{13}]+\phi^2(-Y_{23})+\phi[X_{12}, X_{13}]-\phi[Y_{12}, Y_{13}]\\
&=&-Y_{23}+Y_{23}-\phi(X_{23})+\phi(X_{23})=0
\end{eqnarray*}
since $\phi^2(-Y_{23})=Y_{23}$ as $Y_{23}\in\mc{D}$.  \\
\noindent Take now $X=X_{12}, Y=X_{23}$. Since  $[X_{13}, X_{23}]=-X_{12}\in\mc{D}$, in a similar way as above we get 
 \begin{eqnarray*}
N^{(1)}(X_{12}, X_{23})&=&N_{\phi}(X_{12}, X_{23})=[\phi X_{12}, \phi X_{23}]+\phi^2[X_{12}, X_{23}]-\phi[X_{12}, \phi X_{23}]-\phi[\phi X_{12}, X_{23}]\\
&=&[Y_{12}, Y_{23}]+\phi^2(X_{13})-\phi[X_{12}, Y_{23}]-\phi[Y_{12}, X_{23}]\\
&=&-X_{13}-X_{13}-\phi(Y_{13})-\phi(Y_{13})=-2X_{13}+2X_{13}=0\,.
\end{eqnarray*}
Next take $X=Y_{12}$ and $Y=X_{23}$ with $[Y_{12}, X_{23}]=Y_{13}\in\mc{D}$. Then
 \begin{eqnarray*}
N^{(1)}(Y_{12}, X_{23})&=&N_{\phi}(Y_{12}, X_{23})=[\phi Y_{12}, \phi X_{23}]+\phi^2[Y_{12}, X_{23}]-\phi[Y_{12}, \phi X_{23}]-\phi[\phi Y_{12}, X_{23}]\\
&=&-[X_{12}, Y_{23}]+\phi^2(Y_{13})-\phi[Y_{12}, Y_{23}]+\phi[X_{12}, X_{23}]\\
&=&-Y_{13}-Y_{13}+\phi(X_{13})+\phi(X_{13})=-2Y_{13}+2Y_{13}=0\,.
\end{eqnarray*}
Finally take $X=Y_{13}$ and $Y=Y_{23}$ with $[Y_{13}, Y_{23}]=-X_{12}\in\mc{D}$. The we also see that
\begin{eqnarray*}
N^{(1)}(Y_{13}, Y_{23})&=&N_{\phi}(Y_{13}, Y_{23})=[\phi Y_{13}, \phi Y_{23}]+\phi^2[Y_{13}, Y_{23}]-\phi[Y_{13}, \phi Y_{23}]-\phi[\phi Y_{13}, Y_{23}]\\
&=&[X_{13}, X_{23}]+\phi^2(-X_{12})+\phi[Y_{13}, X_{23}]+\phi[X_{13}, Y_{23}]\\
&=&-X_{12}+X_{12}-\phi(Y_{12})+\phi(Y_{12})=0\,.
\end{eqnarray*}
Similarly are treated the rest pairs.
\pro

Since $N^{(1)}$ vanishes and each $\xi_i$ is Killing, by Proposition \ref{normal_skewtorsion}
we deduce that   $(\U(3), \phi, \xi_i, \eta_i, g)$ $(i=1, 2, 3)$   admits a unique almost metric $\phi$-connection $\nabla$ with skew-torsion given by  
\[
\nabla=\nabla^g+\frac{1}{2}\Big\{\sum_{i=1}^{3}\eta_{i}\wedge\dd\eta_i+\dd^{\phi} F\Big\}\,.
\]
In order to describe $T$ we essentially  need to compute $\dd^{\phi} F$. To do so recall that since $F$ is left-invariant, one has
\begin{eqnarray*}
\dd F(X, Y, Z)&=&-F([X, Y], Z)-F([Z, X], Y)-F([Y, Z], X)\\
&=&-g([X, Y], \phi Z)-g([Y, Z], \phi X)-g([Z, X], \phi Y)\,,
\end{eqnarray*}
where $X, Y, Z$ are arbitrary left-invariant vector fields on  $\U(3)$. 
On the other hand, since $\phi(\xi_i)=0$ for all $i=1, 2, 3$, $\dd F(\phi X, \phi Y, \phi Z)$ may have only horizontal non-trivial components. By definition,  for $X, Y, Z\in\Gamma(\mc{D})$ we have $\phi^2(X)=-X$, $\phi^2(Y)=-Y$ and $\phi^2(Z)=-Z$ and a replacement of $X, Y, Z$ by $\phi X, \phi Y, \phi Z$  in the previous formula gives
\[
\dd F (\phi X, \phi Y, \phi Z)=g([\phi X, \phi Y], Z)+g([\phi Y, \phi Z], X)+g([\phi Z, \phi X], Y)\,.
\]
To apply this formula we only need the Lie brackets   and the definition of $\phi$. Based on this approach  we see that the only non-zero components are given by  $\dd F(\phi X_{12}, \phi X_{13}, \phi X_{23})=-1$, $\dd F(\phi X_{12}, \phi Y_{13}, \phi Y_{23})=-1$,  $\dd F(\phi X_{13}, \phi Y_{12}, \phi Y_{23})=-1$ and $\dd F(\phi X_{23}, \phi Y_{12}, \phi Y_{13})=-1$.  Therefore,  
we deduce that
\[
\dd F(\phi X, \phi Y, \phi Z)
= -\big(
X_{12}^{\flat}\wedge X_{13}^{\flat}\wedge X_{23}^{\flat}
+ X_{12}^{\flat}\wedge Y_{13}^{\flat}\wedge Y_{23}^{\flat}
+ X_{13}^{\flat}\wedge Y_{12}^{\flat}\wedge Y_{23}^{\flat}
+ X_{23}^{\flat}\wedge Y_{12}^{\flat}\wedge Y_{13}^{\flat}
\big)(X,Y,Z)
\]
for any $X, Y, Z\in\Gamma(\mc{D})$ and in particular for any three left-invariant vector fields $X, Y, Z$ of $\U(3)$.
It follows that
\begin{eqnarray*}
T&=&\sum_{i=1}^3 \eta_i \wedge\dd\eta_i  + \dd^{\phi}F\\
&=& -2\,\eta_1 \wedge \big(X_{12}^{\flat}\wedge Y_{12}^{\flat} + X_{13}^{\flat}\wedge Y_{13}^{\flat}\big)\\
&&+ 2\,\eta_2 \wedge \big(X_{12}^{\flat}\wedge Y_{12}^{\flat} - X_{23}^{\flat}\wedge Y_{23}^{\flat}\big)\\
&&+ 2\,\eta_3 \wedge \big(X_{13}^{\flat}\wedge Y_{13}^{\flat} + X_{23}^{\flat}\wedge Y_{23}^{\flat}\big)\\
&&+\big(X_{12}^{\flat}\wedge X_{13}^{\flat}\wedge X_{23}^{\flat}
+ X_{12}^{\flat}\wedge Y_{13}^{\flat}\wedge Y_{23}^{\flat}
+ X_{13}^{\flat}\wedge Y_{12}^{\flat}\wedge Y_{23}^{\flat}
+ X_{23}^{\flat}\wedge Y_{12}^{\flat}\wedge Y_{13}^{\flat}
\big)\,.
\end{eqnarray*}
Using this expression one can  easily  verify the relations  $T(\xi_i, X, Y)=\dd\eta_{i}(X, Y)$ for all $i=1, 2, 3$, and $X, Y\in\Gamma(\mc{D})$, and $T(\xi_i, \xi_j, X)=0$ for all $i, j\in\{1, 2, 3\}$ and $X\in\Gamma(\mc{D})$ (which both can be  extended to the whole of $\U(3)$).

%%%%%%%%%%%%%%%%%%%%%%%%%
%%%%%%%%%bibliography%%%%%%%%%%
%%%%%%%%%%%%%%%%%%%%%%%%%

\medskip
{ }

\begin{thebibliography}{50}

 \bibitem[A06]{Srni}   I.~Agricola,  {\it The Srn\'i lectures on non-integrable 
geometries with torsion}, Arch. Math., 42 (2006), 5--84. With an appendix by M.~Kassuba. 

\bibitem[AD20]{AD20}
 I.~Agricola and G.~Dileo,
 {\it Generalizations of 3-Sasakian manifolds and skew torsion}, 
 Adv. Geom. 20 (2020), 331--374.
 
 \bibitem[AF04]{AF04}
 I.~Agricola and T.~Friedrich, 
{\it On the holonomy of connections with skew-symmetric torsion},
 Math. Ann., 328, (2004) 711--748.
 
\bibitem[AFF15]{AFF15}
   I.~Agricola, A.C.~Ferreira and T.~Friedrich, 
 {\it The classification of naturally reductive homogeneous spaces in dimensions  $\leq 6$},
  Differ. Geom. Appl., 39, (2015), 59--92.

\bibitem[AS53]{AS53}
 A.~Ambrose and I.~Singer, 
 {\it A theorem on holonomy}, 
 Trans. Amer. Math. Soc., 75 (1953),  428--443.

 \bibitem[B03]{Blair} 
  D.~E.~Blair,
 {\it Riemannian geometry of contact and symplectic manifolds}, (Second Edition), 
 Birkh\"auser, Progress in Math., Vol. 203, (2002).

 \bibitem[B70]{Blair70}
 D.~E.~Blair, 
 {\it Geometry of manifolds with structural group $U(n)\times O(s)$}, 
 J. Diff. Geom.,  4 (1970), 155--167.

\bibitem[BL69]{BL69}
D.~E.~Blair and  G.~Ludden,
{\it Hypersurfaces in almost contact manifolds},
T\^ohoku Math. J. 2, (1969), 354--362
  
 \bibitem[BLY73]{BLY73}
 D.~E.~Blair, G.~Ludden and K. Yano,
 {\it Differential geometric structures on principal toroidal bundles.}
  Trans. Amer. Math. Soc., 181 (1973), 175--184.
  
\bibitem[BP08]{BP08}
  L.~Brunetti and A. M.~Pastore,
 {\it Curvature of a class of indefinite globally framed $f$-manifolds},
  Bull. Math. Soc. Sci. Math. Roumanie (N.S.) 51(99) (2008), 183--204.
  
  \bibitem[CFF90]{CFF90}
 J.~Cabrerizo, L.~M.~Fern\'andez and M.~Fern\'andez,
 {\it The curvature tensor fields on $f$-manifolds with complemented
frames},
An. \c{S}tiin\c{t}. Univ. Al. I. Cuza Ia\c{s}i Sec\c{t}. I a Mat. 36,  (1990), 151--161.

 \bibitem[CS04]{CS04}
 R.~Cleyton and A.~Swann, 
 {\it Einstein metrics via intrinsic or parallel torsion}, 
 Math. Z. 247 (2004), 513--528.
 
  \bibitem[CMS21]{CMS21}
 R.~Cleyton, A.~Moroianu and U.~Semmelmann,
 {\it Metric connections with parallel skew-symmetric torsion},
 Adv. Math., 378 (2021), 107519
 
 \bibitem[DT04]{T04}
 L.~Di~Terlizzi, 
 {\it On a generalization of contact metric manifolds}, 
 Publ. Math. Debrecen, 64, 3-4, (2004), 401--413.
 
 \bibitem[DTK07]{TK07}
  L.~Di~Terlizzi and J.~Konderak, 
  {\it Examples of a generalization of contact metric structures on fibre bundles},
 J. Geom. 87, (2007), 31--49.
 
% \bibitem[DTP13]{TP13}
%  L.~Di~Terlizzi and A.~M.~Pastore, 
%  {\it $\mc{K}$-manifolds locally described by Sasaki manifolds}, 
%  An. \c{S}t. Univ. Ovidius Constan\c{t}a,   21(3), (2013), 269--287.
 
 \bibitem[DL05]{DL05}
 G.~Dileo and A.~Lotta, 
 {\it On the structure and symmetry properties of almost $\mc{S}$-manifolds}, 
  Geom. Dedicata (2005) 110, 191--211.
  
 \bibitem[DIP01]{DIP01}
K.~L.~Duggal, S.~Ianus and A.~M.~Pastore, 
{\it Maps interchanging $f$-structures and their harmonicity}, 
 Acta Appl. Math., 67, (2001), 91--115.
 
\bibitem[FI02]{FrIv}
 T.~Friedrich and S.~ Ivanov, 
 {\it Parallel spinors and connections with skew-symmetric torsion in string theory}, 
 Asian J. Math., 6  (2),  (2002), 303--335.  
 
 \bibitem[GYS70]{GY70}
 S.~I.~Goldberg and K.~Yano,
{\it On normal  globally framed $f$-manifolds,}
T\^ohoku Math. J. (2), 22, (1970), 362--370.

\bibitem[G72]{G72}
 S.~I.~Goldberg, 
{\it A generalization of K\"ahler geometry}, 
J.~Diff.~Geom.,   6 (1972), 343--355.
 
 \bibitem[HOA86]{HOA86}
 I.~Hasegawa, Y.~Okuyama and T.~Abe,
 {\it On p-th Sasakian manifolds},
  J. Hokkaido Univ. of Education, Section II A, 37(1), (1986), 1--16.

  \bibitem[IY64]{IsY64}
  S.~Ishihara and K.~Yano, 
{\it On integrability of a structure $f$ satisfying $f^3+f=0$}, 
Quart. J. Math. Oxford (2), 15 (1964), 217--222.

\bibitem[IP01]{IvP01}
S.~Ivanov and G.~Papadopoulos, 
{\it Vanishing theorems and string backgrounds}, 
Class. Quant. Grav., 18 (2001),   1089--1110.

 \bibitem[KT72]{KT72}
M.~Kobayashi and S.~Tsuchiya, 
{\it Invariant submanifolds of an $f$-manifold with complemented frames}, 
Kodai Math. Sem. Rep. 24, (1972), 430--450.

\bibitem[KV83]{KV83}
O.~Kowalski and L.~Vanhecke, 
{\it Four-dimensional naturally reductive homogeneous spaces}, 
Differential geometry on homogeneous spaces, Conf. Torino/Italy 1983, Rend. Semin. Mat., Torino, Fasc. Spec., (1983), 223--232. 

\bibitem[LP04]{LP04}
A.~Lotta and A.~M.~Pastore, 
{\it The Tanaka-Webster connection for almost $\mc{S}$-manifolds and Cartan geometry},
Archivum Math. (Brno), 40, (2004), 47--61.

\bibitem[M93]{Miz93}
R.~I.~Mizner,
{\it  Almost CR structures, $f$-structures, almost product structures and associated connections}, 
Rocky Mount. J. Math., 23 (1993), 1337--1359.

\bibitem[R85]{R85}
J. H. Rawnsley, 
{\it $f$-structures, $f$-twistor spaces and harmonic maps}, in {\it Geometry Seminar ``Luigi
Bianchi'' II--1984}, 85--159, Lecture Notes in Math., 1164, Springer, Berlin, 1985.

\bibitem[S18]{S18}
R.~Storm, 
{\it A new construction of naturally reductive spaces},
 Transform. Groups 23(2) (2018), 527--553.

\bibitem[T89]{T89}
S.~Tanno,
{\it Variational problems on contact Riemannian manifolds},
 Trans. Amer. Math. Soc., 314 (1989), 349--379.
 
\bibitem[Y63]{Yano63}
K.~Yano, 
{\it On a structure defined by a tensor field of type $(1, 1)$  satisfying $f^3+f=0$},
Tensor (N.S.) 14 (1963), 99--109.

\end{thebibliography}
\end{document}